\documentclass[11pt, a4paper]{article}
\usepackage{amsmath}
\usepackage{amsfonts}
\usepackage{amssymb}
\usepackage{ulem}

\usepackage{amsthm}
\usepackage{url}
\usepackage{dcolumn}

\usepackage{cite}
\usepackage{fancyhdr}

\usepackage{graphicx, color}
\usepackage{hyperref}
\usepackage{xcolor}

\newtheorem*{rem*}{Remark}

\newtheorem*{prob}{Problem}

\newtheorem{nameth}{Restriction on The Number of Types of Congruent Diagrams}
\newtheorem*{nameth2}{Tile Combination Specification Condition}

\usepackage{here}

\begin{document}

\title{Aperiodic sets of three types of convex polygons}
\author{ Teruhisa SUGIMOTO$^{ 1), 2)}$ }
\date{}
\maketitle

{\footnotesize

\begin{center}
$^{1)}$ The Interdisciplinary Institute of Science, Technology and Art

$^{2)}$ Japan Tessellation Design Association

E-mail: ismsugi@gmail.com
\end{center}

}

{\small
\begin{abstract}
\noindent
Sets of three types of convex pentagons that are 
aperiodic with no matching conditions on the edges are created from a chiral 
aperiodic monotile \mbox{Tile$(1, 1)$}. This method divides the interior of \mbox{Tile$(1, 1)$} 
into five convex polygons with five or more edges, and we have so far 
identified four methods.
\end{abstract}
}

\textbf{Keywords:} Tiling, aperiodic, non-periodic, polygon, tile

%%%%%%%%%%%%%%%%%%%%%%%%%%%%%%%%%%%%%%%%%%%%%%%%%%%%%%%%%%%%%%%%%%%%%%
%%%%%%%%%%%%%%%%%%%%%%%%%%%%%%%%%%%%%%%%%%%%%%%%%%%%%%%%%%%%%%%%%%%%%%
\section{Introduction}
\label{section1}

A \textit{tiling} (or \textit{tessellation}) of the plane is a collection of sets, called tiles, 
that covers the plane without gaps or overlaps, except for the boundaries of the tiles. 
The term ``tile" refers to a topological disk, whose boundary is a simple closed curve.

A tiling exhibits \textit{periodicity} if its translation by a non-zero vector coincides with 
itself; a tiling is considered periodic if it coincides with its translation by two 
linearly independent vectors. However, in this study, a tiling with periodicity is 
referred to as \textit{periodic}, and a tiling without periodicity is referred to as 
\textit{non-periodic}.

In general, when considering tilings, there are only a finite number of tile shapes 
that can be used. These finite-shape diagrams are called \textit{prototiles}, 
and a set of prototiles admits tilings of the plane. The tiling generated by a set of 
prototiles can cover the plane infinitely. A set of prototiles is said to be \textit{aperiodic} 
if copies of the prototiles can be assembled into tilings of the plane such that 
all tilings with the prototiles are non-periodic \cite{ Hallard_1991, G_and_S_1987, 
wiki_aperiodic_set}.

The most well-known problem regarding a set of prototiles (set of tiles, , tile set) that 
are aperiodic, that is, ``aperiodic set of prototiles (aperiodic set of tiles)" (hereafter 
abbreviated as ``ASP"), is ``Is there a single aperiodic prototile (with or without a 
matching condition\footnote{
Matching conditions specify how tiles must connect to form a valid tiling, which can sometimes 
be represented by assigning colors or orientations to specific edges of the prototiles.
}), that is, one that admits only non-periodic tilings by congruent 
copies?" \cite{Hallard_1991}. Smith et al. presented solutions to this problem for 
concave diagram tiles without matching conditions \cite{Smith_2024a, Smith_2024b}. 
In \cite{Hallard_1991}, another major problem concerning ASP was presented as 
follows:

\begin{prob}[] 
\label{pro1}
There is a set of three convex polygons that are aperiodic with no matching condition 
on the edges, but are two or even one such convex prototiles possible?
\end{prob}

In this study, our discussion focuses on the above problem. Therefore, we 
consider a set of tiles (prototiles) that are polygons and have no matching 
conditions on the edges.

\renewcommand{\figurename}{{\small Figure}}
\begin{figure}[htbp]
\centering\includegraphics[width=14.5cm,clip]{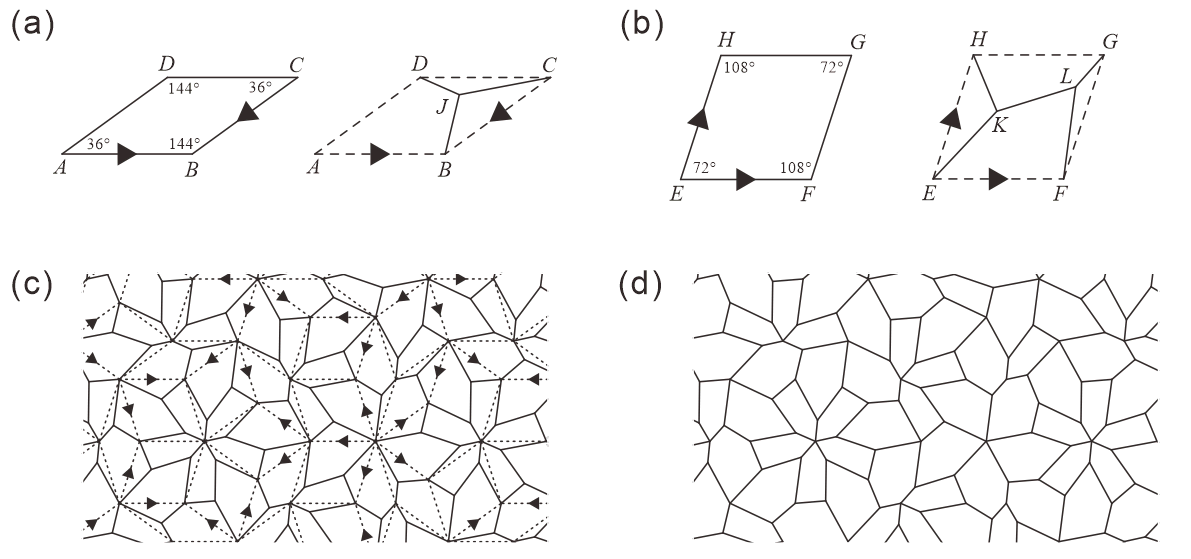} 
\caption{{\small 
Tiling (d) generated by an aperiodic set of prototiles created by Ammann, which used 
three convex polygons with no matching condition on the edges. As shown in (c), 
this is a recomposition by Penrose tiles (rhombuses) and is obtained by 
markings of using the points $J$, $K$, and $L$ as shown in (a) and (b). These 
markings are completely determined by the choice of the point $J$, as we must 
have \textit{GL} = \textit{DJ}, \textit{FL} = \textit{CJ} =\textit{EK}, 
and \textit{HK} = \textit{BJ}. For aperiodicity, the point $J$ must be chosen so 
that \textit{DJ}, \textit{CJ}, \textit{BJ}, and \textit{KL }are of different 
lengths \cite{G_and_S_1987, Sugimoto_2017}.
}
\label{fig01}
}
\end{figure}

The ``ASP comprising three convex polygons (= the set of three convex 
polygons that are aperiodic with no matching condition on the edges)" in the 
above problem is known for the two convex pentagons and one convex hexagon 
created by Ammann, as shown in Figure~\ref{fig01}, based on the Penrose 
tiles\footnote{
The Penrose tiles are considered to be an ASP comprising two convex 
polygons (rhombuses) with matching condition on the edges 
\cite{Hallard_1991, G_and_S_1987, Sugimoto_2017}.} 
\cite{G_and_S_1987, Sugimoto_2017}. 
In the tiling of the two convex pentagons and one 
convex hexagon shown in Figure~\ref{fig01}, the reflected (flipped) convex pentagons 
and hexagon were not used. The tiling is edge-to-edge\footnote{ 
A tiling by convex polygons is \textit{edge-to-edge} if any 
two convex polygons in a tiling are either disjoint or share one vertex or an entire 
edge in common. Then other case is \textit{non-edge-to-edge} \cite{G_and_S_1987, 
Sugimoto_2016, Sugimoto_2017}.}.

In tiling problems, it is important to note how the reflected tiles 
(posterior side tiles) are handled. For example, if all the tiles in a 
tiling are of the same size and shape, then the tiling is \textit{monohedral}, which 
allows the use of reflected tiles in monohedral tilings \cite{G_and_S_1987}. 
In other words, in monohedral tilings, we treat the anterior and posterior side 
tiles as the same type (the concept is that there is only one type of tile). In this 
study, a tile that admits monohedral tiling is referred to as a \textit{monotile}. 
Note that there are monotiles that require the use of reflected 
tiles during the tiling generation process (see \ref{appA}). 
In some tiling problems, congruent diagrams, including reflection images, are treated 
as tiles of the same type as in monohedral tiling.

According to Smith et al. \cite{Smith_2024a}, \mbox{Tile$(a, b)$}, which can 
generate only non-periodic tiling, is an ``aperiodic monotile," in which all tiles are of 
the same size and shape; however, the tiling generation process requires the use of 
reflected tiles, and both anterior and posterior side tiles must be treated as 
the same type. However, according to Smith et al. \cite{Smith_2024b}, the tiles 
``Spectres," which can generate only non-periodic tiling, do not require the use 
of reflected tiles during the tiling generation process\footnote{
Spectres can generate only non-periodic tiling using only one side, either anterior 
or posterior, even if it allows the use of reflected tiles during the tiling generation 
process.}. The tiles are referred to as ``chiral aperiodic monotiles."

Let us consider ``ASP comprising $n$ types of polygons ($n$ is an integer 
greater than or equal to one)" by focusing on ``a set of three convex 
polygons that are aperiodic" in the above-mentioned problem. 
For the polygons (tiles) contained in the ASP and the set of original problems, 
if the polygons are congruent, they are treated as the same type; however, 
from the text (or symbols) alone, it is not clear whether the anterior or 
posterior side polygons are treated as the same type. Further, similar to 
the problem of monohedral tiling, it is possible that polygons in those cases 
may sometimes be treated as the same type of tiles in the reflected image.

Therefore, in this study, we proceed by treating the anterior and posterior 
side tiles as different types, even if they are congruent diagrams. 
In other words, we assume that the reflected tiles cannot be 
used (the tiles in the tile set cannot be reflected) during all tiling generation 
process with the tiles in a tile set. However, when we use the terms ``viewpoint 
that does not distinguish between anterior and posterior sides" or ``monotile," 
we treat the anterior and posterior side tiles as the same type.

Following this policy, the case in \cite{Smith_2024a} is ``ASP comprising two types of 
concave polygons," and the case in \cite{Smith_2024b} is ``ASP comprising one 
type of concave diagram\footnote{
\mbox{Tile$(1, 1)$} corresponding to $a=b=1$ of \mbox{Tile$(a, b)$} shown 
in \cite{Smith_2024a} is a concave polygon and a concave diagram, and Spectres 
in \cite{Smith_2024b} are concave diagrams. From a viewpoint that does not 
distinguish between anterior and posterior sides, both cases of \cite{Smith_2024a} 
and \cite{Smith_2024b} are considered ``ASP comprising one type of concave diagram." 
Note that \mbox{Tile$(1, 1)$} with 13-edges is considered an equilateral 14-edges polygon.}". 
Moreover, the case shown by Ammann in Figure~\ref{fig01} is ``ASP comprising 
three types of convex polygons (two types of convex pentagons and one type 
of convex hexagon)."

This study presents examples of an ``ASP comprising three types of convex 
polygons," similar to the case shown in Figure~\ref{fig01}.

Note that if we extract the convex polygons of $n - 1$ or fewer types 
contained in the set, we assume that they cannot admit a tiling. It is because it is 
natural that the convex polygons of $n - 1$ or fewer types extracted from the set 
cannot generate periodic tiling where $n \ge 2$; however, if they can generate 
non-periodic tiling, it is replaced by the result that there exists an ASP 
comprising fewer types of convex polygons.

The two types of convex pentagons and one type of convex 
hexagon shown in Figure~\ref{fig01} can be identified as not belonging to known Type 
families of convex pentagonal monotiles (see \ref{appA}) or convex hexagonal 
monotiles (see \ref{appB}), respectively\footnote{ 
If Rao's claim in \cite{Rao_2017} is formally established as correct, there will be no ASP 
for convex polygonal monotiles (i.e., there does not exist ASP comprising one type of convex 
polygon with no matching condition on the edges), and furthermore, it can be 
said that ``a convex pentagon does not belong to known Type families" and 
``a convex pentagon cannot admit tiling" are equivalent. At the stage where 
these are not truly established \cite{wiki_pentagon_tiling}, it is necessary to confirm 
whether tiling is truly impossible for a convex pentagon, including whether it is 
possible to generate only non-periodic tilings. In this study, we assume that 
Rao's claim is correct and proceed with the discussion. 
However, in Section~\ref{section2} (and \ref{appD}) of this study, to check whether 
or not a convex polygon is a monotile, we do not judge based on only the tile 
conditions of known Types, and are checking whether tiling is actually possible.
Note that if convex polygons are convex polygonal monotiles that can generate 
edge-to-edge tilings, it has been proved that they can always generate periodic 
tilings \cite{Sugimoto_2016, Sugimoto_2017}.
}. Because the set of these polygons is the ASP 
comprising three types of convex polygons, if two types of polygons are extracted 
from three types of convex polygons, those two types of polygons cannot 
admit tilings (but the author has not confirmed this).

%%%%%%%%%%%%%%%%%%%%%%%%%%%%%%%%%%%%%%%%%%%%%%%%%%%%%%%%%%%%%%%%%%%%%%
%%%%%%%%%%%%%%%%%%%%%%%%%%%%%%%%%%%%%%%%%%%%%%%%%%%%%%%%%%%%%%%%%%%%%%

\section{Dividing \mbox{Tile$(1, 1)$} into ASP comprising convex polygons}
\label{section2}

\begin{figure}[htbp]
\centering\includegraphics[width=15cm,clip]{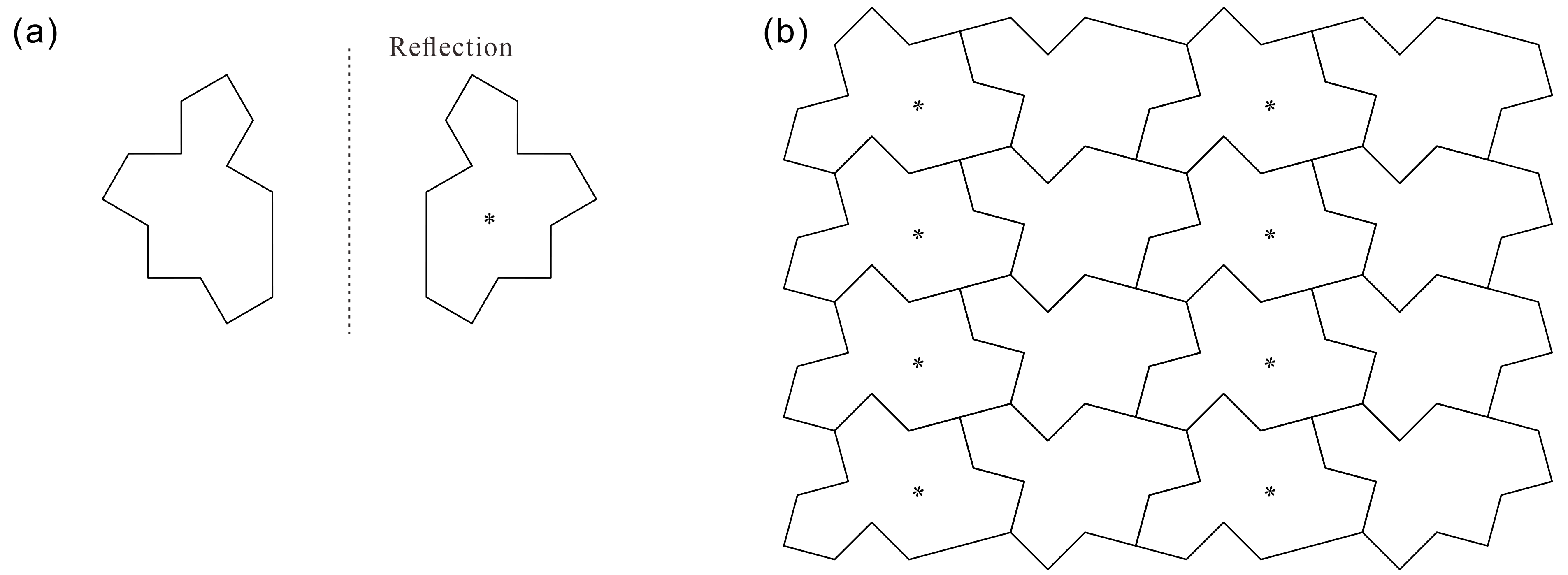} 
\caption{{\small 
\mbox{Tile$(1, 1)$}, and periodic tiling with anterior and posterior sides of \mbox{Tile$(1, 1)$}.
}
\label{fig02}
}
\end{figure}

\begin{figure}[htbp]
\centering\includegraphics[width=15cm,clip]{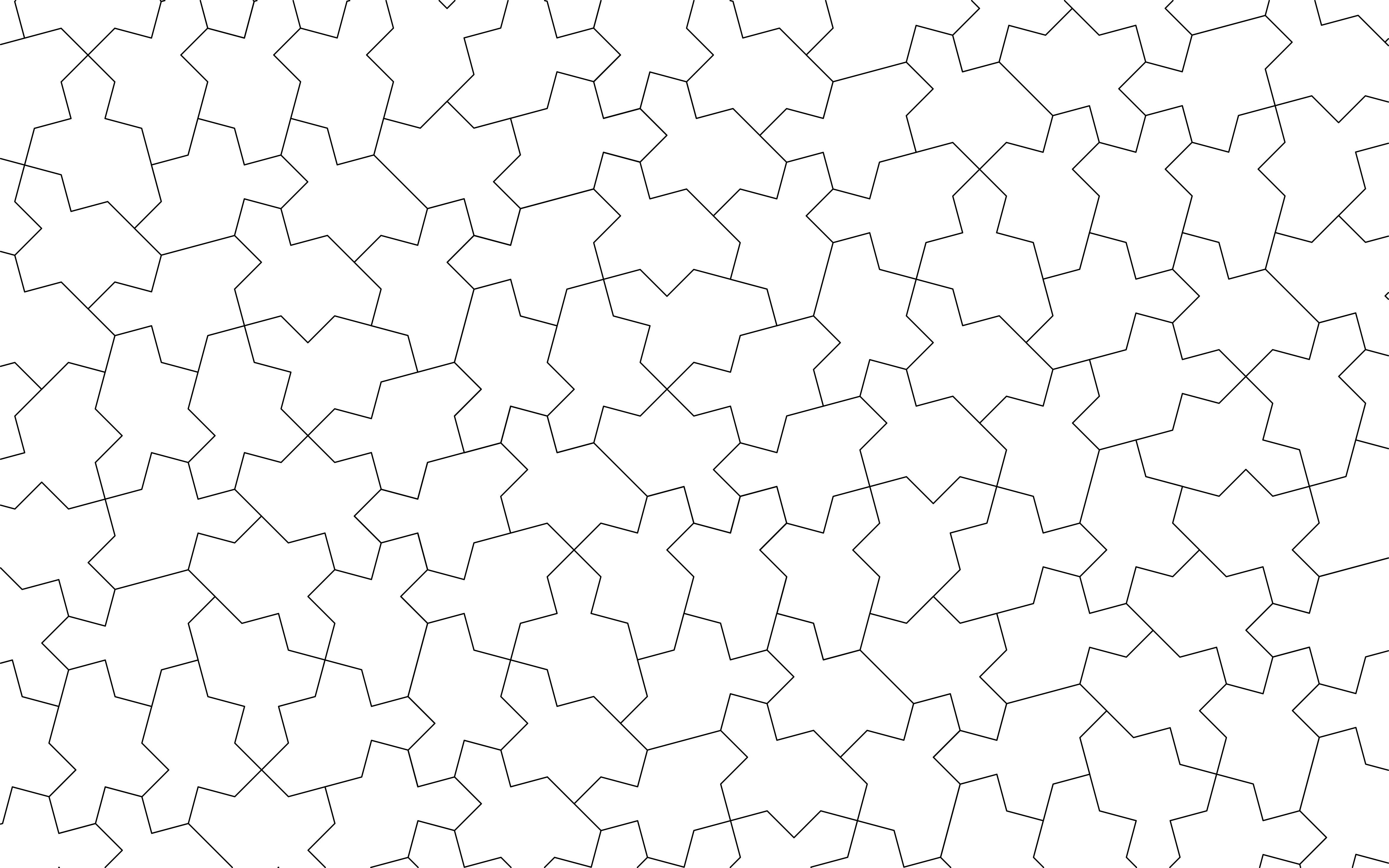} 
\caption{{\small 
Non-periodic tiling generated by using only one side of \mbox{Tile$(1, 1)$}.
}
\label{fig03}
}
\end{figure}

\mbox{Tile$(1, 1)$}, corresponding to $a=b=1$ in \mbox{Tile$(a, b)$} shown in 
\cite{Smith_2024a}, can generate a periodic tiling if the use of reflected tiles is 
allowed during the tiling generation process (see Figure~\ref{fig02}). However, 
as shown in \cite{Smith_2024b}, \mbox{Tile$(1, 1)$} can only generate non-periodic 
tiling if and only if it does not allow the use of reflected tiles (see Figure~\ref{fig03}). 
Using this property, ASP comprising three types of convex polygons can be 
obtained from \mbox{Tile$(1, 1)$}. This method divides the interior of 
\mbox{Tile$(1, 1)$} into five convex polygons with five or more edges\footnote{
Because all triangles and quadrilaterals are monotiles (which do not require the 
use of reflected tiles to form periodic tilings)
the number of edges of convex polygons contained in ASP is at least five. 
Then, there are no convex polygonal monotiles with seven or more edges 
\cite{Gardner_1975, G_and_S_1987, Klam_1980, Niven_1978, Reinhardt_1918, Zong_2020}.
}, and we have so far identified four methods 
(Methods 1--4). The tile sets created by the division of the four methods all 
contain the convex pentagon $P_{1}$, which cannot generate a tiling (that is, no 
monotile) in the upper-left corner of Figure~\ref{fig04} ($P_{1}$ is convex 
pentagons with line symmetry and cannot distinguish between the anterior 
and posterior sides). Methods 1 (Figure~\ref{fig04}), 2 (Figure~\ref{fig05}), 
and 4 (Figure~\ref{fig18v2}) divide \mbox{Tile$(1, 1)$} to create five convex 
polygons containing three $P_{1}$. In Method 3 (Figure~\ref{fig15}), 
one $P_{1}$ (= $CP_{1}(X_{3})$) is contained in the five convex polygons formed by 
dividing \mbox{Tile$(1, 1)$}. We note that the tilings generated by the ASP, 
comprising three types of convex polygons created by Methods 1--4, are 
non-edge-to-edge.

Here again, it is emphasized that each convex polygon in the tile set 
created by dividing \mbox{Tile$(1, 1)$} satisfies the condition that the reflected 
tiles cannot be used during all tiling generation processes. As described 
above, \mbox{Tile$(1, 1)$} can generate periodic tilings if the reflected tiles can be 
used during the tiling generation process; therefore, it is obvious that the 
three convex polygons in the tile set can generate periodic tilings if they can 
be reflected during the tiling generation process. ``All tiling generation 
processes" means that even if we try to generate tiling using only one or 
two types of convex polygons in a tile set, we cannot use the reflected 
convex polygons (tiles) during the tiling generation process.

As described in Section~\ref{section1}, some monotiles require the use of reflected 
tiles during tiling generation process. It is important to note that this property 
means that some tiling patterns cannot be formed without using reflected 
tiles, but the tiles forming the patterns do not necessarily have to use 
reflected tiles to generate tiling (see \ref{appA} for details) \cite{G_and_S_1987,
Sugimoto_2016, Sugimoto_2017}.

Note that in this study, we assumed that the reflected tiles cannot be used 
during the tiling generation process; therefore, even if the convex polygon 
of a tile set created by dividing \mbox{Tile$(1, 1)$} is a monotile, we do not 
consider that the tile set is not aperiodic. This is because if the tile set 
contains convex polygonal monotiles that require the use of reflected tiles 
during the tiling generation processes and does not contain the reflected 
convex polygons of the monotiles; then, the possibility that the tile set is 
aperiodic cannot be ruled out.

\renewcommand{\figurename}{{\small Figure}}

\renewcommand{\figurename}{{\small Figure}}

%%%%%%%%%%%%%%%%%%%%%%%%%%%%%%%%%%%%%%%%%%%%%%%%%%%%%%%%%%%%%%%%%%%%%%
\subsection{Method 1 and Method 2}
\label{subsec2_1}

\begin{figure}[htbp]
\centering\includegraphics[width=15cm,clip]{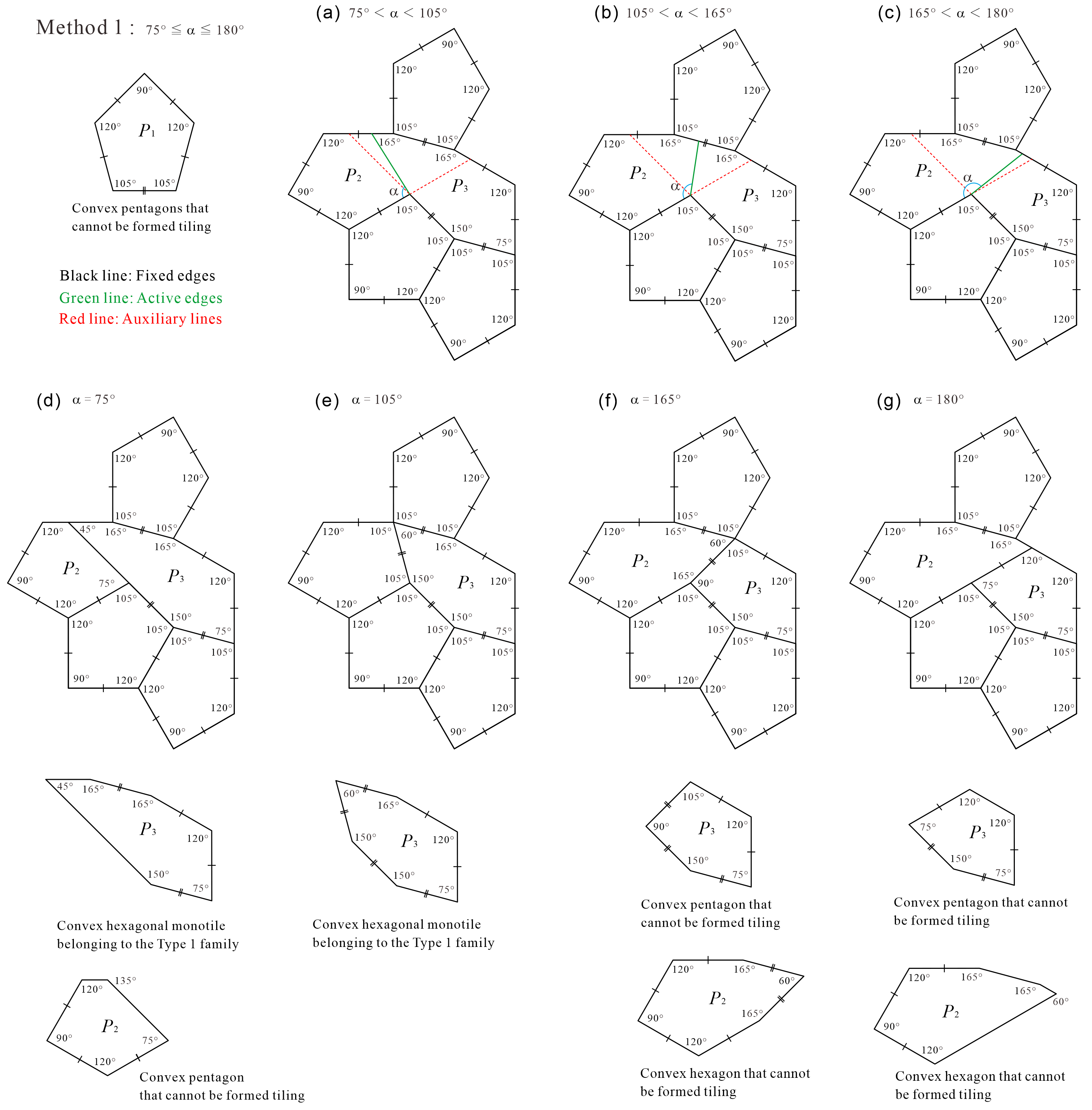} 
\caption{{\small 
Cases in Method 1.
}
\label{fig04}
}
\end{figure}

\begin{figure}[htbp]
\centering\includegraphics[width=15cm,clip]{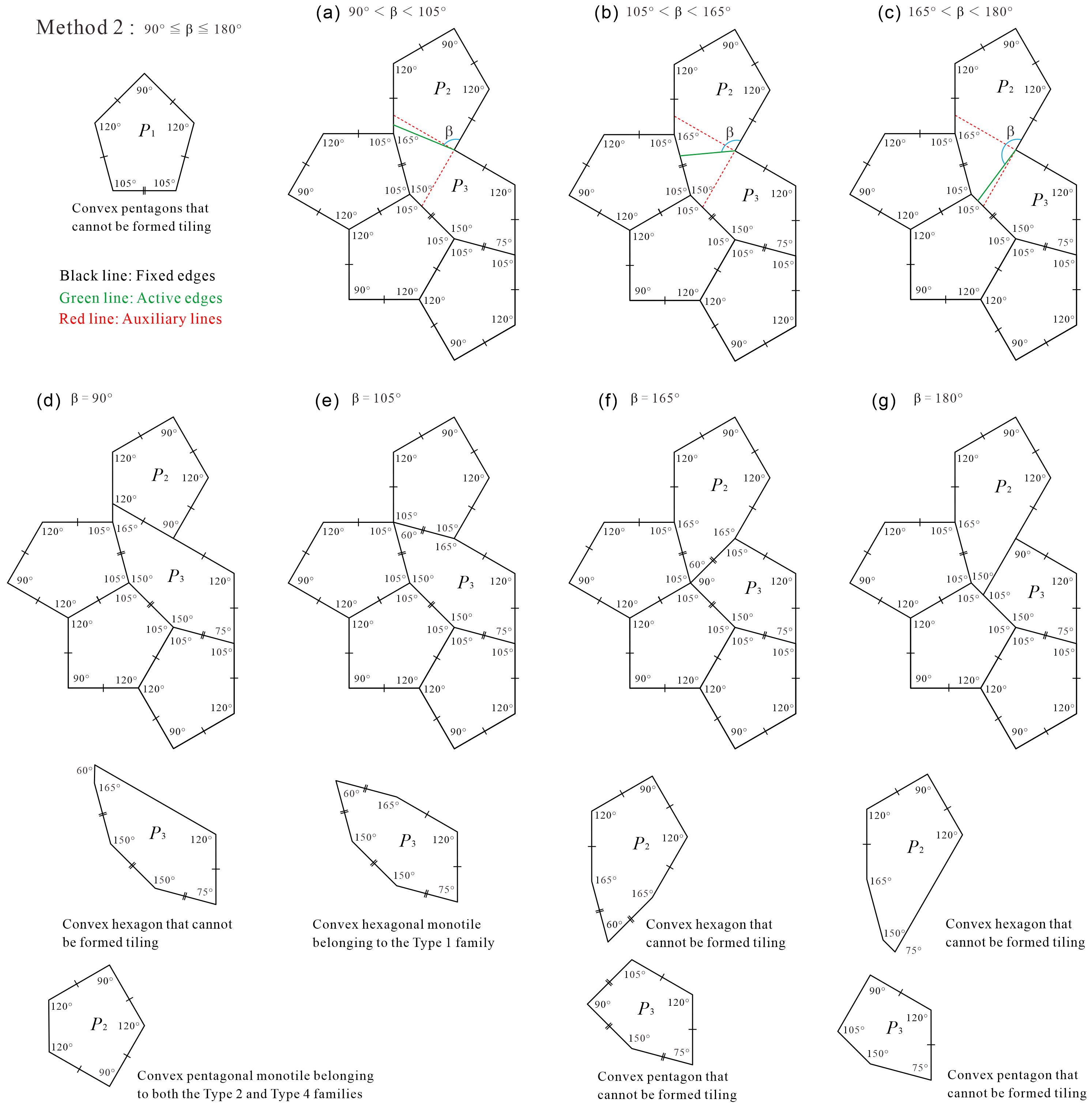} 
\caption{{\small 
Cases in Method 2.
}
\label{fig05}
}
\end{figure}

Method 1 in Figure~\ref{fig04} and Method 2 in Figure~\ref{fig05} are based on the division 
of \mbox{Tile$(1, 1)$} into five convex polygons containing three $P_{1}$. The division 
in Method 1 depends on the value of the parameter $\alpha$ shown in Figure~\ref{fig04}, 
and in Method 2, it depends on the value of the parameter $\beta$ shown 
in Figure~\ref{fig05}. In the Methods, there are clear cases in which the tile set is 
not an ASP. For example, there are the following cases.

\begin{itemize}
\item In the case of $\alpha = 75^ \circ $ in Method 1 (see Figure~\ref{fig04}(d)), 
there is the convex hexagonal monotile $P_{3}$ belonging to the Type 1 family that 
can form periodic tiling without using reflected tiles (see \ref{appB}); therefore, 
the set of three types of convex polygons by this division is not ASP.
\end{itemize}

\begin{itemize}
\item In the case of $\beta = 90^ \circ $ in Method 2 (see Figure~\ref{fig05}(d)), 
there is the convex pentagonal monotile $P_{2}$ belonging to both the Type 2 and 
Type 4 families with the property of line symmetry that cannot distinguish between 
the anterior and posterior sides (see \ref{appA}); therefore, the set of 
three types of convex polygons by this division is not an ASP.
\end{itemize}

\begin{itemize}
\item In the cases of $\alpha = 105^ \circ $ in Method 1 (see Figure~\ref{fig04}(e)) 
and $\beta = 105^ \circ $ in Method 2 (see Figure~\ref{fig05}(e)), the convex 
polygons created by division are identical. The set of two types of convex polygons 
created by this division is not an ASP, because the convex hexagon $P_{3}$ created 
by this division is a convex hexagonal monotile belonging to the Type 1 
family, which can form periodic tiling without using reflected tiles (see \ref{appB}).
\end{itemize}

\begin{figure}[htbp]
\centering\includegraphics[width=14.5cm,clip]{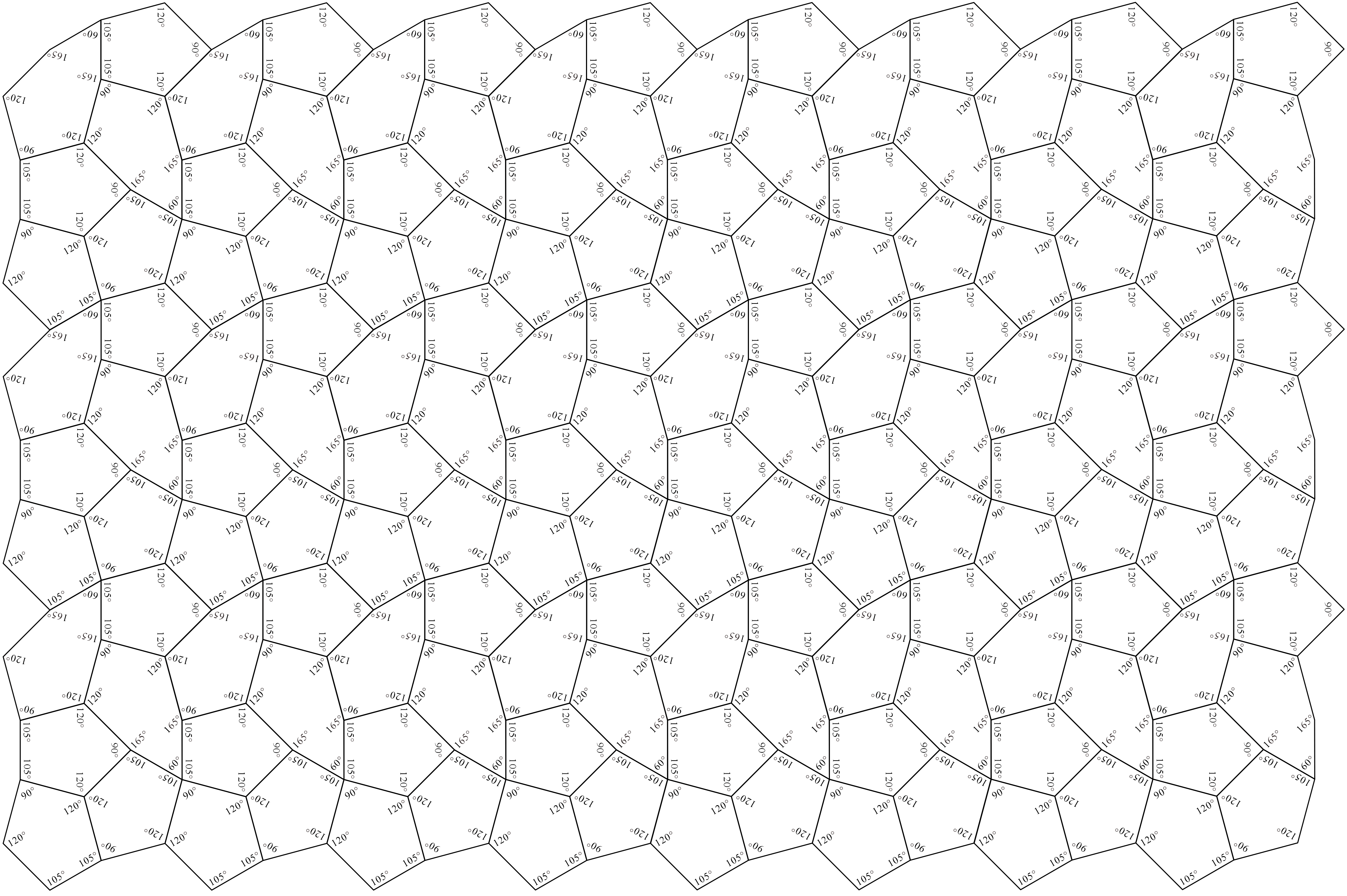} 
\caption{{\small 
Periodic tiling with $P_{1}$ and $P_{2}$ with $\alpha = 165^ \circ $ 
in Method 1 ($\beta = 165^ \circ $ in Method 2).
}
\label{fig06}
}
\end{figure}

In addition, in the cases of $\alpha = 165^ \circ $ in Method 1 (see Figure~\ref{fig04}(f)) 
and $\beta = 165^ \circ $ in Method 2 (see Figure~\ref{fig05}(f)), the three 
convex polygons created by division are identical. We confirmed that each 
convex polygon in these cases cannot generate a tiling even if the use of reflected 
tiles is allowed during the tiling generation process. In other words, 
$P_{1}$, $P_{2}$, and $P_{3}$ in this case are not monotiles. However, the 
convex pentagon $P_{1}$ and convex hexagon $P_{2}$ can form periodic tilings 
as shown in Figure~\ref{fig06}. Therefore, the set of three convex polygons in the 
case of $\alpha = 165^ \circ $ in Method 1 ($\beta = 165^ \circ $ in Method 
2) is not an ASP.

\begin{figure}[htbp]
\centering\includegraphics[width=14.5cm,clip]{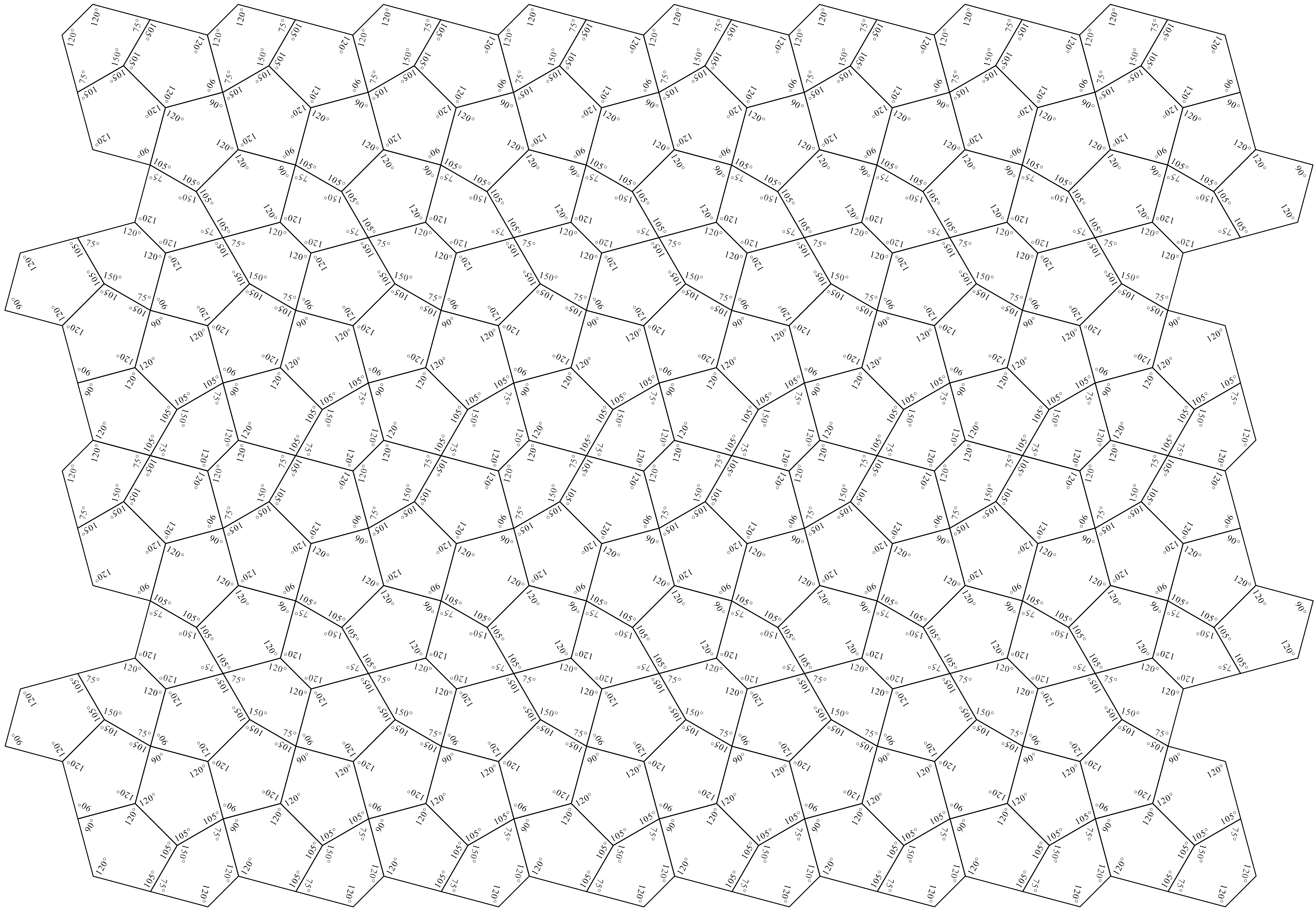} 
\caption{{\small 
Periodic tiling with $P_{1}$ and $P_{3}$ with $\alpha = 180^ \circ $ 
in Method 1.
}
\label{fig07}
}
\end{figure}

In the case of $\alpha = 180^ \circ $ in Method 1 (see Figure~\ref{fig04}(g)), we 
confirmed that each convex polygon created by the division cannot generate a 
tiling, even if the use of reflected tiles is allowed during the tiling generation 
process. In other words, $P_{1}$, $P_{2}$, and $P_{3}$ in this case are not 
monotiles. However, the convex pentagons $P_{1}$ and $P_{3}$ can form periodic 
tilings, as shown in Figure~\ref{fig07}. Therefore, the set of three convex polygons 
in the case of $\alpha = 180^ \circ $ in Method 1 is not an ASP. However, 
outside the conditions, if convex pentagon $P_{1}$, convex hexagon $P_{2}$, 
and reflected $P_{2}$ are used, a periodic tiling can be formed (see \ref{appC}). 
Because such cases are possible, the condition of this study that 
``reflected tiles cannot be used during all tiling generation processes" is 
important.

\renewcommand{\figurename}{{\small Figure}}

\renewcommand{\figurename}{{\small Figure}}

\renewcommand{\figurename}{{\small Figure}}

\renewcommand{\figurename}{{\small Figure}}

\begin{figure}[htbp]
\centering\includegraphics[width=12.5cm,clip]{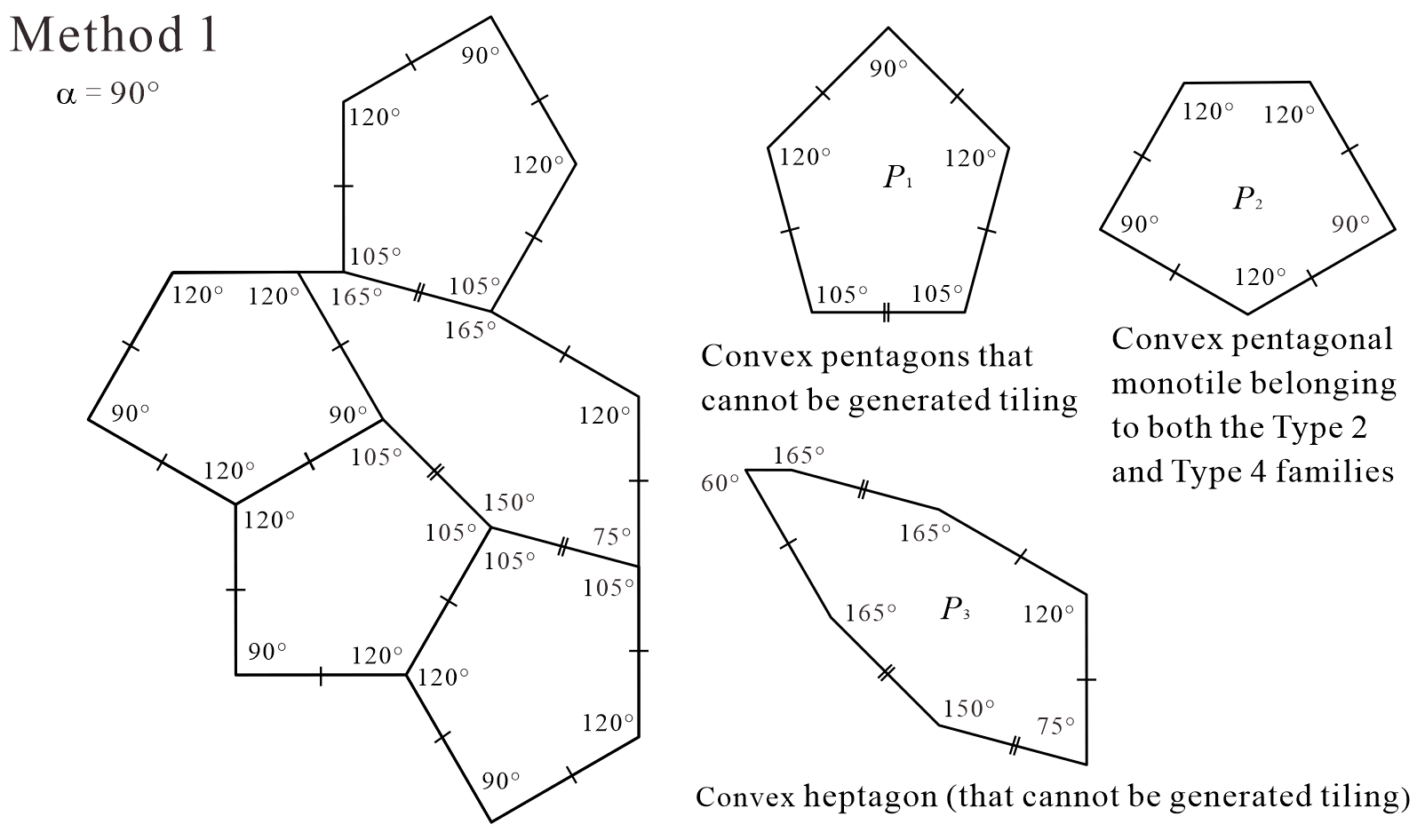} 
\caption{{\small 
Case where the tile set is not aperiodic in Method 1.
}
\label{fig08}
}
\end{figure}

From the cases ``(a) $75^ \circ < \alpha < 105^ \circ $, 
(b) $105^ \circ < \alpha < 165^ \circ $, (c) $165^ \circ < \alpha < 180^ \circ $" 
of Method 1 shown in Figure~\ref{fig04} and ``(a) $90^ \circ < \beta < 105^ \circ $, 
(b) $105^ \circ < \beta < 165^ \circ $, (c) $165^ \circ < \beta < 180^ \circ $" 
of Method 2 shown in Figure~\ref{fig05}, candidates of ASP comprising three types 
of convex polygons can be created. We refer to it as ``candidate" because if we 
generate tilings with convex polygons contained in the tile set, we have to 
make sure that there exist only non-periodic tilings using the shape of 
\mbox{Tile$(1, 1)$}, if and only if all types of convex polygons contained in the tile 
set are used. For example, the division of $\alpha = 90^ \circ $ in 
``(a) $75^ \circ < \alpha < 105^ \circ $" of Method 1 creates a convex 
pentagonal monotile $P_{2}$ belonging to both the Type 2 and Type 4 families 
with line symmetry, as shown in Figure~\ref{fig08}, which can form periodic tiling 
(see \ref{appA}). Thus, the set of three types of convex polygons created by 
the division of $\alpha = 90^ \circ $ in Method 1 is not an ASP.

\begin{figure}[htbp]
\centering\includegraphics[width=12.5cm,clip]{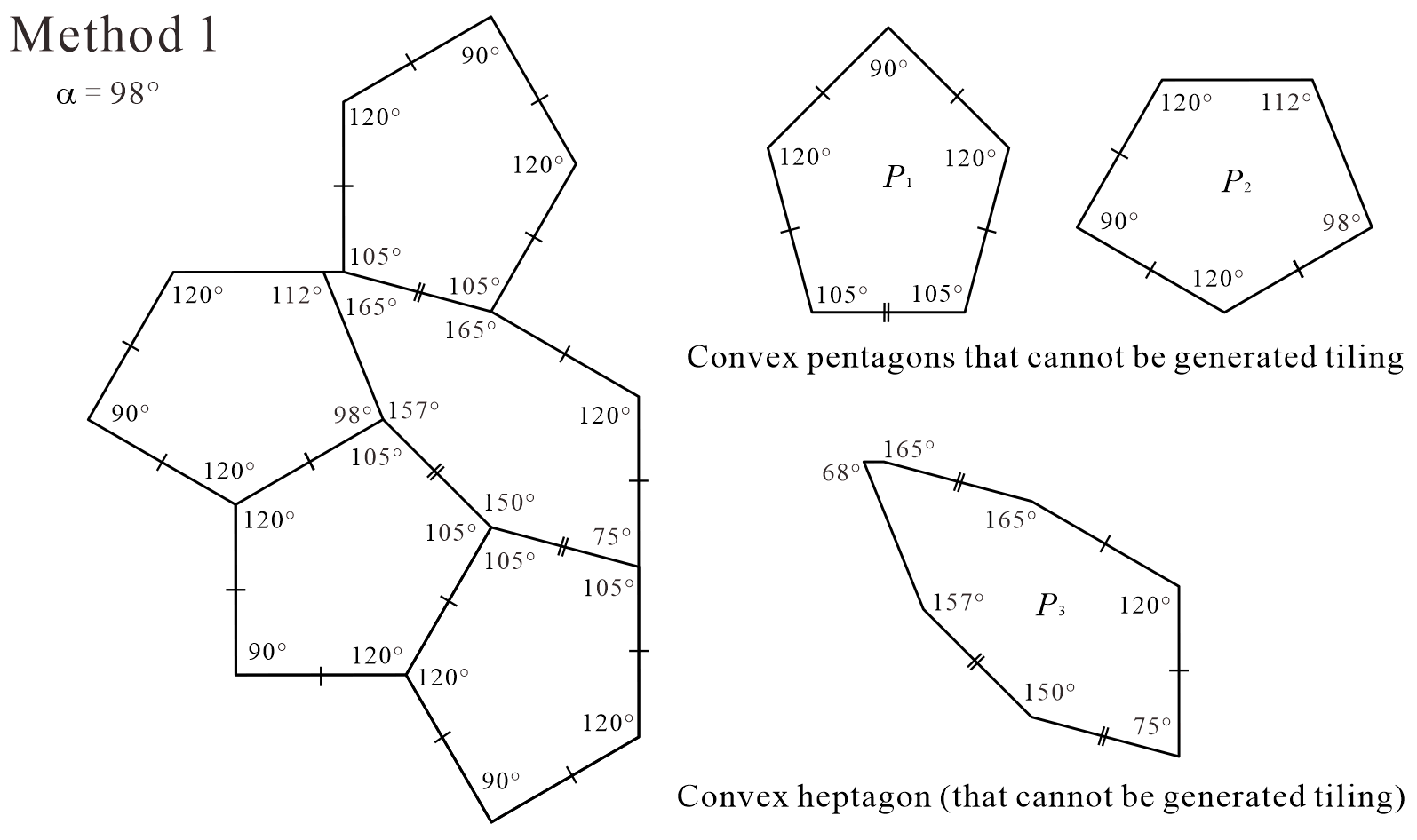} 
\caption{{\small 
Concrete example of an ASP that can be created using Method 1 
(Part 1).
}
\label{fig09}
}
\end{figure}

\begin{figure}[htbp]
\centering\includegraphics[width=12.5cm,clip]{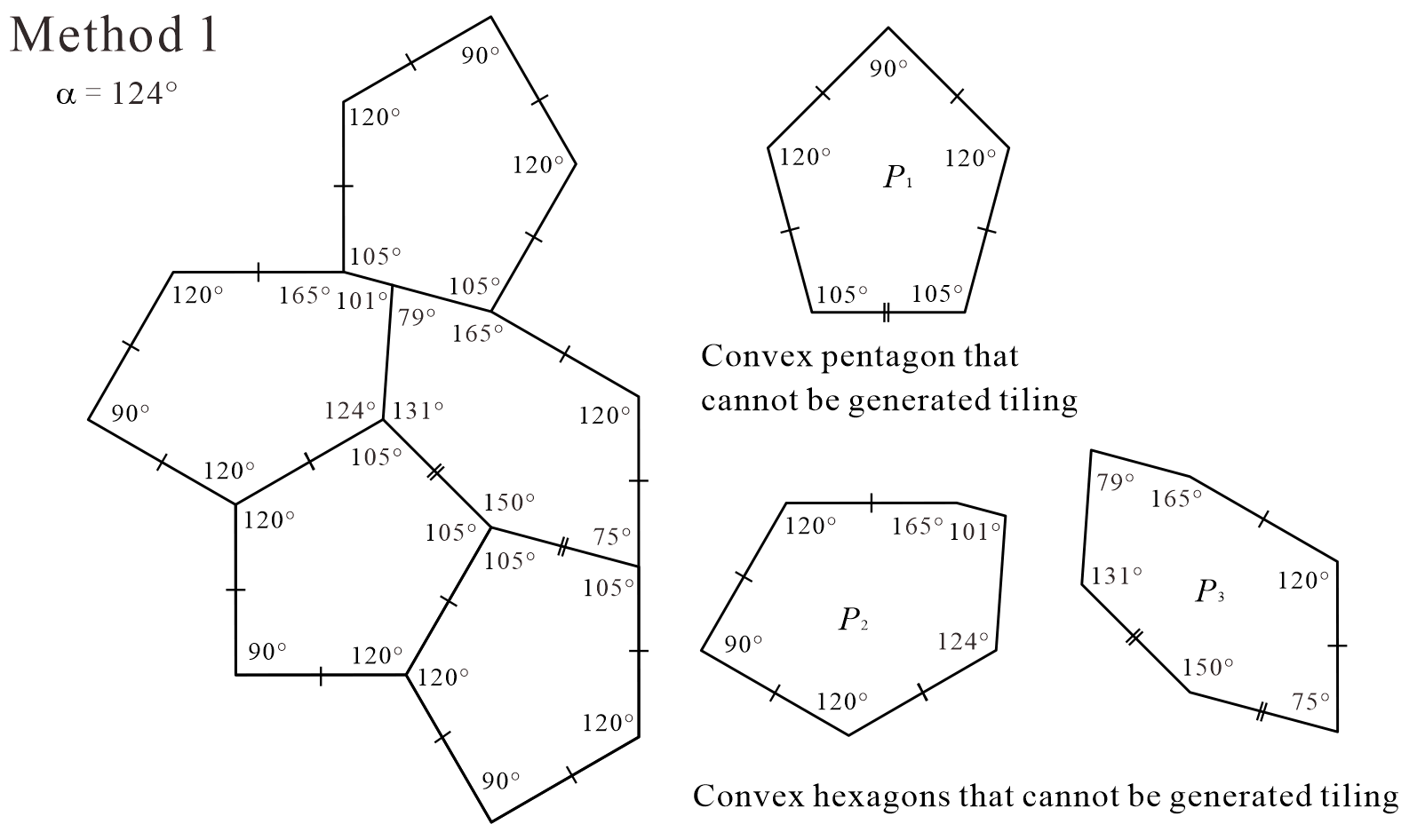} 
\caption{{\small 
Concrete example of an ASP that can be created using Method 1 
(Part 2).
}
\label{fig10}
}
\end{figure}

\begin{figure}[htbp]
\centering\includegraphics[width=15cm,clip]{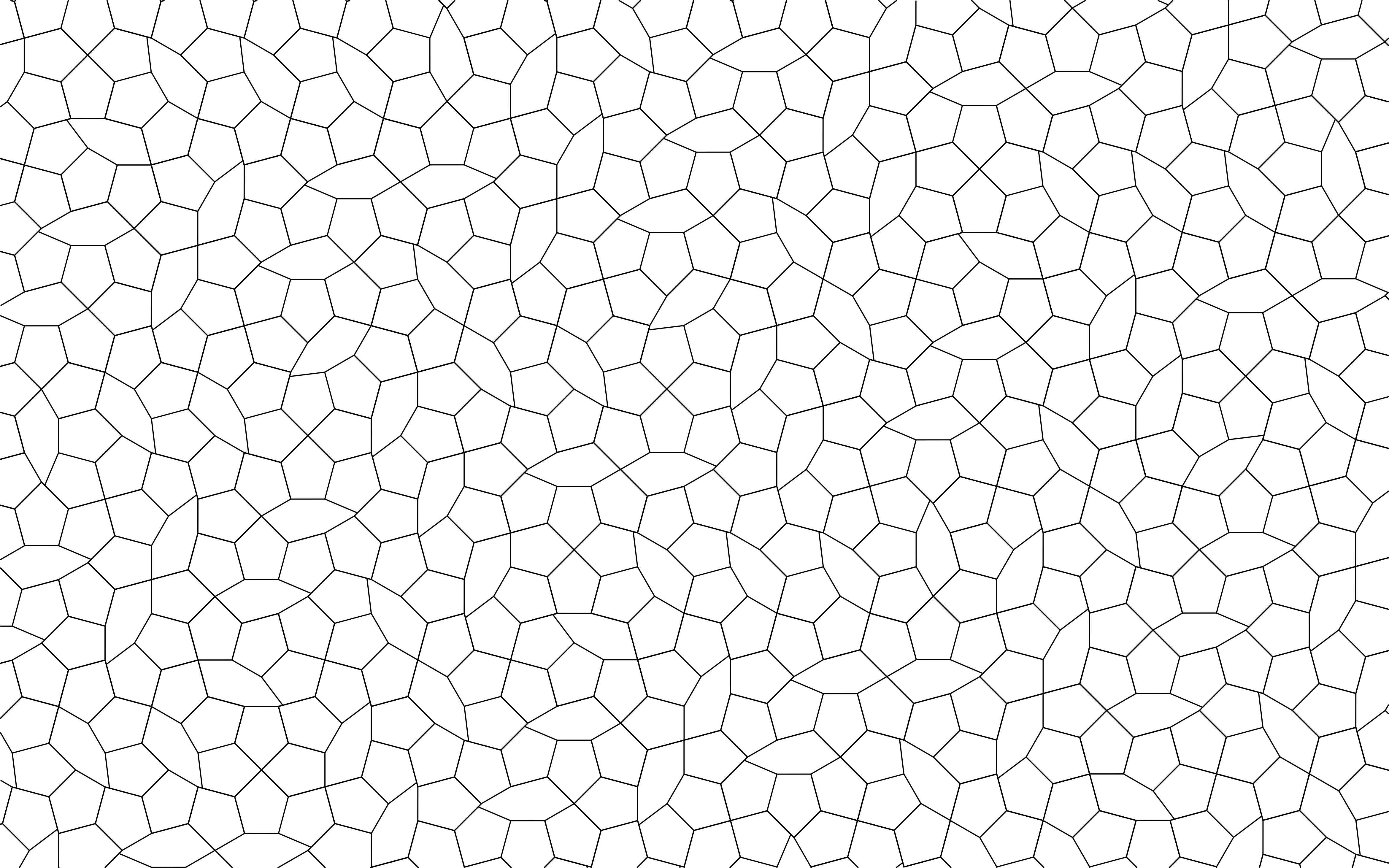} 
\caption{{\small 
Non-periodic tiling generated by the ASP comprising two types of 
convex pentagons and one type of convex heptagon in the case of 
$\alpha = 98^ \circ $ in Method 1.
}
\label{fig11}
}
\end{figure}

\begin{figure}[htbp]
\centering\includegraphics[width=15cm,clip]{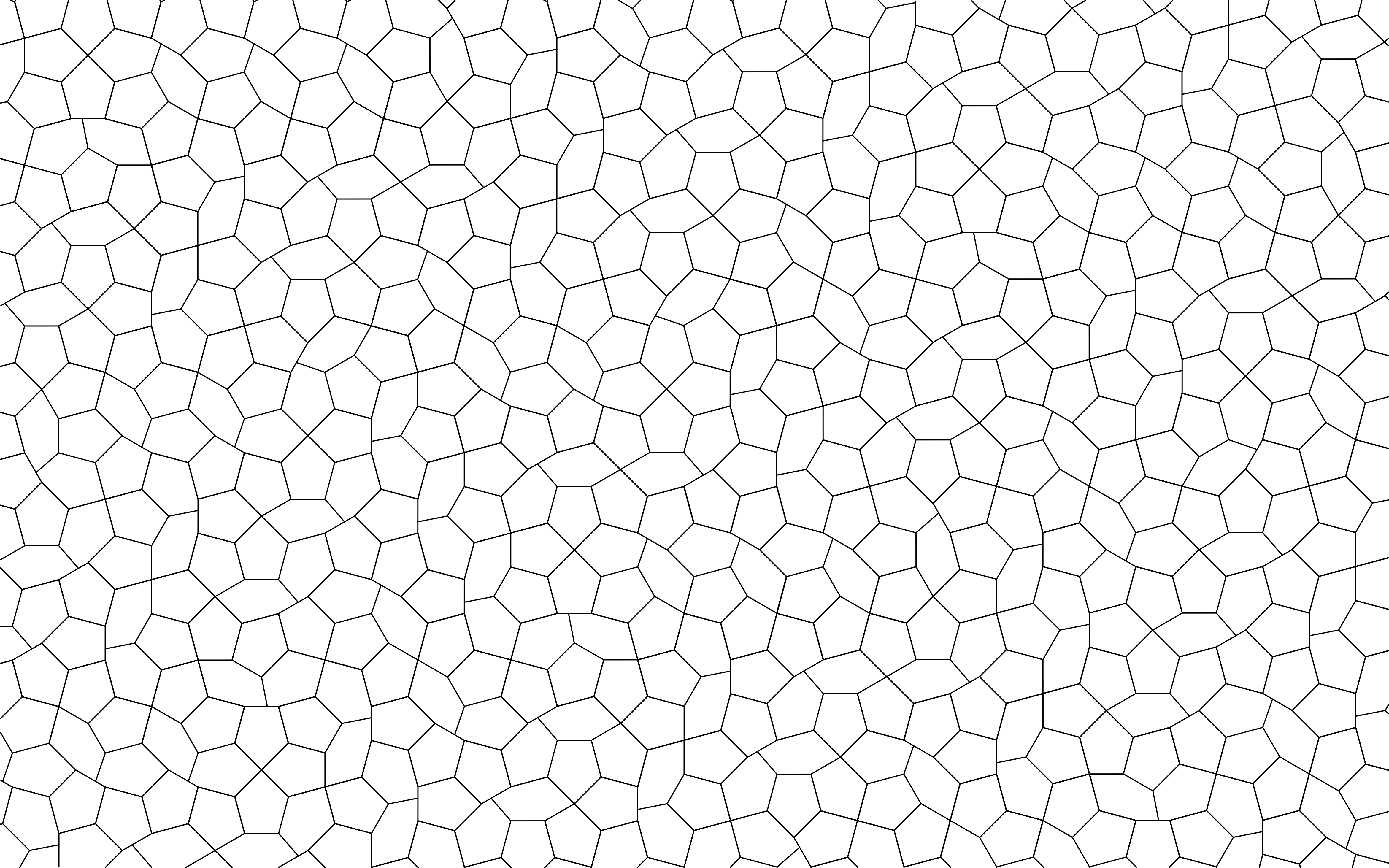} 
\caption{{\small 
Non-periodic tiling generated by the ASP comprising one type of 
convex pentagon and two types of convex hexagons in the case of 
$\alpha = 124^ \circ $ in Method 1.
}
\label{fig12}
}
\end{figure}

As concrete examples of the ``candidate" cases, we investigated the sets of 
three types of convex polygons in the case of $\alpha = 98^ \circ $ in 
Method 1 shown in Figure~\ref{fig09} and in the case of $\alpha = 124^ \circ $ in 
Method 1 shown in Figure~\ref{fig10}. We confirmed that these tile sets are ASPs 
comprising three types of convex polygons. Figures~\ref{fig11} and ~\ref{fig12} show the 
non-periodic tilings generated by the ASP comprising three types of convex 
polygons in the cases of $\alpha = 98^ \circ $ and $\alpha = 124^ \circ $ in 
Method 1.

\renewcommand{\figurename}{{\small Figure}}

\begin{figure}[htbp]
\centering\includegraphics[width=10cm,clip]{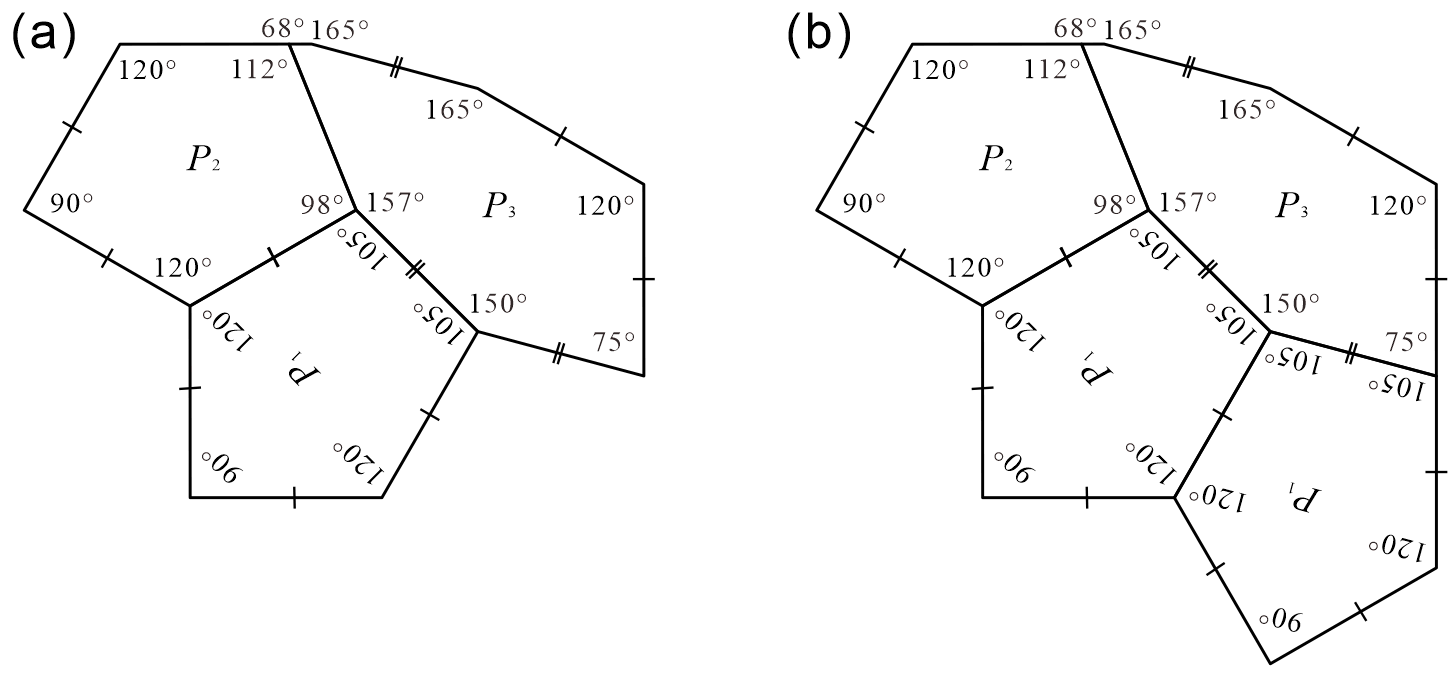} 
\caption{{\small 
Combination focusing on vertices with internal angles of $157^ \circ $ 
and $150^ \circ $ in the case of $\alpha = 98^ \circ $ in Method 1.
}
\label{fig13}
}
\end{figure}

Using the case of the three convex polygons created by the division of 
$\alpha = 98^ \circ $ in Method 1, as shown in Figure~\ref{fig09}, some of the key 
points of the investigation (identification) method are briefly described. 
To form tilings, two or more tile vertices must be concentrated at one point 
and the sum of the interior angles of the vertices at the concentration 
point must be $180^ \circ$ or $360^ \circ$. Therefore, as there were no convex 
pentagons with interior angles of $65^ \circ$ or less, $120^ \circ$, $150^ \circ$, 
$157^ \circ$, and $165^ \circ$ cannot be used in a concentration with a sum of 
$180^ \circ$. From the balance with other interior angles, we observed that an 
even number of vertices with interior angles of $68^ \circ$, $98^ \circ$, $112^ \circ$, 
and $157^ \circ$ must always be concentrated at one point, because the place 
value of one should be zero or five to make $180^ \circ$ or $360^ \circ$. Based on 
these considerations, it was concluded that the concentration of vertices 
with an interior angle of $157^ \circ$ of $P_{3}$ was impossible, except for the 
arrangements using $P_{1}$ and $P_{2}$, as shown in Figure~\ref{fig13}(a). 
Subsequently, note the point in Figure~\ref{fig13}(a) where $P_{3}$ has an internal 
angle of $150^ \circ$. If a tiling is formed, the only method to obtain a sum of 
$360^ \circ$ is to place the vertex of the tile with an interior angle of $105^ \circ$ 
at the concentrating point. Figure~\ref{fig13}(b) shows the only combination that 
can form a structure in which the tiling does not collapse. Then, for the 
positions with vertices having internal angles of $165^ \circ$ of $P_{3}$ in 
Figure~\ref{fig13}(b), the only combination that can construct a structure in which 
the tiling does not collapse using $P_{1}$, $P_{2}$, and $P_{3}$ in Figure~\ref{fig09} is 
that of dividing \mbox{Tile$(1, 1)$} (using the shape of \mbox{Tile$(1, 1)$}), as shown in 
Figure~\ref{fig09}.

If an ASP comprising three types of convex polygons can be constructed from the 
cases ``(a) $75^ \circ < \alpha < 105^ \circ $, (b) $105^ \circ < \alpha < 165^ \circ $, 
(c) $165^ \circ < \alpha < 180^ \circ $" of Method 1 in Figure~\ref{fig04} 
or ``(a) $90^ \circ < \beta < 105^ \circ $, (b) $105^ \circ < \beta < 165^ \circ $, 
(c) $165^ \circ < \beta < 180^ \circ $" of Method 2 in Figure~\ref{fig05}, 
the combination of three types of convex polygons in the set is as follows.

\begin{itemize}
\item The tile set comprising two types of convex pentagons and 
one type of convex heptagon.

\begin{itemize}
\item[$\circ$] Method 1: $75^ \circ < \alpha < 105^ \circ $ and $165^ \circ < \alpha < 180^ \circ $

\item[$\circ$] Method 2: $90^ \circ < \beta < 105^ \circ $ and $165^ \circ < \beta < 180^ \circ $
\end{itemize}
\end{itemize}

\renewcommand{\figurename}{{\small Figure}}

\renewcommand{\figurename}{{\small Figure}}

\renewcommand{\figurename}{{\small Figure}}

\renewcommand{\figurename}{{\small Figure}}

\renewcommand{\figurename}{{\small Figure}}

 {\color{white}. }

\renewcommand{\figurename}{{\small Figure}}
\begin{figure}[htbp]
\centering\includegraphics[width=15cm,clip]{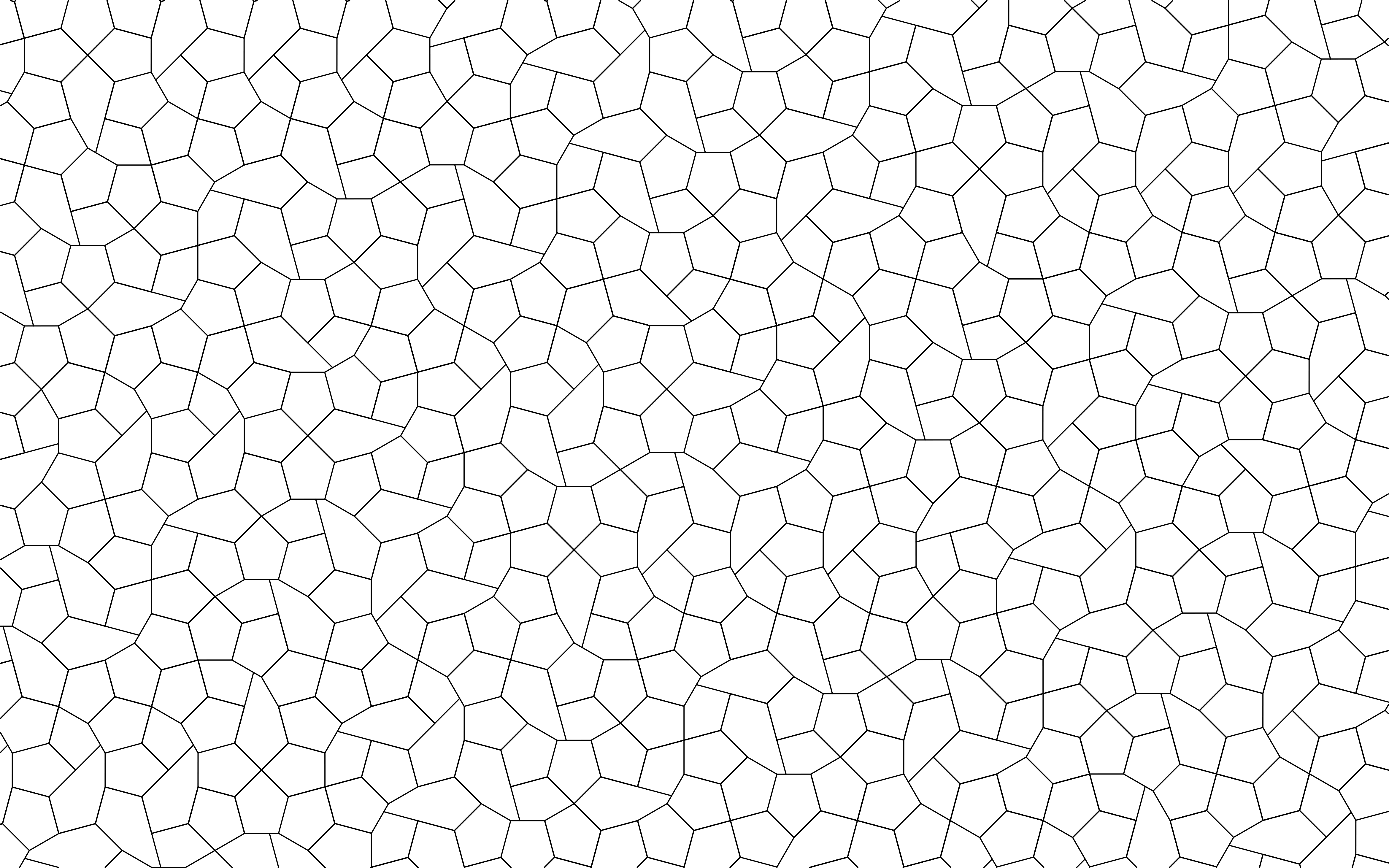} 
\caption{{\small 
Non-periodic tiling generated by the ASP comprising two types of 
convex pentagons and one type of convex hexagon in the case of 
$\beta = 180^ \circ $ in Method 2.
}
\label{fig14}
}
\end{figure}

\begin{itemize}
\item The tile set comprising one type of convex pentagon and two types 
of convex hexagons.

\begin{itemize}
\item[$\circ$] Method 1: $105^ \circ < \alpha < 165^ \circ $

\item[$\circ$] Method 2: $105^ \circ < \beta < 165^ \circ $
\end{itemize}
\end{itemize}

However, in the case of $\beta = 180^ \circ $ in Method 2 (see Figure~\ref{fig05}(g)), 
two types of convex pentagons ($P_{1}$ and $P_{3})$ and one type of convex 
hexagon ($P_{2})$ were created, and we confirmed that this tile set is an ASP 
comprising three types of convex polygons. Figure~\ref{fig14} shows the non-periodic 
tiling generated by the ASP comprising three types of convex polygons (two 
types of convex pentagons and one type of convex hexagon) in the case of 
$\beta = 180^ \circ $ in Method 2.

%%%%%%%%%%%%%%%%%%%%%%%%%%%%%%%%%%%%%%%%%%%%%%%%%%%%%%%%%%%%%%%%%%%%%%

\subsection{Method 3}
\label{subsecl2_2}

\begin{figure}[htbp]
\centering\includegraphics[width=12.5cm,clip]{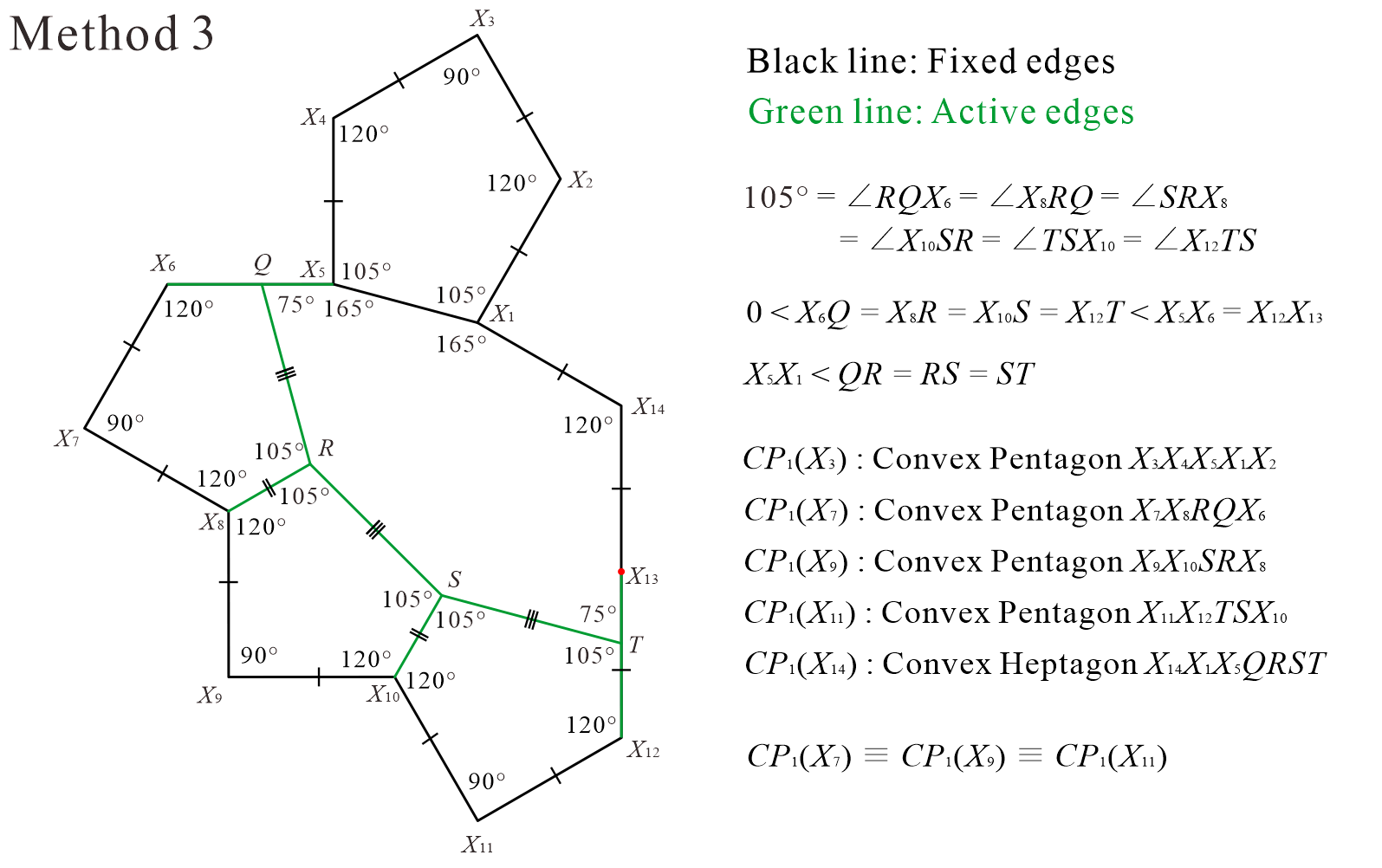} 
\caption{{\small 
Case in Method 3.
}
\label{fig15}
}
\end{figure}

The third division method is Method 3, as shown in Figure~\ref{fig15}. This can be 
regarded as a modified version of the division shown in Figure~\ref{fig04}(e). The 
convex pentagon $P_{1}$ in Figure~\ref{fig04}(e), which corresponds to the vertices 
$X_{7}$, $X_{9}$, and $X_{11}$ in Figure~\ref{fig15}, can be considered as transformed, 
maintaining the same values of the interior angles.

In Figure~\ref{fig15}, the point $Q$ moves on the edge $X_{5}X_{6}$ and the point $T$ 
moves on the edge $X_{12}X_{13}$ of \mbox{Tile$(1, 1)$} in the range 
``$0 < X_{6}Q = X_{10}T < X_{5}X_{6} $." In Figure~\ref{fig15}, the convex pentagon 
$CP_{1}(X_{3})$ with the points ``$X_{3}$, $X_{4}$, $X_{5}$, $X_{1}$, $X_{2}$" as vertices is 
fixed shape, which is the same as the convex pentagon $P_{1}$ in Figure~\ref{fig04}. 
Additionally, the convex pentagon $CP_{1}(X_{7})$ with the points ``$X_{7}$, $X_{8}$, $R$, 
$Q$, $X_{6}$" as vertices, the convex pentagon $CP_{1}(X_{9})$ with the points 
``$X_{9}$, $X_{10}$, $S$, $R$, $X_{8}$" as vertices, and the convex pentagon 
$CP_{1}(X_{11})$ with the points ``$X_{11}$, $X_{12}$, $T$, $S$, $X_{10}$" as vertices 
are congruent convex pentagons with interior angles always at 
``$90^ \circ , 120^ \circ , 105^ \circ , 105^ \circ , 120^ \circ $" and with line 
symmetry. The convex polygon $CP_{1}(X_{14})$ with the points ``$X_{14}$, $X_{1}$, 
$X_{5}$, $Q$, $R$, $S$, $T$" as vertices is always a convex heptagon.

\begin{figure}[htbp]
\centering\includegraphics[width=12.5cm,clip]{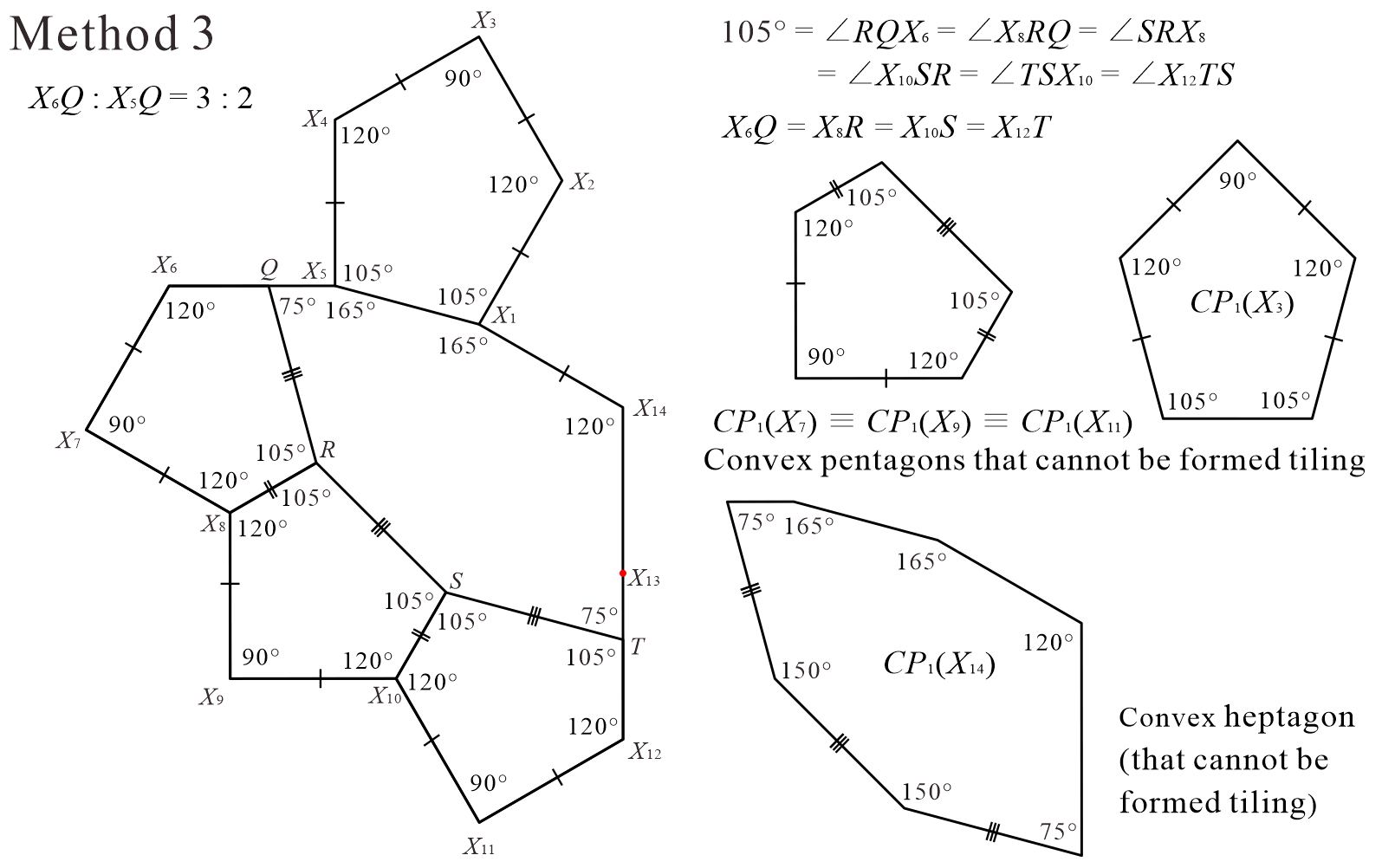} 
\caption{{\small 
Concrete example of an ASP that can be created using Method 3.
}
\label{fig16}
}
\end{figure}

\begin{figure}[htbp]
\centering\includegraphics[width=15cm,clip]{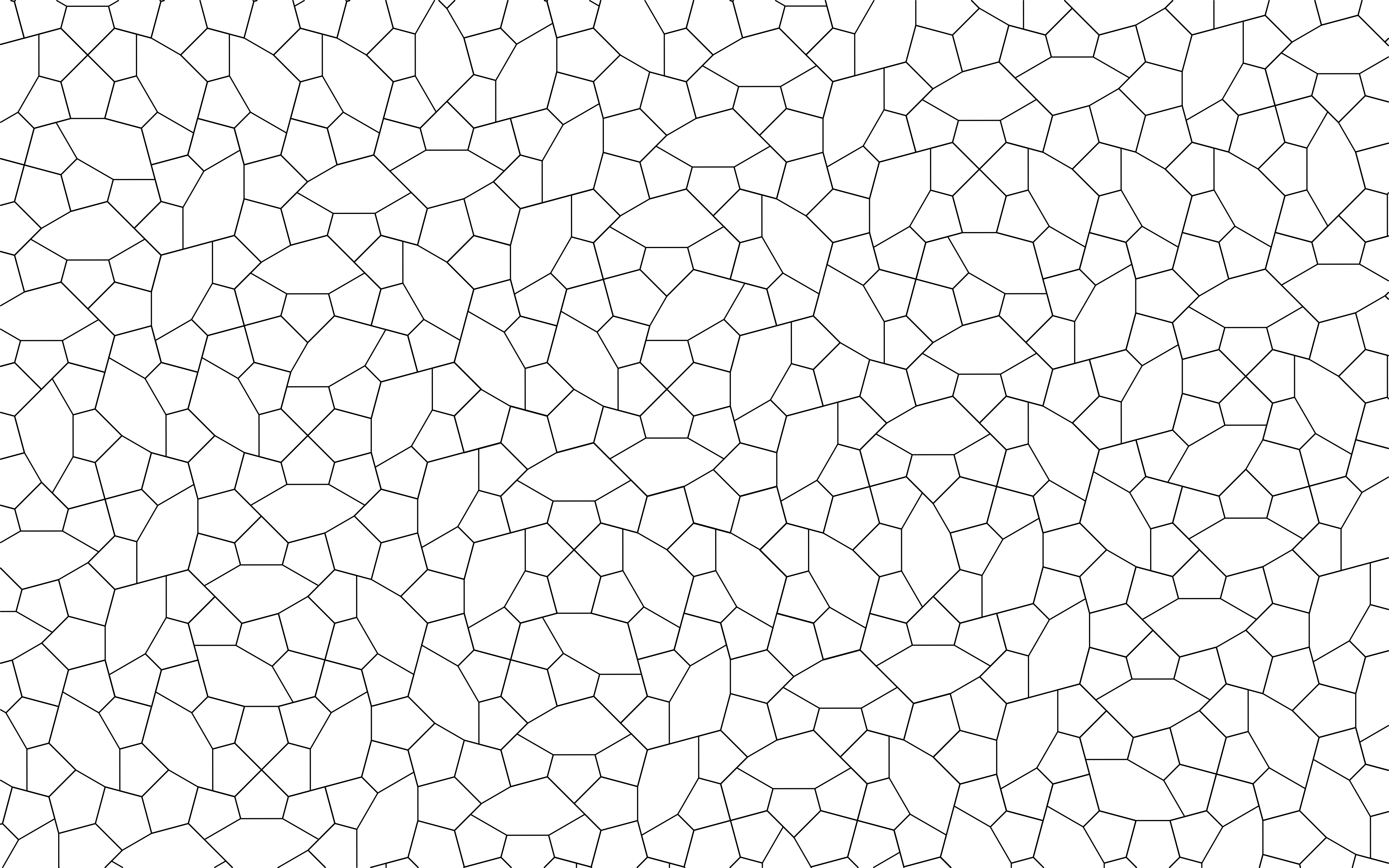} 
\caption{{\small 
Non-periodic tiling generated by the ASP comprising two types of 
convex pentagons and one type of convex heptagon in Figure~\ref{fig16}.
}
\label{fig17}
}
\end{figure}

As a concrete example, we investigated using the case where the ratio of the 
edges is ``$X_{6}Q:X_{5}Q = 3:2$" as shown in Figure~\ref{fig16}. As a result, we 
confirmed that this tile set is an ASP comprising three types of convex 
polygons (two types of convex pentagons and one type of convex heptagon). 
Figure~\ref{fig17} shows the non-periodic tiling generated by the ASP comprising two 
types of convex pentagons and one type of convex heptagon in Figure~\ref{fig16}.

\bigskip
\bigskip
\renewcommand{\figurename}{{\small Figure}}

\renewcommand{\figurename}{{\small Figure}}

\renewcommand{\figurename}{{\small Figure}}

From the process of investigation (identifying) for the case shown in Figure~\ref{fig16}, 
we conjecture that in the case of $X_{6}X_{7} \ne QR$, the tile set of 
three convex polygons in Figure~\ref{fig15} is an ASP, comprising two types of convex 
pentagons and one type of convex heptagon. When the conjecture of $X_{6}X_{7} \ne QR$ 
was not directly applicable, we individually confirmed that the tile set in the case 
of $X_{6}X_{7} = QR$ was an ASP comprising two types of convex 
pentagons and one type of convex heptagon.

%%%%%%%%%%%%%%%%%%%%%%%%%%%%%%%%%%%%%%%%%%%%%%%%%%%%%%%%%%%%%%%%%%%%%%

\subsection{Method 4}
\label{subsec2_3}

\begin{figure}[htbp]
\centering\includegraphics[width=15cm,clip]{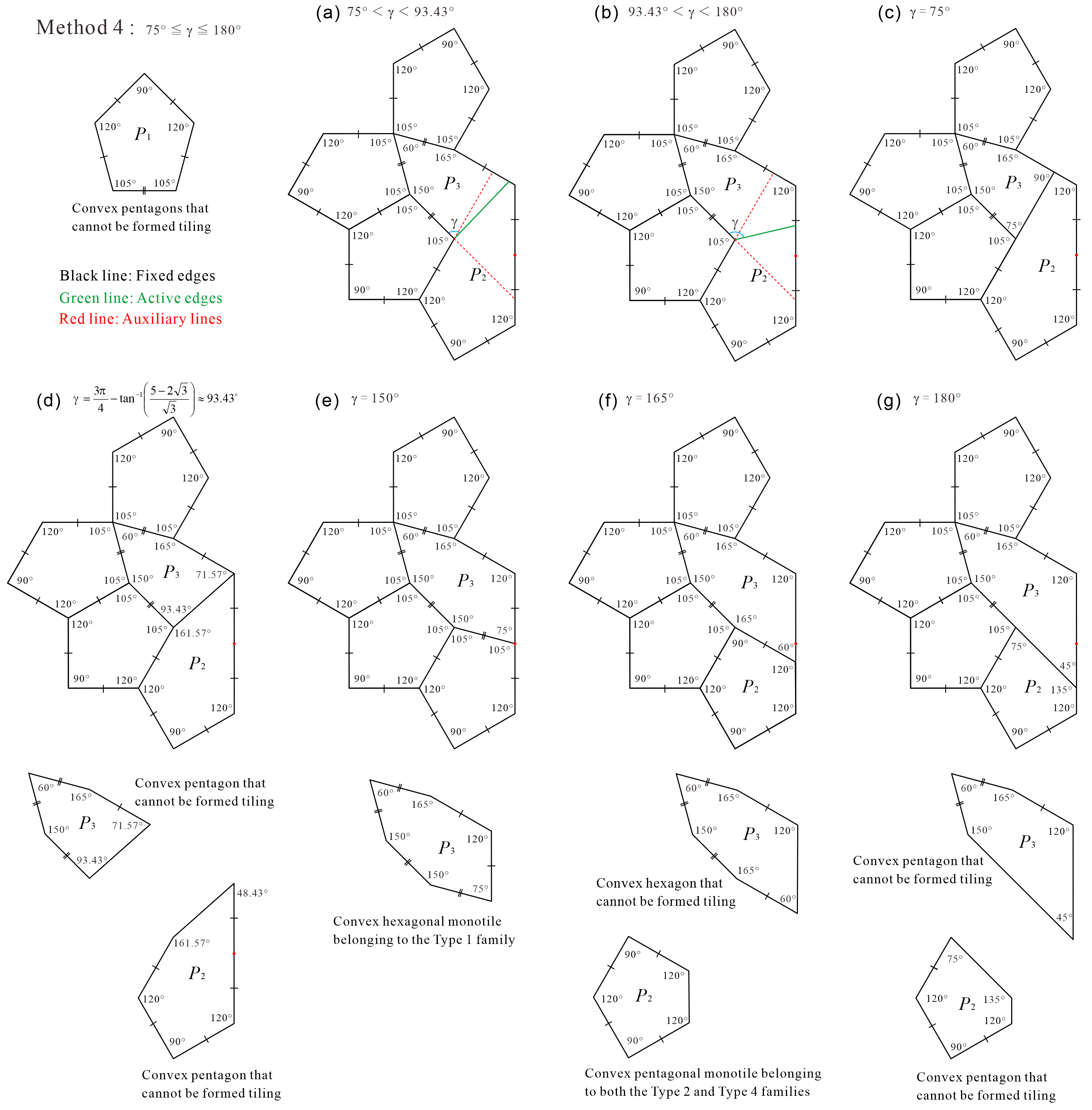} 
\caption{{\small 
Cases in Method 4.
}
\label{fig18v2}
}
\end{figure}

The fourth division method is Method 4, as shown in Figure~\ref{fig18v2}. 
The results of the tile sets (five convex polygons containing three $P_{1}$) created 
by the division of Method 4 differ depending on the value of the parameter  
$\gamma$ shown in Figure~\ref{fig18v2}.

First, from the cases ``(a) $75^ \circ < \gamma < 93.43^ \circ $, 
(b) $93.43^ \circ < \alpha < 180^ \circ $"\footnote{ 
When $\gamma = (3\pi / 4) + \tan ^{ - 1}\left( {{\left( {5 - 2\sqrt 3 } \right)} 
\mathord{\left/ {\vphantom {{\left( {5 - 2\sqrt 3 } \right)} 3}} \right. 
\kern-\nulldelimiterspace} \sqrt 3} \right) \approx 93.43^\circ $  in Figure~\ref{fig18v2}(d),  
$P_{2}$ and $P_{3}$ are both convex pentagons.}
of Method 4 shown in Figure~\ref{fig18v2}, candidates of ASP comprising 
three types of convex polygons can be created. 
We refer to it as ``candidate" because if we generate tilings with convex polygons 
contained in the tile set, we have to make sure that there exist only non-periodic 
tilings using the shape of \mbox{Tile$(1, 1)$}, if and only if all types of convex polygons 
contained in the tile set are used. If an ASP comprising three types of convex 
polygons can be constructed from the cases ``(a) $75^ \circ < \gamma < 93.43^ \circ $, 
(b) $93.43^ \circ < \alpha < 180^ \circ $" of Method 4, it is the tile set comprising two types 
of convex pentagons and one type of convex hexagon.

The cases ``(e) $\gamma = 150^ \circ $, (f) $\gamma = 165^ \circ $" of Method 4 
shown in Figure~\ref{fig18v2} are not ASPs. The reasons are shown below.

\begin{itemize}
\item In the case of $\gamma = 150^ \circ $ in Method 4 (see Figure~\ref{fig18v2}(e)), 
there is the convex hexagonal monotile $P_{3}$ belonging to the Type 1 family that 
can form periodic tiling without using reflected tiles (see \ref{appB}); therefore, 
the set of two types of convex polygons by this division is not ASP. 
Note that this division is the same as for $\alpha = 105^ \circ $ in Method 1 
(see Figure~\ref{fig04}(e)).
\end{itemize}

\begin{itemize}
\item In the case of $\gamma = 165^ \circ $ in Method 4 (see Figure~\ref{fig18v2}(f)), 
there is the convex pentagonal monotile $P_{2}$ belonging to both the Type 2 and 
Type 4 families with the property of line symmetry that cannot distinguish between 
the anterior and posterior sides (see \ref{appA}); therefore, the set of 
three types of convex polygons by this division is not an ASP.
\end{itemize}

By contrast, for the cases ``(c)  $\gamma = 75^ \circ $, 
(d) $\gamma = (3\pi / 4) + \tan ^{ - 1}\left( {{\left( {5 - 2\sqrt 3 } \right)} 
\mathord{\left/ {\vphantom {{\left( {5 - 2\sqrt 3 } \right)} 3}} \right. 
\kern-\nulldelimiterspace} \sqrt 3} \right) \approx 93.43^\circ $, (g) $\gamma = 180^ \circ $" 
of Method 4 shown in Figure~\ref{fig18v2}, under the condition that the tiles in the tile set 
cannot be reflected (the use of reflected tiles is not allowed) during the tiling generation 
process, we confirmed that these tile sets are ASPs comprising three types of 
convex pentagons.

\renewcommand{\figurename}{{\small Figure}}

\begin{figure}[htbp]
\centering\includegraphics[width=12.5cm,clip]{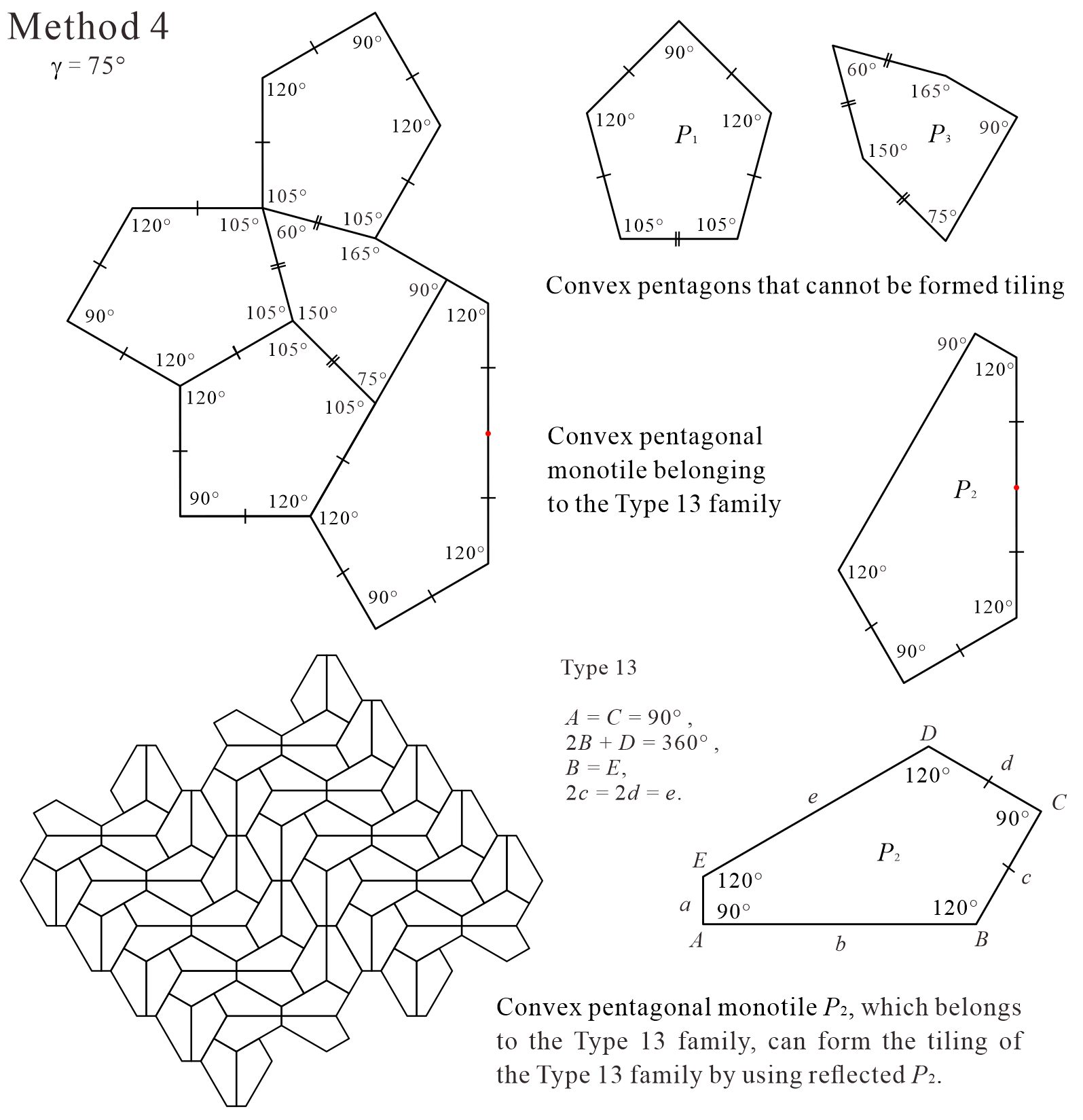} 
\caption{{\small 
Case of $\gamma = 75^ \circ $ in Method 4.
}
\label{fig19v2}
}
\end{figure}

It can be considered that the tile set of $P_{1}$, $P_{2}$, and $P_{3}$ in the case of 
$\gamma = 75^ \circ $ in Method 4 (see Figure~\ref{fig19v2}) is not an ASP, 
because the set contains the convex pentagonal monotile $P_{2}$ belonging 
to the Type 13 family. However, in this study, we assumed that the reflected 
tiles cannot be used during all tiling generation processes. $P_{2}$, which belongs to 
the Type 13 family, can form the (representative) tiling of the Type 13 family 
shown in Figure~\ref{fig19v2} by using reflected $P_{2}$. However, the tile set of the 
three types of convex pentagons created by this division does not contain the 
reflected $P_{2}$. We actually confirmed that tiling cannot be generated with only 
this $P_{2}$ under the condition that ``reflected tiles cannot be used during all 
tiling generation processes"\footnote{
If a convex pentagonal monotile belongs to at least one of the Type 1 or 3--6 
families, it can generate tilings without the use of reflected tiles (see \ref{appA}). 
The convex pentagonal monotiles belonging to the Type 13 families and those 
belonging to the Type 1 or 3--6 families are disjoint \cite{G_and_S_1987, Sugimoto_2016, 
Sugimoto_2017, Sugimoto_2018, Sugimoto_2025, wiki_pentagon_tiling}.
}. We examined whether tiling is possible using $P_{1}$, $P_{2}$, and $P_{3}$, as 
shown in Figure~\ref{fig19v2}, and confirmed that the only combination that can 
construct a structure in which the tiling does not collapse using $P_{1}$, 
$P_{2}$, and $P_{3}$ is that of dividing \mbox{Tile$(1, 1)$} (using the shape of 
\mbox{Tile$(1, 1)$}), as shown in Figure~\ref{fig19v2}. Therefore, under the condition 
that the tiles in the tile set cannot be reflected during the tiling generation process, 
the tile set of the three types of convex pentagons that can be created 
using the case of $\gamma = 75^ \circ $ in Method 4 is an ASP.

\begin{figure}[htbp]
\centering\includegraphics[width=15cm,clip]{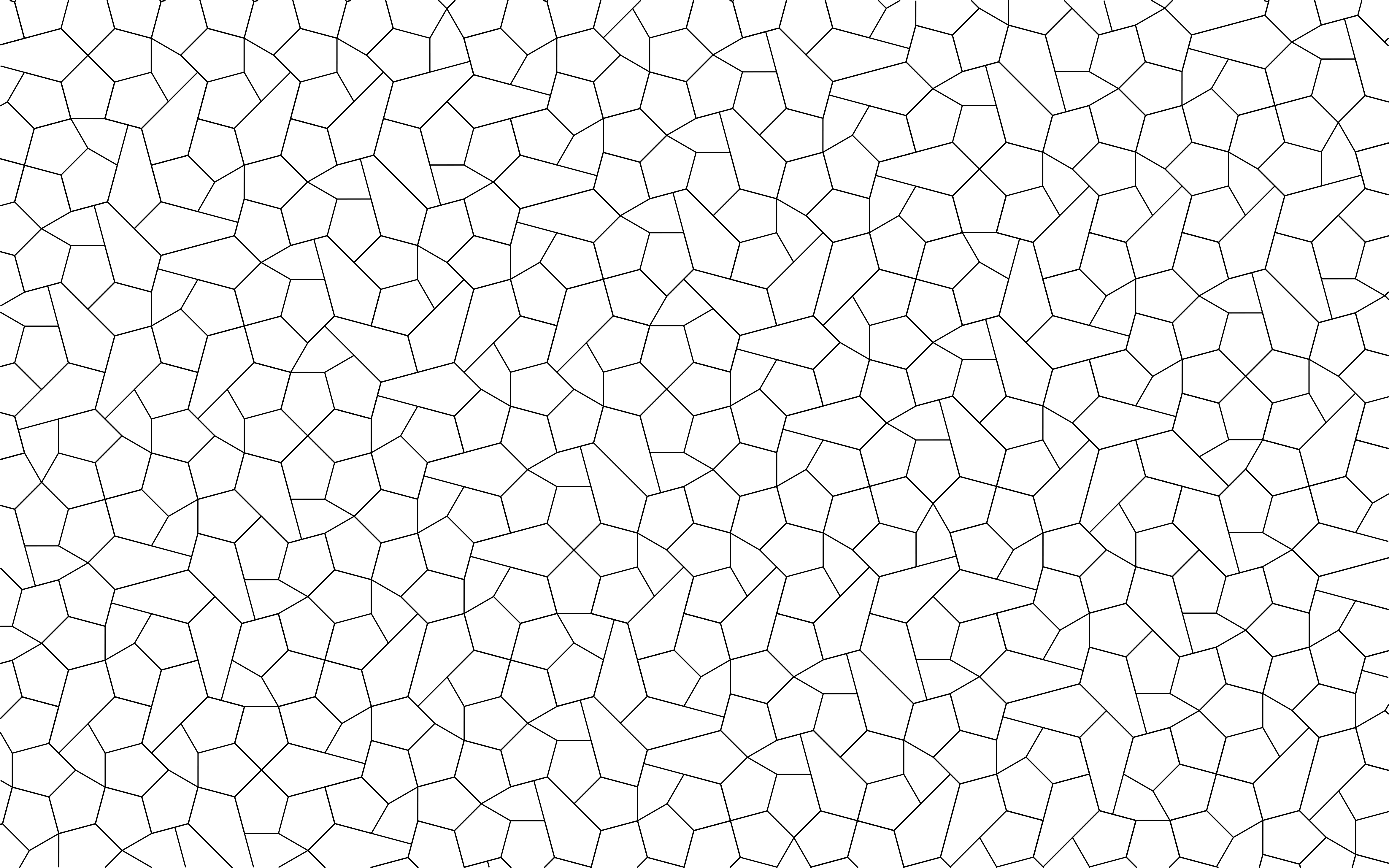} 
\caption{{\small 
Non-periodic tiling generated by the ASP comprising three types 
of convex pentagons in  in the case of 
$\gamma = 75^ \circ $ 
in Method 4.
}
\label{fig20v2}
}
\end{figure}

\begin{figure}[htbp]
\centering\includegraphics[width=15cm,clip]{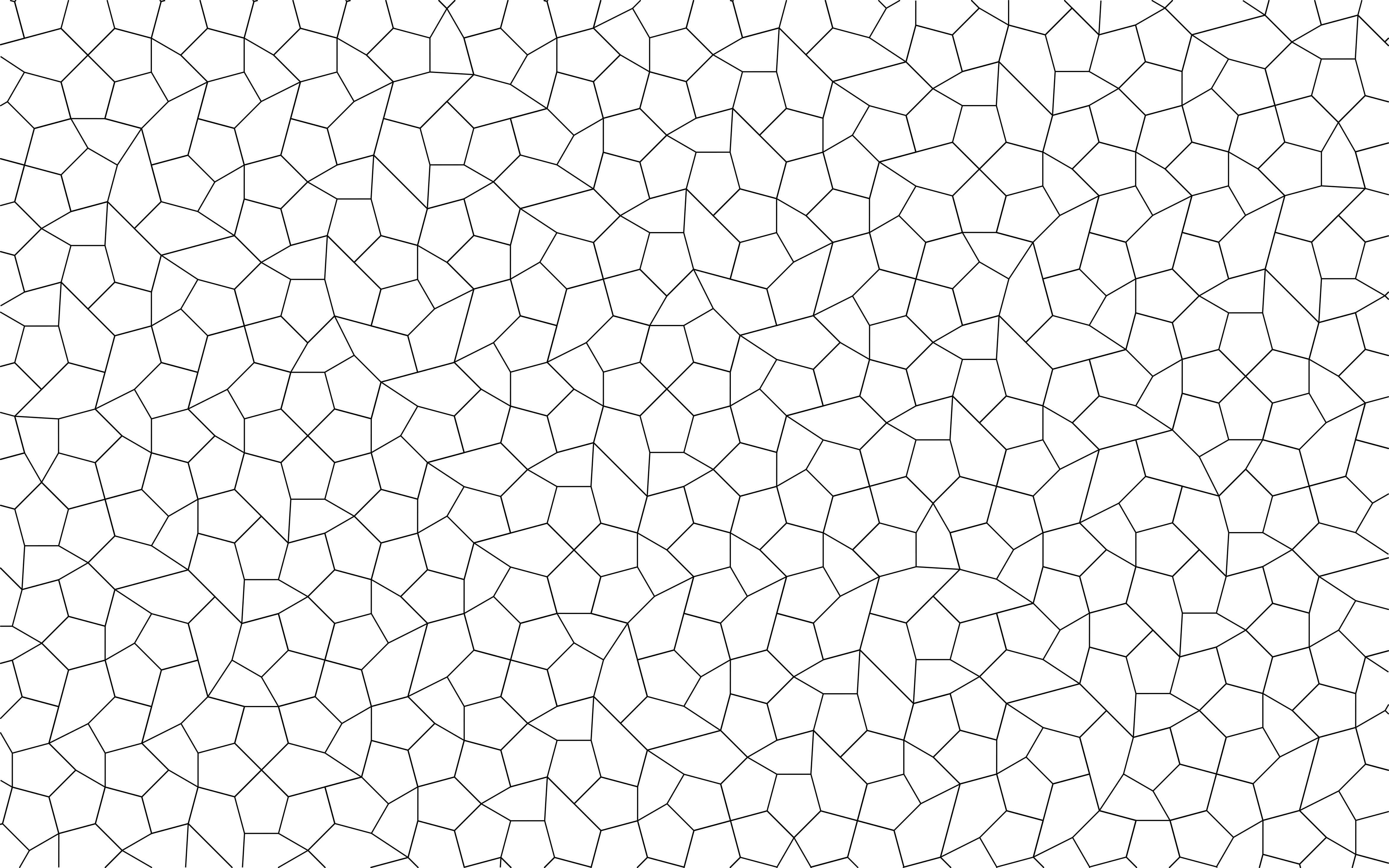} 
\caption{{\small 
Non-periodic tiling generated by the ASP comprising three types 
of convex pentagons in  in the case of 
$\gamma = (3\pi / 4) + \tan ^{ - 1}\left( {{\left( {5 - 2\sqrt 3 } \right)} 
\mathord{\left/ {\vphantom {{\left( {5 - 2\sqrt 3 } \right)} 3}} \right. 
\kern-\nulldelimiterspace} \sqrt 3} \right) \approx 93.43^\circ $  
in Method 4.
}
\label{fig21v2}
}
\end{figure}

\begin{figure}[htbp]
\centering\includegraphics[width=15cm,clip]{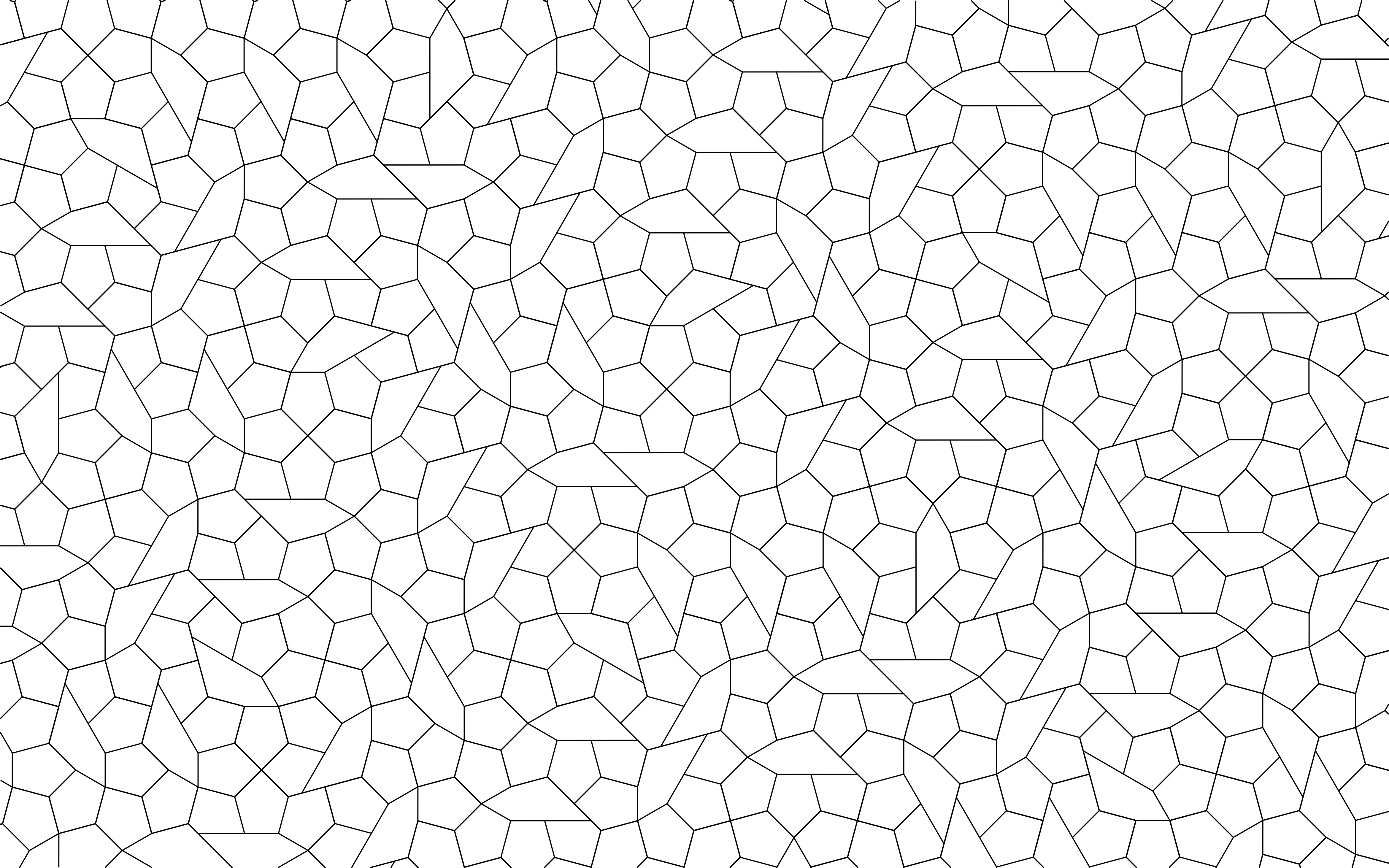} 
\caption{{\small 
Non-periodic tiling generated by the ASP comprising three types 
of convex pentagons in  in the case of 
$\gamma = 180^ \circ $ 
in Method 4.
}
\label{fig22v2}
}
\end{figure}

Figures~\ref{fig20v2}, ~\ref{fig21v2}, and ~\ref{fig22v2} show the non-periodic tilings 
generated by the ASP comprising three types of convex pentagons in the cases of 
$\gamma = 75^ \circ $, $\gamma = (3\pi / 4) + \tan ^{ - 1}\left( {{\left( {5 - 2\sqrt 3 } \right)} 
\mathord{\left/ {\vphantom {{\left( {5 - 2\sqrt 3 } \right)} 3}} \right. 
\kern-\nulldelimiterspace} \sqrt 3} \right) \approx 93.43^\circ $,  and $\gamma = 180^ \circ $ 
in Method 4.

\bigskip
\bigskip
\renewcommand{\figurename}{{\small Figure}}

\renewcommand{\figurename}{{\small Figure}}

\renewcommand{\figurename}{{\small Figure}}

\renewcommand{\figurename}{{\small Figure}}

%%%%%%%%%%%%%%%%%%%%%%%%%%%%%%%%%%%%%%%%%%%%%%%%%%%%%%%%%%%%%%%%%%%%%%
%%%%%%%%%%%%%%%%%%%%%%%%%%%%%%%%%%%%%%%%%%%%%%%%%%%%%%%%%%%%%%%%%%%%%%

\section{Conclusion}
\label{section3}

In this study, we presented the ASPs comprising three types of convex 
polygons other than the ``the set of three convex polygons that are 
aperiodic with no matching condition on the edges" shown in Figure~\ref{fig01}. 
Note that, in this study, we treated the anterior and posterior side tiles as 
different types, even if they are congruent diagrams, and assumed that the 
reflected tiles cannot be used during all tiling generation process with the 
tiles in a tile set.

There appears to be no known ASP comprising two types of convex polygons 
with no matching conditions on the edges (there also seems to be no proof 
that there is no such ASP) \cite{wiki_aperiodic_set, wiki_list_aperiodic}.

Because the ASPs comprising three types of convex polygons shown in this study are 
based on \mbox{Tile$(1, 1)$} shown in \cite{Smith_2024a} and assumes that \mbox{Tile$(1, 1)$} 
is the chiral aperiodic monotile shown in \cite{Smith_2024b}, it is not confirmed that these sets 
are actually ASPs. 

If the types and number of tiles in a tile set can be freely increased under 
the condition that the tiles in the tile set cannot be reflected during the tiling 
generation process, several ASPs comprising multiple types of convex polygons 
can be created. There will be various methods for dividing \mbox{Tile$(1, 1)$} into 
four, five, or more types of convex polygons that are likely to be ASP (see \ref{appD}). 
This is because the convex pentagon can be divided into six types. 
Thus, by such a division, it would be easy to create convex 
pentagons that are not monotiles, or a convex pentagon containing only one 
side of a convex pentagonal monotile that requires the use of reflected tiles during 
the tiling generation process, as in the case of $\gamma = 75^ \circ$ in Method 4. 
However, even if the tile sets that increase the number of tiles in this manner are 
aperiodic, they are not considered interesting. Therefore, we did not 
discuss the effects of increasing the type or number of such tiles in this 
study.

%%%%%%%%%%%%%%%%%%%%%%%%%%%%%%%%%%%%%%%%%%%%%%%%%%%%%%%%%%%%%%%%%%%%%%
%%%%%%%%%%%%%%%%%%%%%%%%%%%%%%%%%%%%%%%%%%%%%%%%%%%%%%%%%%%%%%%%%%%%%%

%%%%%%%%%%%%%%%%%%%%%%%%%%%%%%%%%%%%%%%%%%%%%%%%%%%%%%%%%%%%%%%%%%%%%%
%%%%%%%%%%%%%%%%%%%%%%%%%%%%%%%%%%%%%%%%%%%%%%%%%%%%%%%%%%%%%%%%%%%%%%

\bigskip
\noindent
{\textbf{Acknowledgments.}
The authors would like to express my deep gratitude to Professor Emeritus Yoshio Agaoka, 
Hiroshima University, for his helpful comments and advice during the preparation of 
this manuscript. The authors also received constructive comments from 
Professor Shigeki Akiyama, University of Tsukuba, and Yoshiaki Araki.
\bigskip

%%%%%%%%%%%%%%%%%%%%%%%%%%%%%%%%%%%%%%%%%%%%%%%%%%%%%%%%%%%%%%%%%%%%%%
%%%%%%%%%%%%%%%%%%%%%%%%%%%%%%%%%%%%%%%%%%%%%%%%%%%%%%%%%%%%%%%%%%%%%%

%%%%%%%%%%%%%%%%%%%%%%%%%%%%%%%%%%%%%%%%%%%%%%%%%%%%%%%%%%%%%%%%%%%%%%
%%%%%%%%%%%%%%%%%%%%%%%%%%%%%%%%%%%%%%%%%%%%%%%%%%%%%%%%%%%%%%%%%%%%%%

\bigskip
\appendix
\def\thesection{Appendix \Alph{section}}
\section{Convex pentagonal monotiles}
\label{appA}

Currently, there are 15 known families of convex pentagonal monotiles, each labeled as a  
``Type," as shown in Figure~\ref{fig23v2} \cite{Gardner_1975, G_and_S_1987, Kershner_1968, 
Mann_2015, Rao_2017, Schattschneider_1978, Stein_1985, Sugimoto_2016, Sugimoto_2017, 
Sugimoto_2018, Sugimoto_2025, wiki_pentagon_tiling, Zong_2020}. 
The tilings depicted with each convex pentagon in Figure~\ref{fig23v2} are representative 
of each Type family. A representative tiling can be formed using only the relationships 
that can be derived from the tile conditions of each Type. 
As shown in Figure~\ref{fig23v2}, the representative tilings of the Types 2 and 7--15 families 
use reflected tiles (convex pentagons). By contrast, as shown in Figure~\ref{fig23v2}, 
the representative tilings of the Types 1 and 3--6 families do not use reflected tiles. 
Thus, if a convex pentagonal monotile belongs to at least one of the Type 1 or 3--6 
families, it can generate tilings without the use of reflected tiles \cite{Sugimoto_2025}.

\renewcommand{\figurename}{{\small Figure}}
\begin{figure}[http]
\centering\includegraphics[width=15cm,clip]{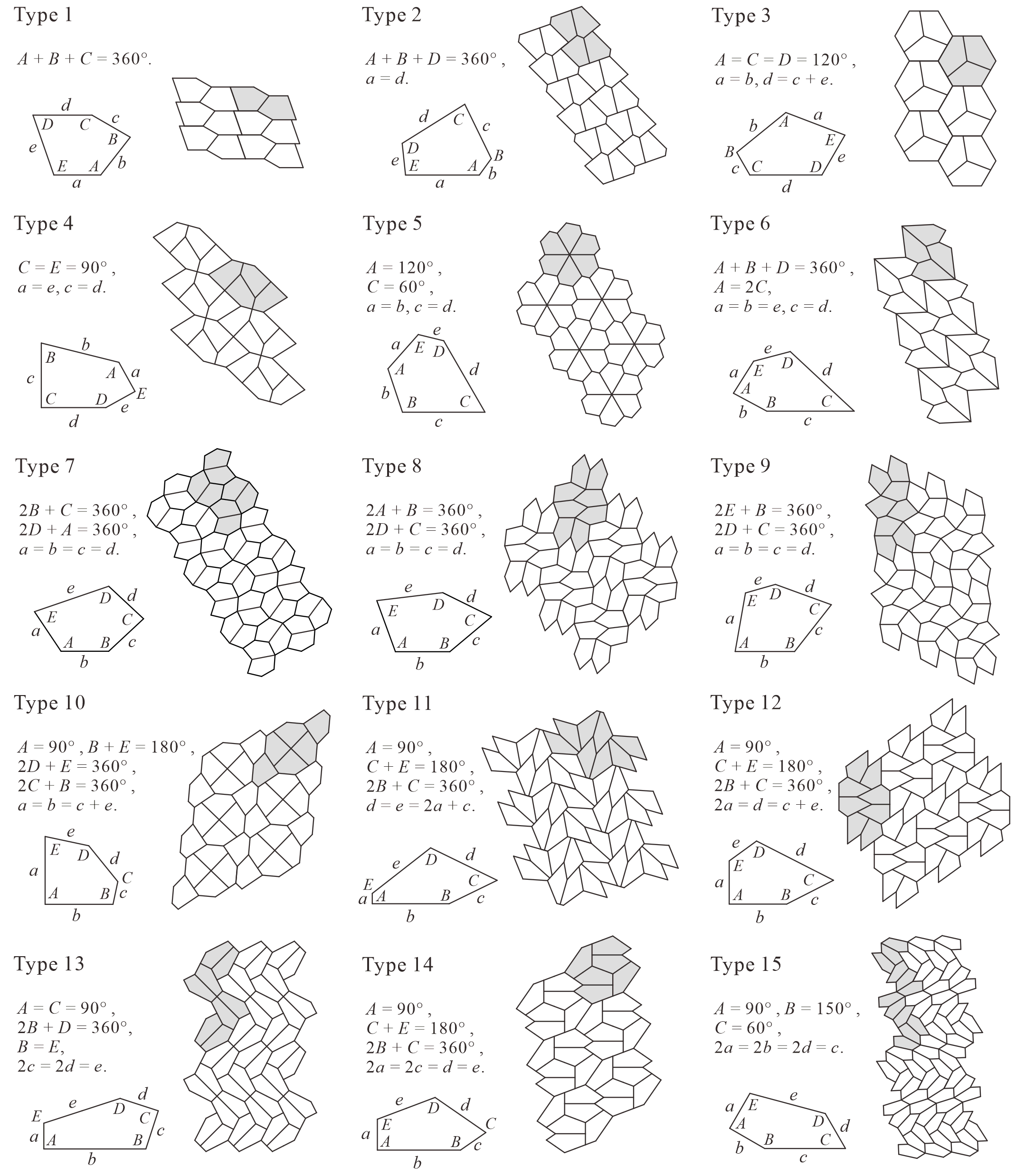} 
\caption{{\small 
Fifteen Type families of convex pentagonal monotiles. Each of the 
convex pentagonal monotiles is defined by some conditions between the 
lengths of the edges and the magnitudes of the angles, but some degrees of freedom 
remain. For example, convex pentagonal monotiles belonging to the Type 1 family 
(also referred to as convex pentagonal monotiles belonging to Type 1, or simply 
as convex pentagonal monotiles of Type 1) satisfy the property the sum of 
three consecutive angles is equal to $360^ \circ$. In this figure, this property is 
expressed as $A + B + C = 360^ \circ $, which represents the tile conditions 
of Type 1 (the Type 1 family). The pentagonal monotiles of Types 14 and 15 
have one degree of freedom, that of size. For example, the value of angle $C$ in the 
pentagonal monotile of Type 14 is $\cos ^{ - 1} \bigl( \bigl( 3\sqrt {57} 
- 17 \big) / 16 \big)  \approx 1.2099 \;\mbox{rad} \approx 69.32^ \circ $. 
The gray region in each tiling indicates a translation unit (a unit that can generate a 
periodic tiling through translation alone).
}
\label{fig23v2}
}
\end{figure}

The Type families of convex pentagonal monotiles are not necessarily 
``disjoint." That is, there are some convex pentagonal monotiles that belong to more 
than one Type family  \cite{Sugimoto_2016, Sugimoto_2020a, Sugimoto_2020b, Sugimoto_2025, 
Sugi_Araki_2017a}. For example, the convex pentagonal monotile in Figure~\ref{fig24v2} 
belongs to both the Type 1 and Type 7 families \cite{Sugimoto_2020a, Sugimoto_2020b}. 
If $A = 86^ \circ $ for the tile conditions of Type 7 in Figure~\ref{fig23v2}, the convex 
pentagonal monotile with $A=86^ \circ$ belongs only to the Type 7 family and can only 
form the (representative) tiling of the Type 7 family, as depicted in Figure~\ref{fig23v2}. By contrast, 
the convex pentagonal monotiles in Figure~\ref{fig24v2}, which belong to both the Type 1 
and Type 7 families, can form the tiling of the Type 7 family, but they can 
also form the tiling of the Type 1 family without using reflected tiles. In other words, 
the convex pentagonal monotiles belonging to both the Type 1 and Type 7 
families in Figure~\ref{fig24v2} can generate tiling without using reflected tiles. 
 It should be noted that some convex pentagonal monotiles can form tilings 
other than the representative tilings of each Type \cite{Sugimoto_2020a, 
Sugimoto_2020b, Sugimoto_2025, Sugi_Araki_2017a}.

Some tiling patterns cannot be formed without using reflected tiles. 
However, the tile forming the pattern does not necessarily have to be 
reflected to generate the tilings. There are also monotiles that cannot 
generate tilings if the use of reflected tiles is not allowed during the tiling 
generation process.

Figure~\ref{fig25v2} shows the periodic tiling formed by convex pentagonal monotiles 
belonging to both the Type 2 and Type 4 families with line symmetry in the cases of 
$\beta = 90^ \circ $ in Method 2 (Figure~\ref{fig05}(d)), $\alpha = 90^ \circ $ 
in Method 1 (Figure~\ref{fig08}), and $\gamma = 165^ \circ $ 
in Method 4 (Figure~\ref{fig18v2}(f)).

\bigskip
\renewcommand{\figurename}{{\small Figure}}
\begin{figure}[htbp]
\centering\includegraphics[width=15cm,clip]{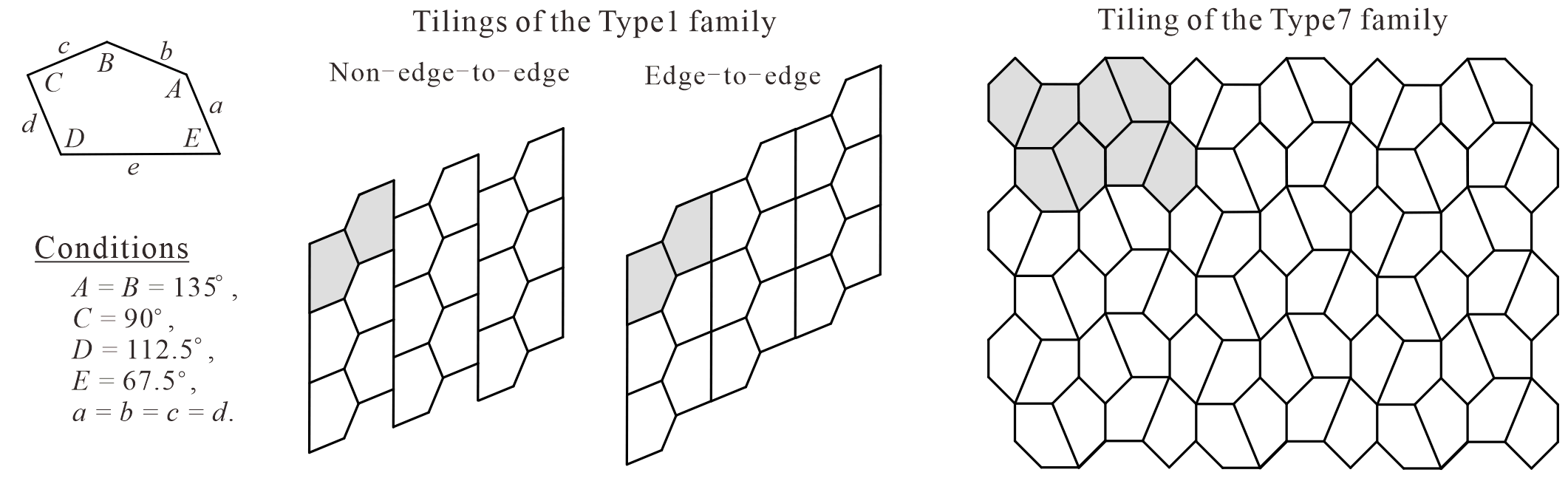} 
\caption{{\small 
Convex pentagonal monotile belonging to both the Type 1 and Type 7 families, 
the tilings of the Type 1 family and the tiling of the Type 7 family with 
the convex pentagon.
}
\label{fig24v2}
}
\end{figure}

\bigskip
\renewcommand{\figurename}{{\small Figure}}
\begin{figure}[htbp]
\centering\includegraphics[width=13.5cm,clip]{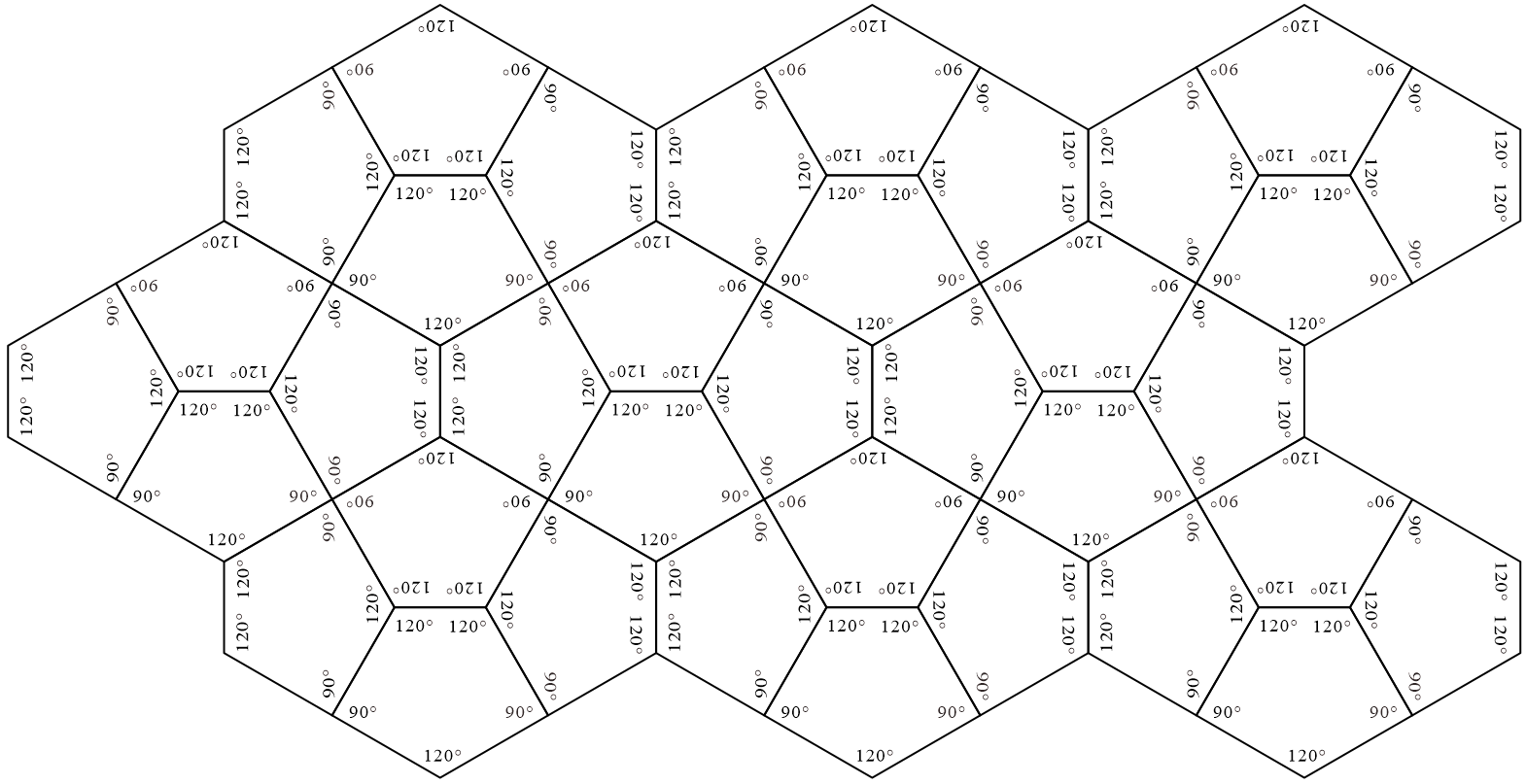} 
\caption{{\small 
Periodic tiling with convex pentagonal monotile belonging to both 
the Type 2 and Type 4 families.
}
\label{fig25v2}
}
\end{figure}

%%%%%%%%%%%%%%%%%%%%%%%%%%%%%%%%%%%%%%%%%%%%%%%%%%%%%%%%%%%%%%%%%%%%%%
%%%%%%%%%%%%%%%%%%%%%%%%%%%%%%%%%%%%%%%%%%%%%%%%%%%%%%%%%%%%%%%%%%%%%%

\section{Convex hexagonal monotiles}
\label{appB}

Three families of convex hexagonal monotiles are known, called ``Type," as 
shown in Figure~\ref{fig26v2} \cite{Bollo_1963, Gardner_1975, G_and_S_1987, Kershner_1968, 
Reinhardt_1918, Sugimoto_2017, Zong_2020}. The tilings depicted with each convex 
hexagon in Figure~\ref{fig26v2} are representative tilings of each Type family. Representative 
tiling can be formed using only those relations that can be derived from the 
tile conditions of each Type. As shown in Figure~\ref{fig26v2}, the representative 
tiling of the Type 2 family uses reflected tiles (convex hexagons). If convex 
hexagonal monotiles do not belong to the Type 1 or Type 3 families, then 
the convex hexagonal monotiles cannot generate tilings without using 
reflected tiles. The Type families of convex hexagonal monotiles are not 
``disjoint." For example, a regular hexagon belongs to all Type 
families.

Figures~\ref{fig27v2}(a) and \ref{fig27v2}(b) show the periodic tilings formed by convex 
hexagonal monotiles belonging to the Type 1 family in the cases of $\alpha = 75^ \circ $ 
(Figure~\ref{fig04} (d)) and $\alpha = 105^ \circ $ (Figure~\ref{fig04}(e)) in Method 1, 
respectively.

\bigskip
\renewcommand{\figurename}{{\small Figure}}
\begin{figure}[htbp]
\centering\includegraphics[width=14cm,clip]{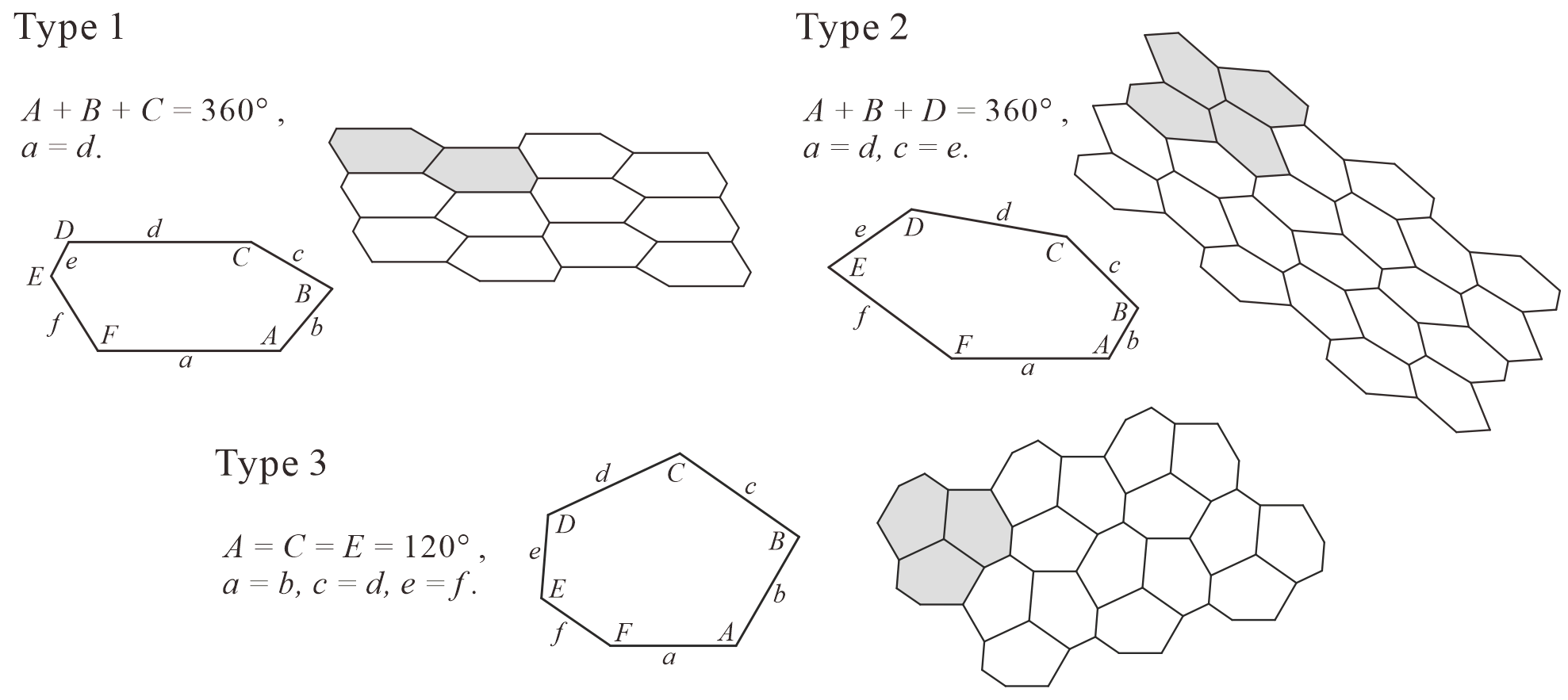} 
\caption{{\small 
Three Type families of convex hexagonal monotiles. The gray 
region in each tiling indicates a translation unit (a unit that can generate 
a periodic tiling by translation only).
}
\label{fig26v2}
}
\end{figure}

\renewcommand{\figurename}{{\small Figure}}
\begin{figure}[htbp]
\centering\includegraphics[width=14cm,clip]{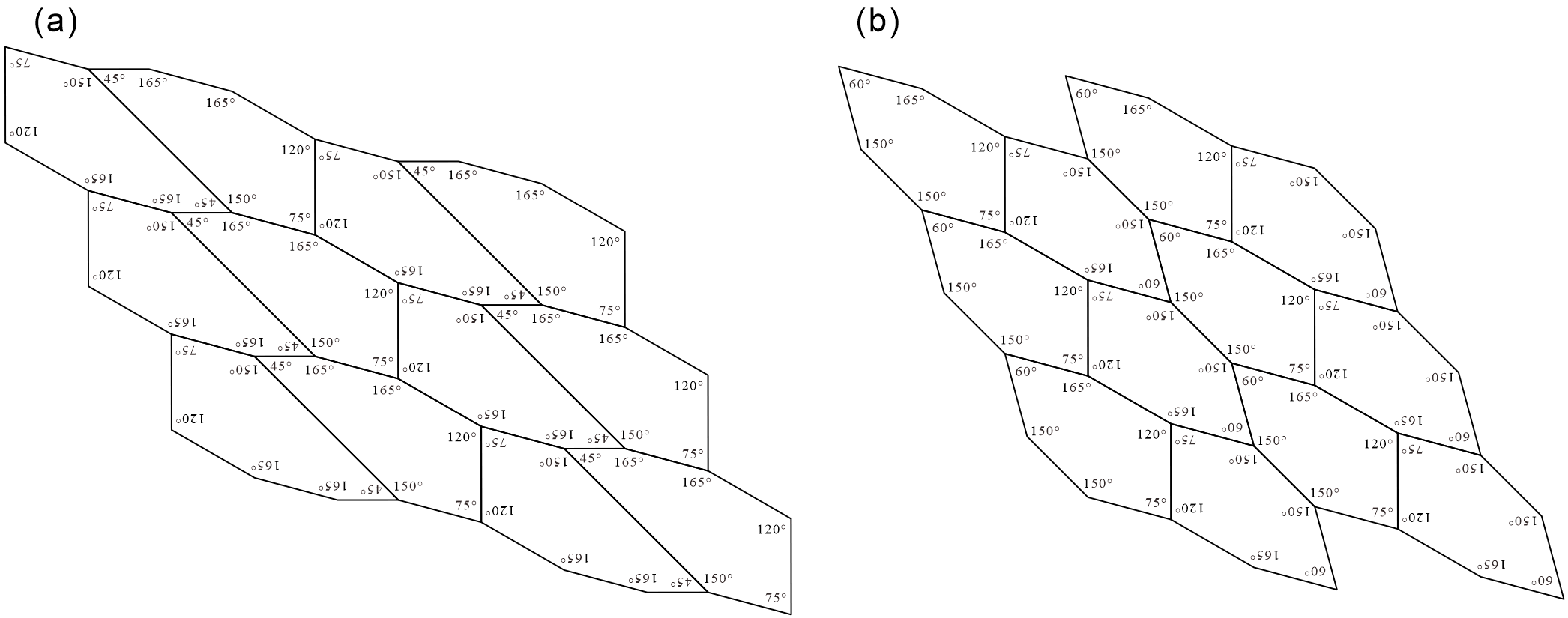} 
\caption{{\small 
Periodic tilings with convex hexagonal monotiles belonging to the 
Type 1 family, created by the division of Method 1.
}
\label{fig27v2}
}
\end{figure}

%%%%%%%%%%%%%%%%%%%%%%%%%%%%%%%%%%%%%%%%%%%%%%%%%%%%%%%%%%%%%%%%%%%%%%
%%%%%%%%%%%%%%%%%%%%%%%%%%%%%%%%%%%%%%%%%%%%%%%%%%%%%%%%%%%%%%%%%%%%%%

\section{Periodic tiling with $P_{1}$, $P_{2}$, and reflected $P_{2}$}
\label{appC}

As introduced in Section~\ref{section2}, if $P_{1}$, $P_{2}$, and reflected $P_{2}$ in the 
case of $\alpha = 180^ \circ $ in Method 1 (see Figure~\ref{fig04}(g)) are used, it is 
possible to form a periodic tiling, as shown in Figure~\ref{fig28v2}. In this case, the 
merged diagram of $P_{2}$ and reflected $P_{2}$ can be regarded as forming a 
convex octagon. In other words, Figure~\ref{fig28v2} can be regarded as periodic tiling 
using anterior and posterior convex octagons and convex pentagons without 
any distinction between the anterior and posterior sides (two types of 
convex pentagons and octagons that are not monotiles from a viewpoint that 
does not distinguish between the anterior and posterior sides).

\renewcommand{\figurename}{{\small Figure}}
\begin{figure}[htbp]
\centering\includegraphics[width=14cm,clip]{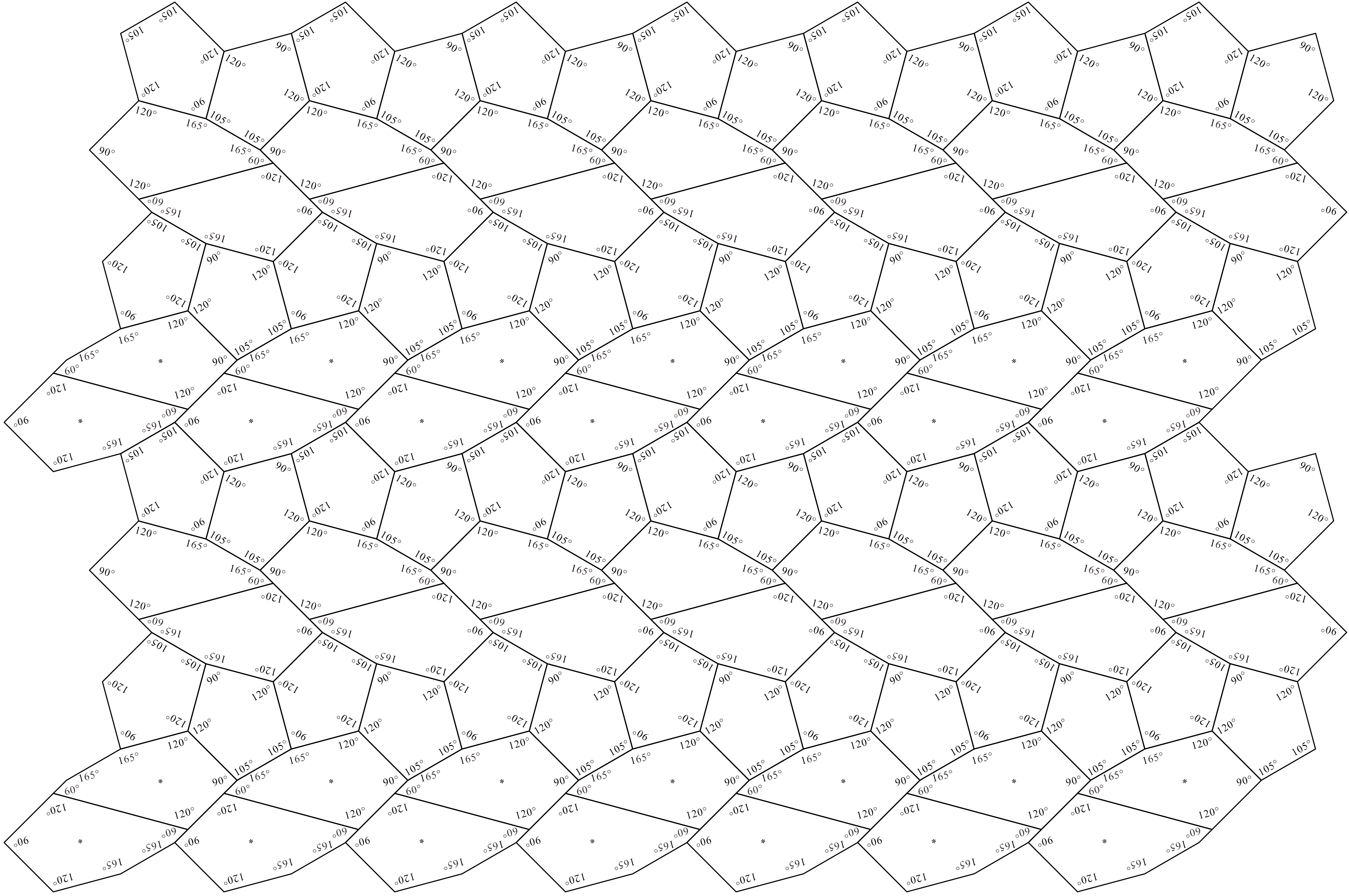} 
\caption{{\small 
Periodic tiling with $P_{1}$, $P_{2}$, and reflected $P_{2}$ in the case of $\alpha = 180^ \circ $ 
in Method 1 (see Figure~\ref{fig04}(g)). The reflected $P_{2}$ are given an asterisk mark ``*.''
}
\label{fig28v2}
}
\end{figure}

%%%%%%%%%%%%%%%%%%%%%%%%%%%%%%%%%%%%%%%%%%%%%%%%%%%%%%%%%%%%%%%%%%%%%%
%%%%%%%%%%%%%%%%%%%%%%%%%%%%%%%%%%%%%%%%%%%%%%%%%%%%%%%%%%%%%%%%%%%%%%

\section{ASPs comprising four or five types of convex polygons\protect\footnote{
The contents (the results of study) in this appendix were obtained prior to the 
main results of this manuscript.}}
\label{appD}

\subsection{Preparation}
\label{appD_1}

In this appendix (i.e., \ref{appD}), we present examples of ``ASP comprising four 
types of convex polygons" and ``(candidates of) ASP comprising five types of convex 
polygons." Let $ASPcvx(n)$ denote an ASP comprising $n$ types of convex 
polygons. However, we limit ourselves to $ASPcvx(4)$ and $ASPcvx(5)$ that 
satisfy the ``restriction on the number of types of congruent diagrams" below 
and show examples of them.

\begin{nameth}
For a viewpoint that does not distinguish between the anterior and posterior sides 
\rm{(}\it{a viewpoint that treats congruent diagrams, including reflected images, as 
tiles of the same type}\rm{)}\it{, the number of prototile types in a set must be less 
than or equal to three.}
\end{nameth}

In this study, the ``prototiles" mentioned in the above restriction 
corresponds to convex polygons; let $ASPcvx(n, m)$ be ASP comprising $n$ 
types of convex polygons such that there are $m$ types of convex polygons from 
a viewpoint that does not distinguish between anterior and posterior sides 
($m$ is an integer greater than or equal to one). Thus, it is always $n \ge m.$ 
Let $ASPcvx(n, m \le 3)$ denote $ASPcvx(n)$ satisfying the restriction on the 
number of types of congruent diagrams. 
Because these ``ASP comprising four types of convex polygons" and ``candidate 
of ASP comprising five types of convex polygons" satisfy the restriction on the 
number of types of congruent diagrams, they can be regarded as ``ASP comprising 
three types of convex polygons" from the viewpoint that does not distinguish 
between anterior and posterior sides.

The restriction on the number of types of congruent diagrams is imposed because 
the method of creating the ASP introduced in the study of this appendix would create 
a variety of $ASPcvx(n)$ with $n \ge 4$ ($ASPcvx(n \ge 4, m \ge 4)$) that do 
not satisfy this restriction. In addition, considering that only a small number of 
types of prototype tiles contained in an ASP has been explored so far, we 
thought that it would be meaningful to consider cases that satisfy this restriction.

%%%%%%%%%%%%%%%%%%%%%%%%%%%%%%%%%%%%%%%%%%%%%%%%%%%%%%%%%%%%%%%%%%%%%%

\subsection{Dividing \mbox{Tile$(1, 1)$} into convex pentagons and convex hexagons}
\label{appD_2}

\mbox{Tile$(1, 1)$}, corresponding to $a=b=1$ in \mbox{Tile$(a, b)$} shown in \cite{Smith_2024a}, 
can generate a periodic tiling if the use of reflected tiles is allowed during the tiling 
generation process (see Figure~\ref{fig02}). However, as shown in \cite{Smith_2024b}, 
\mbox{Tile$(1, 1)$} can only generate non-periodic tiling if and only if it does not allow the use of reflected 
tiles (see Figure~\ref{fig03}). Using this property, $ASPcvx(n)$ can be obtained from \mbox{Tile$(1, 1)$}.

This method divides the interior of \mbox{Tile$(1, 1)$} into convex polygons with at 
least five edges. However, the purpose of the study of this appendix was to conduct 
$ASPcvx(4, m \le 3)$ and \mbox{$ASPcvx(5, m \le 3)$}. The method for dividing 
\mbox{Tile$(1, 1)$} such that the restriction on the number of types of congruent 
diagrams is satisfied is to divide the interior into five convex polygons, 
as shown in Figure~\ref{fig03AD}. (Refer to the \ref{appD_5} and \ref{appD_6} for 
other methods of dividing \mbox{Tile$(1, 1)$} that satisfy the restriction on the 
number of types of congruent diagrams that we are currently identifying.) It is inferred 
that there are various methods to divide \mbox{Tile$(1, 1)$} that are likely to result 
in ASP that do not satisfy the restriction on the number of types of congruent diagrams 
(see \ref{appD_5}).

Figure~\ref{fig03AD} shows how to divide the interior into five convex polygons so that 
the parts with vertices ``$X_{3}, X_{7}, X_{9}, X_{11}$" at $90^ \circ $ in \mbox{Tile$(1, 1)$} 
are convex pentagons or convex hexagons. Therefore, a line segment 
$X_{8}R$ is added to bisect interior angle $X_{8}$ of \mbox{Tile$(1, 1)$}, and its length 
is the same as that of its edge when \mbox{Tile$(1, 1)$} is considered an equilateral 
14-edges polygon. Subsequently, to satisfy the restriction on the number of 
types of congruent diagrams, the convex polygon $CP(X_{3})$ with the points 
``$X_{3}$, $X_{4}$, $X_{5}$, $P$, $X_{1}$, $X_{2}$" as vertices and the convex polygon 
$CP(X_{7})$ with the points ``$X_{7}$, $X_{8}$, $R$, $P$, $X_{5}$, $X_{6}$" as vertices are 
congruent convex polygons with reflection relation. In addition, the convex 
polygon $CP(X_{9})$ with the points ``$X_{9}$, $X_{10}$, $T$, $R$, $X_{8}$" as vertices and 
the convex polygon $CP(X_{11})$ with the points ``$X_{11}$, $X_{12}$, $X_{13}$, $T$, 
$X_{10}$" are congruent convex pentagons with reflection relation.

As shown in Figure~\ref{fig03AD}, by changing the parameters $\delta$ and $\epsilon$, the 
shape of the five convex polygons inside \mbox{Tile$(1, 1)$} changes as point $P$ moves 
on the bisector $X_{5}Q$ of the interior angle $X_{5}$ and point $T $ moves on the 
bisector $X_{10}S$ of the interior angle $X_{10}$ (the edge corresponding to 
the green line in Figure~\ref{fig03AD} has its length changed by changing $\delta $ and 
$\epsilon$). Based on the change, the following cases appear inside \mbox{Tile$(1, 1)$} 
from a viewpoint that does not distinguish between anterior and 
posterior sides: ``one type of convex pentagon and two types of convex 
hexagons," ``two types of convex pentagons and one type of convex hexagon," 
and ``three types of convex pentagons."

The convex polygon $CP(X_{14})$ with the points ``$X_{14}$, $X_{1}$, $P$, $R$, $T$, 
$X_{13}$" as vertices in Figure~\ref{fig03AD} must be $105^ \circ \le \delta \le 165^ \circ $, 
$65.104^ \circ < \epsilon \le 120^ \circ $, and $\delta + \epsilon \ge 180^ \circ $ to 
be convex\footnote{ 
When $\epsilon = (\pi / 3) + \tan ^{ - 1}\left( {{\left( {2 - \sqrt 3 } \right)} 
\mathord{\left/ {\vphantom {{\left( {2 - \sqrt 3 } \right)} 3}} \right. 
\kern-\nulldelimiterspace} 3} \right) \approx 65.104^\circ $ (i.e., $T=X_{10})$ 
in Figure~\ref{fig03AD}, $CP(X_{9})$ and $CP(X_{11})$ are convex 
quadrilaterals, therefore $65.104^ \circ < \epsilon \le 120^ \circ $.}. 
Moreover, $CP(X_{14})$ is a convex pentagon in the following three cases: 
``$\delta + \epsilon = 180^ \circ $ (Figure~\ref{fig04AD}(a))," 
``$\epsilon = 120^ \circ $ ($T=S)$ and $105^ \circ \le \delta < 165^ \circ $ (Figure~\ref{fig04AD}(b))," 
and ``$\delta = 165^ \circ $ and $105^ \circ \le \epsilon < 120^ \circ $ (Figure~\ref{fig04AD}(c))."
If $\delta = 105^ \circ $, $P=X_{5}$ holds, $CP(X_{3})$ and $CP(X_{7})$ are convex pentagons 
(otherwise, convex hexagons). Note that during $\delta = 165^ \circ $ and 
$\epsilon = 120^ \circ $, $CP(X_{14})$ is a convex quadrangle; therefore, it is 
beyond the scope of the study of this appendix.

Let $S_t (\delta, \epsilon )$ be the tile set of {\{}$CP(X_{14})$, $CP(X_{3})$, 
$CP(X_{7})$, $CP(X_{9})$, $CP(X_{11})${\}} for parameters $\delta$ and 
$\epsilon$. The convex polygons contained in $S_t (\delta, \epsilon )$ become 
convex pentagonal or hexagonal monotiles if $\delta$ and $\epsilon$ satisfy 
certain conditions, as listed in Table~\ref{tab1} \cite{G_and_S_1987, Sugimoto_2016, 
Sugimoto_2017, wiki_pentagon_tiling}. 
All cases in in Table~\ref{tab1} contain the convex pentagonal or hexagonal monotiles 
that can form periodic tilings without using reflected tiles (see \ref{appA} and \ref{appB})\footnote{
The case of $\epsilon = 120^ \circ $ and $105^ \circ \le \delta < 165^ \circ $ also contains 
the convex pentagonal monotiles belonging only to the Type 2 family which require the use 
of reflected tiles during the tiling generation process. However, the tile set in the case contains 
the anterior and posterior sides of convex pentagonal monotiles belonging to the Type 2 family. 
Thus, they can form the representative tiling of the Type 2 family.
}.

\renewcommand{\figurename}{{\small Figure}}
\begin{figure}[htbp]
\centering\includegraphics[width=12cm,clip]{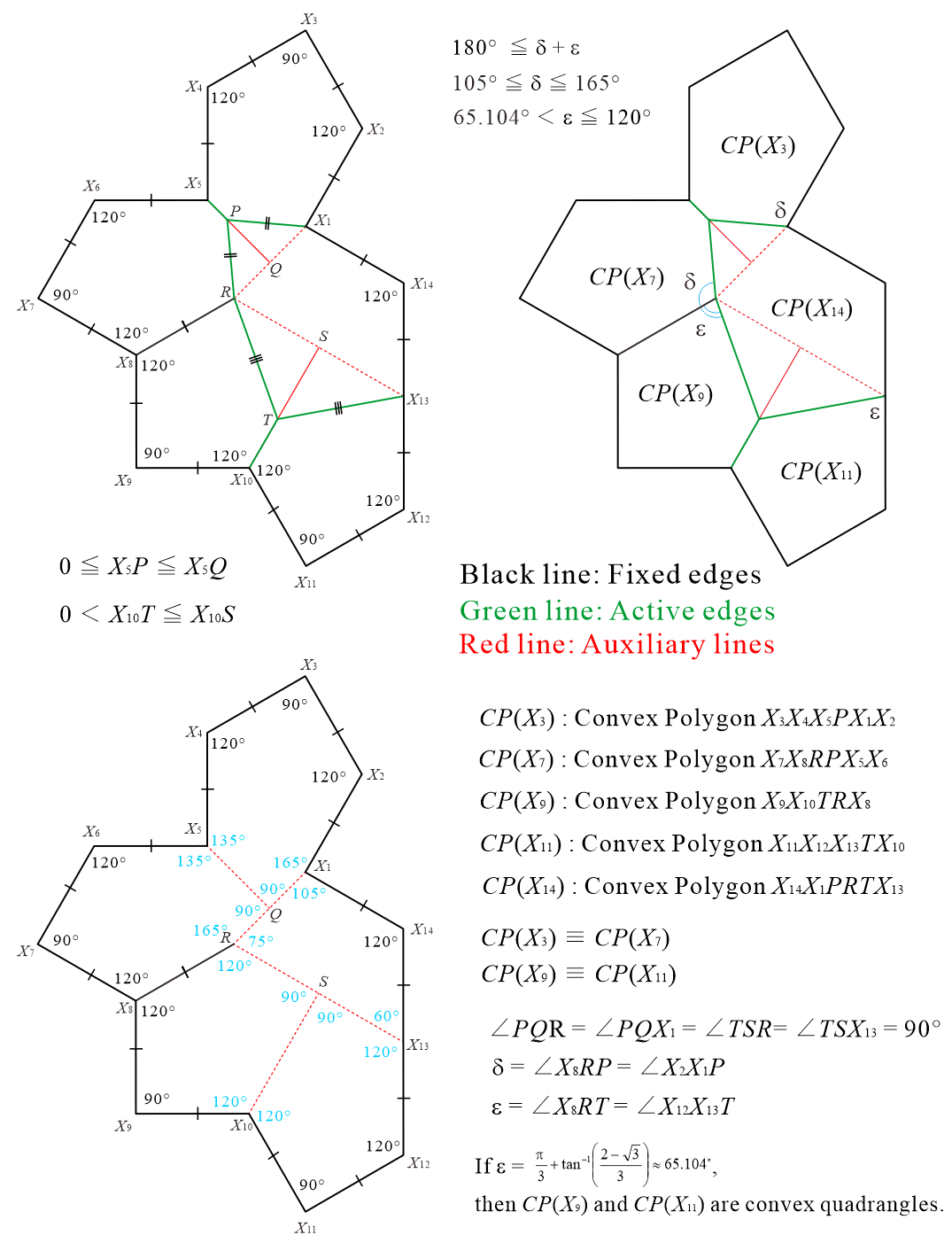} 
\caption{{\small 
The method to divide \mbox{Tile$(1, 1)$} such that the restriction on the 
number of types of congruent diagrams is satisfied.
}
\label{fig03AD}
}
\end{figure}

{\color{white}. }
\renewcommand{\figurename}{{\small Figure}}
\begin{figure}[htbp]
\centering\includegraphics[width=12cm,clip]{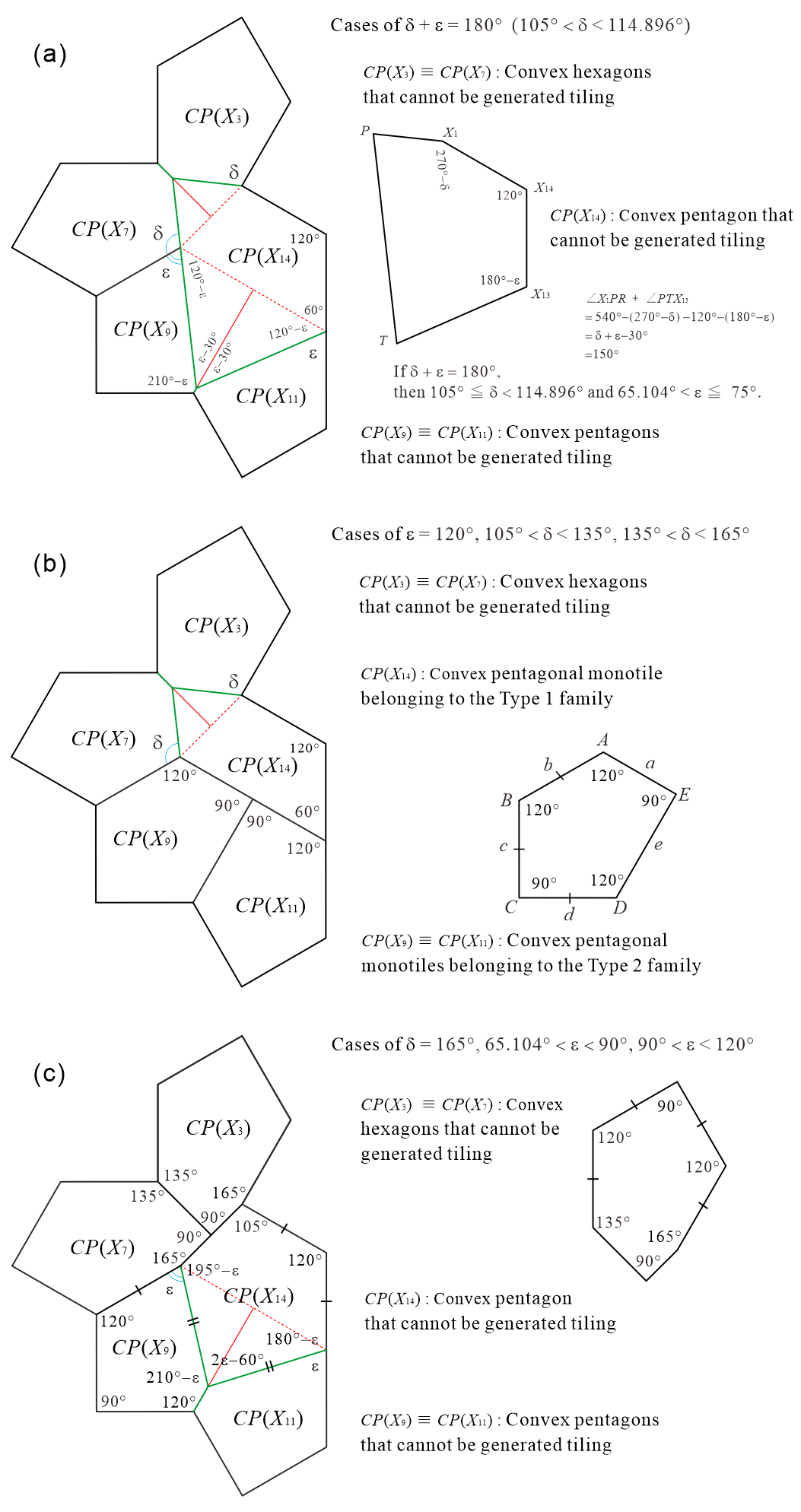} 
\caption{{\small 
Cases of $\delta + \epsilon = 180^ \circ $, $\epsilon = 120^ \circ $ or 
$\delta = 165^ \circ $.
}
\label{fig04AD}
}
\end{figure}

\renewcommand{\figurename}{{\small Figure}}
\begin{figure}[htbp]
\centering\includegraphics[width=12cm,clip]{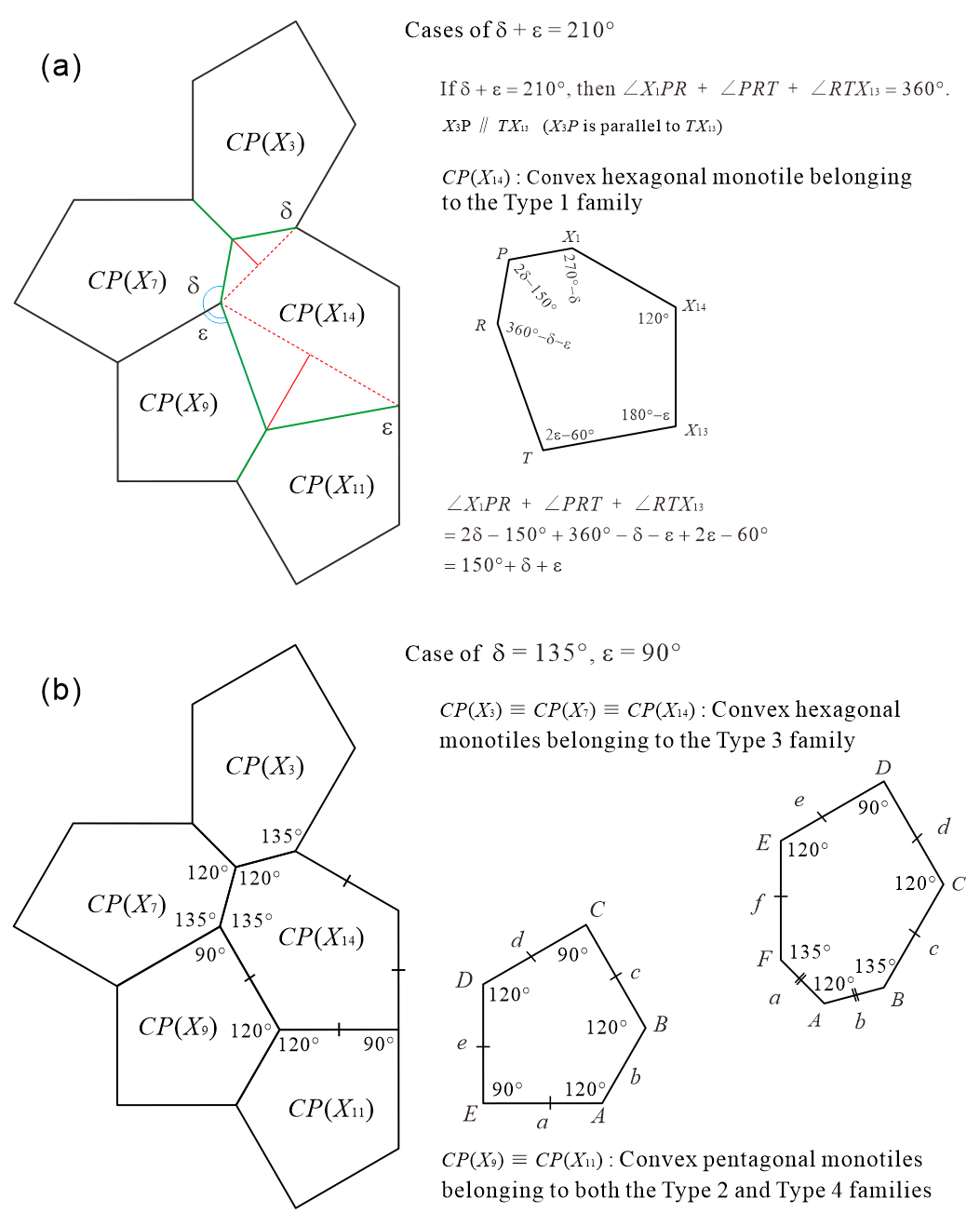} 
\caption{{\small 
Cases of $\delta + \epsilon = 210^ \circ $, $\delta = 135^ \circ $ or 
$\epsilon = 90^ \circ $.
}
\label{fig05AD}
}
\end{figure}

\begin{table}[htbp]
{\small
\caption[Table 1]{{\small Cases where $CP(X_{14})$, $CP(X_{3})$, $CP(X_{7})$, $CP(X_{9})$, and 
$CP(X_{11})$ are monotiles.}}
\label{tab1}
}
\begin{tabular}
{p{68pt}p{89pt}p{89pt}p{90pt}p{46pt}}
\hline
\footnotesize{}& 
\footnotesize{$CP(X_{14})$}& 
\footnotesize{$CP(X_{3})$ and $CP(X_{7})$}& 
\footnotesize{$CP(X_{9})$ and $CP(X_{11})$}& 
\footnotesize{Examples} \\
\hline
\raisebox{-2.5ex}[0cm][0cm]{\footnotesize{$\delta + \epsilon = 210^ \circ $} }& 
\footnotesize{Convex hexagonal \par monotile belonging \par to the Type 1 family}& 
\footnotesize{}& 
\footnotesize{}& 
\raisebox{-2.5ex}[0cm][0cm]{\footnotesize{Figure~\ref{fig05AD}(a)}} \\
\hline
\raisebox{-2.5ex}[0cm][0cm]{\footnotesize{$\delta = 135^ \circ $}}& 
\footnotesize{}& 
\footnotesize{Convex hexagonal \par monotiles belonging \par to the Type 3 family}& 
\footnotesize{}& 
\raisebox{-2.5ex}[0cm][0cm]{\footnotesize{Figure~\ref{fig05AD}(b) }}\\
\hline
\raisebox{-4ex}[0cm][0cm]{\footnotesize{$\epsilon = 90^ \circ $}}& 
\footnotesize{}& 
\footnotesize{}& 
\footnotesize{Convex pentagonal \par monotiles belonging \par to both the Type 2 \par and Type 4 families}& 
\raisebox{-3ex}[0cm][0cm]{\footnotesize{Figure~\ref{fig05AD}(b) }}\\
\hline
\raisebox{-5ex}[0cm][0cm]{\footnotesize{\shortstack[l]{$\delta = 135^ \circ $ and \\ 
$\epsilon = 90^ \circ $}}}& 
\raisebox{-6ex}[0cm][0cm]{\footnotesize{\shortstack[l]{Convex hexagonal \\ monotile belonging 
\\ to the Type 3 family}}}& 
\raisebox{-6ex}[0cm][0cm]{\footnotesize{\shortstack[l]{Convex hexagonal \\ monotiles belonging 
\\ to the Type 3 family}}}& 
\footnotesize{Convex pentagonal \par monotiles belonging \par to both the Type 2 \par and Type 4 families }& 
\raisebox{-3ex}[0cm][0cm]{\footnotesize{Figure~\ref{fig05AD}(b) } }\\
\hline
\raisebox{-4ex}[0cm][0cm]{\footnotesize{\shortstack[l]{$\epsilon = 120^ \circ $ and \\
$105^ \circ \le \delta < 165^ \circ $}}}& 
\footnotesize{Convex pentagonal \par monotile belonging \par to the Type 1 family}& 
\footnotesize{}& 
\footnotesize{Convex pentagonal \par monotiles belonging \par to the Type 2 family}& 
\raisebox{-2.5ex}[0cm][0cm]{\footnotesize{Figure~\ref{fig04AD}(b)}}\\
\hline
\end{tabular}

\end{table}
\bigskip

%%%%%%%%%%%%%%%%%%%%%%%%%%%%%%%%%%%%%%%%%%%%%%%%%%%%%%%%%%%%%%%%%%%%%%

\subsection{ASPs satisfying the restriction on the number of types of congruent 
diagrams created from \mbox{Tile$(1, 1)$}}
\label{appD_3}

If the convex polygons contained in $S_t (\delta, \epsilon )$ are not convex 
pentagonal or hexagonal monotiles, we can create candidates of $ASPcvx(4, 3)$ 
and $ASPcvx(5,3)$ from \mbox{Tile$(1, 1)$}. We refer to it as ``candidate" because if 
we generate tilings with convex polygons contained in $S_t (\delta, \epsilon )$, 
we have to make sure that there exist only non-periodic tilings using 
the shape of \mbox{Tile$(1, 1)$}, if and only if all types of convex polygons 
contained in $S_t (\delta, \epsilon)$ are used. Therefore, we will confirm it 
with the following steps:

\begin{itemize}
\item \textbf{Step 1}. Check that the convex polygons contained in 
$S_t (\delta, \epsilon)$ are not monotile.

\item \textbf{Step 2}. If it uses more than one type (but not all) of 
the convex polygons contained in $S_t (\delta, \epsilon)$, it cannot admit tilings.

\item \textbf{Step 3}. Check that with all the convex polygons contained in 
$S_t (\delta, \epsilon)$, they cannot admit tilings in any combination other than 
the shape of \mbox{Tile$(1, 1)$}, as shown in Figure~\ref{fig03AD}.
\end{itemize}

For example, we confirmed that $S_t (\delta = 105^ \circ ,\epsilon = 75^ \circ )$, 
which contains five convex pentagons created by dividing \mbox{Tile$(1, 1)$} shown 
in Figure~\ref{fig06AD}, satisfies Steps 1--3. As shown in Figure~\ref{fig06AD}, when $\delta = 105^ \circ $ 
and $\epsilon = 75^ \circ $ are used, $CP(X_{3})$ and $CP(X_{7})$ are convex 
pentagons with line symmetry and cannot distinguish between the anterior and 
posterior sides. Therefore, $S_t (\delta = 105^ \circ , \epsilon = 75^ \circ )$ 
is $ASPcvx(4, 3)$. In $S_t (\delta, \epsilon )$, we conjecture that Figure~\ref{fig06AD} is 
the only case where the ASP comprises of four types of convex pentagons, and 
the ASP comprises of three types of convex pentagons from a viewpoint that 
does not distinguish between the anterior and posterior sides.

In Steps 1--3, we briefly describe some of the key points for identifying 
cases where the convex polygons in Figure~\ref{fig06AD} cannot form tilings. We focused 
on the interior angles of the convex polygons at $75^ \circ$, $105^ \circ$, 
$135^ \circ$, and $165^ \circ$. To form tilings, two or more tile vertices must be 
concentrated at one point, and the sum of the interior angles of the 
vertices at the concentration point must be $180^ \circ$ or $360^ \circ$. Therefore, 
as there are no convex pentagons with interior angles of $60^ \circ$ or less, 
$135^ \circ$ and $165^ \circ$ cannot be used in a concentration with a sum of 
$180^ \circ$. We also observe that an even number of vertices with interior angles 
of $75^ \circ$, $105^ \circ$, $135^ \circ$, and $165^ \circ$ must always be concentrated 
at one point because the place value of one is ``five." Based on these considerations, it was 
concluded that the concentration of vertices with an interior angle of $135^ \circ$ was 
impossible, except for the arrangements using $CP(X_{9})$ and $CP(X_{11})$, as shown in 
Figure~\ref{fig07AD}(a). Subsequently, note the point in Figure~\ref{fig07AD}(a) where the two 
vertices have an internal angle of $135^ \circ$. If a tiling is formed, the only method to 
obtain a sum of $360^ \circ$ is to place the vertex of the tile with an interior angle 
of $90^ \circ$ at the concentrating point. Figure~\ref{fig07AD}(b) shows the only 
combination that can form a structure in which the tiling does not collapse. 
Then, for the positions with vertices having internal angles of $75^ \circ$ and 
$165^ \circ$ circled in blue in the cluster in Figure~\ref{fig07AD}(b), the only combination 
that can construct a structure in which the tiling does not collapse using 
$S_t (\delta = 105^ \circ , \epsilon = 75^ \circ )$ is that of dividing 
\mbox{Tile$(1, 1)$} (using the shape of \mbox{Tile$(1, 1)$}), as shown in Figure~\ref{fig06AD}.

Figure~\ref{fig08AD} shows non-periodic tiling generated by an ASP comprising the four types 
of convex pentagons shown in Figure~\ref{fig06AD}. The tiling is non-edge-to-edge.

\bigskip
\bigskip

\renewcommand{\figurename}{{\small Figure}}
\begin{figure}[htbp]
\centering\includegraphics[width=14.5cm,clip]{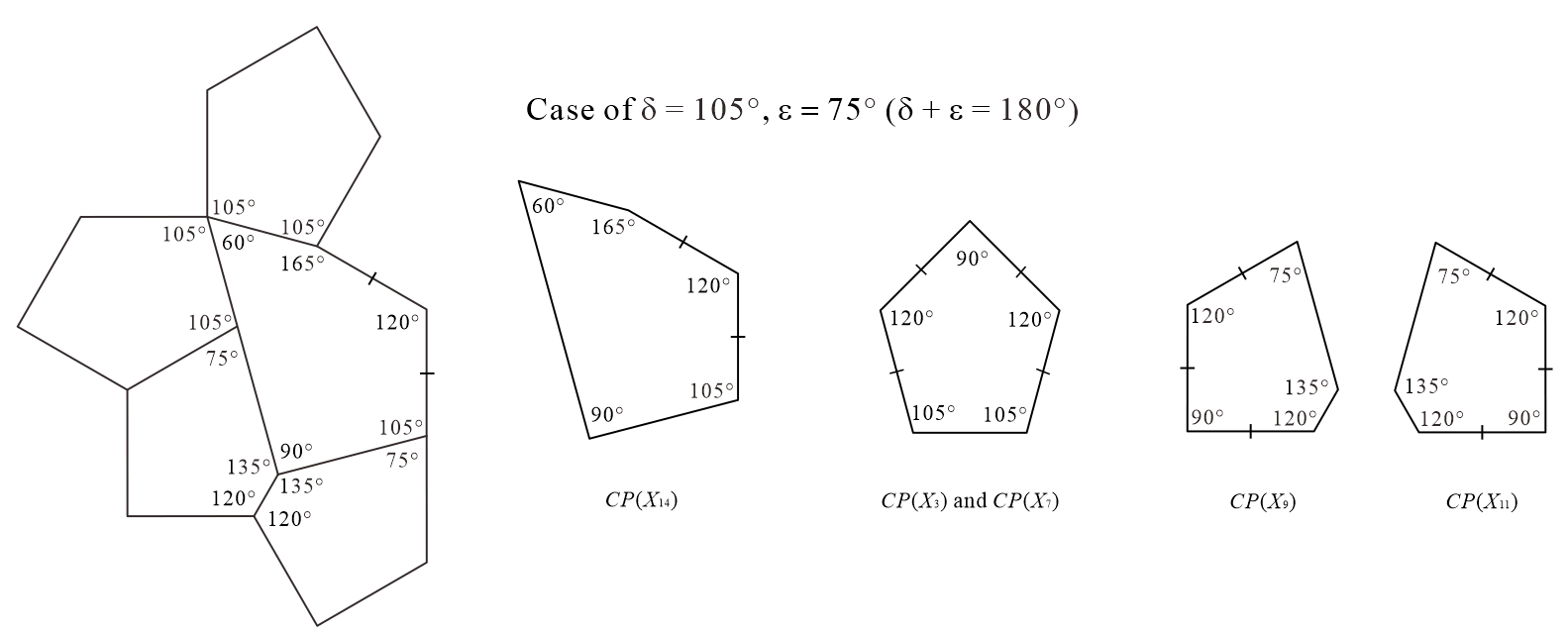} 
\caption{{\small 
How to divide \mbox{Tile$(1, 1)$} to create an ASP comprising four types 
of convex pentagons.
}
\label{fig06AD}
}
\end{figure}

{\color{white}. }

\bigskip
\bigskip
\renewcommand{\figurename}{{\small Figure}}
\begin{figure}[htbp]
\centering\includegraphics[width=8.5cm,clip]{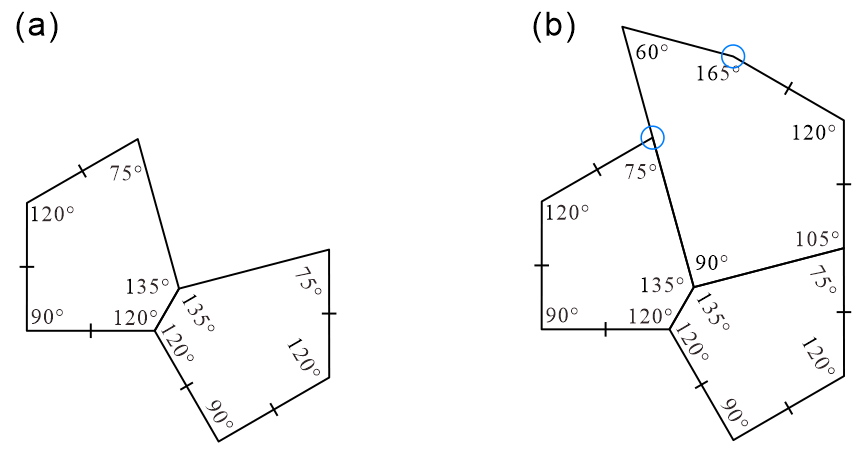} 
\caption{{\small 
Example of combination focusing on vertices with an internal angle 
of $135^ \circ$.
}
\label{fig07AD}
}
\end{figure}

\renewcommand{\figurename}{{\small Figure}}
\begin{figure}[htbp]
\centering\includegraphics[width=15cm,clip]{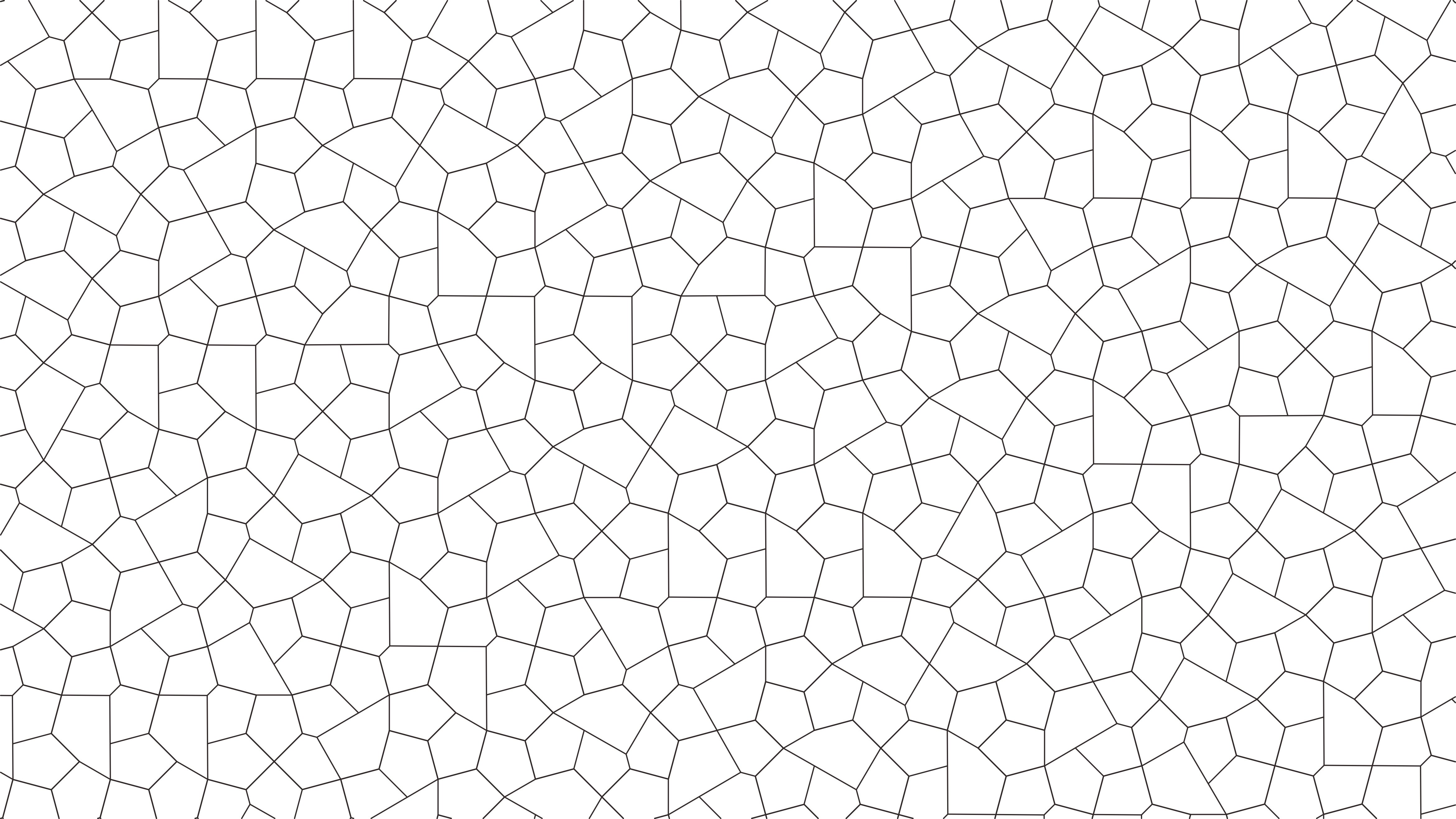} 
\caption{{\small 
Non-periodic tiling generated by the ASP comprising four types of convex 
pentagons.
}
\label{fig08AD}
}
\end{figure}

\bigskip
\bigskip
\bigskip
\bigskip
\bigskip
\bigskip
\bigskip
\bigskip
Other candidates for $ASPcvx(4,3)$ and $ASPcvx(5,3)$ related to 
$S_t (\delta, \epsilon )$ are as follows:

\begin{itemize}
\item Candidates for $ASPcvx(4, 3)$
\begin{itemize}
\item[$\circ$] \textbf{Case 1}. When $\delta = 105^ \circ $, $75^ \circ < \epsilon < 120^ \circ $, 
$\epsilon \ne 90^ \circ $, and $\delta + \epsilon \ne 210^ \circ $
(Figures~\ref{fig03AD} and \ref{fig06AD}), 
there is a ``set of three types of convex pentagons ($CP(X_{3})$, $CP(X_{7})$, $CP(X_{9})$, 
$CP(X_{11})$) and one type of convex hexagon ($CP(X_{14})$)," and, from the 
viewpoint that does not distinguish between anterior and posterior 
sides, it is a ``set of two types of convex pentagons and one type of convex 
hexagon."
\end{itemize}
\end{itemize}

\begin{itemize}
\item Candidates for $ASPcvx(5,3)$
\begin{itemize}
\item[$\circ$] \textbf{Case 2}. When $\delta + \epsilon = 180^ \circ $ and $\delta \ne 105^ \circ $ 
(Figure~\ref{fig04AD}(a)), there is a ``set of three types of convex pentagons 
($CP(X_{14})$, $CP(X_{9})$, $CP(X_{11})$) and two types of convex hexagons 
($CP(X_{3})$, $CP(X_{7})$),'' and, from the viewpoint that does not distinguish 
between anterior and posterior sides, it is a ``set of two types of 
convex pentagons and one type of convex hexagon.''

\item[$\circ$] \textbf{Case 3}. When $\delta = 165^ \circ $, $65.104^ \circ < \epsilon < 120^ \circ $, 
and $\epsilon \ne 90^ \circ $ (Figure~\ref{fig04AD}(c)), there is a ``set of three 
types of convex pentagons ($CP(X_{14})$, $CP(X_{9})$, $CP(X_{11})$) and two types of 
convex hexagons ($CP(X_{3})$, $CP(X_{7})$),'' and, from the viewpoint that does 
not distinguish between anterior and posterior sides, it is ``set of two 
types of convex pentagons and one type of convex hexagon.''

\item[$\circ$] \textbf{Case 4}. When $105^ \circ < \delta < 165^ \circ $, 
$\delta \ne 135^ \circ $, $65.104^ \circ < \epsilon < 120^ \circ $, $\epsilon \ne 90^ \circ $, 
$\delta + \epsilon > 180^ \circ $, and $\delta + \epsilon \ne 210^ \circ $ 
(Figure~\ref{fig03AD}), there is a ``set of two types of convex pentagons 
($CP(X_{9})$, $CP(X_{11})$) and three types of convex hexagons 
($CP(X_{14})$, $CP(X_{3})$, $CP(X_{7})$),'' and, from the viewpoint that does not distinguish 
between anterior and posterior sides, it is a ``set of one type of 
convex pentagon and two types of convex hexagon.''
\end{itemize}
\end{itemize}

As can be observed from the parameters $\delta$ and $\epsilon$, in contrast 
to the case shown in Figure~\ref{fig06AD}, the shapes of the convex polygons in the tile sets 
in Cases 1--4 are not uniquely determined. Tilings using the shape of 
\mbox{Tile$(1, 1)$} in the tile sets of Cases 1 and 4 were edge-to-edge tilings with convex 
polygons. Tilings using the shape of \mbox{Tile$(1, 1)$} in the tile sets of Cases 2 and 3 
were non-edge-to-edge tilings with convex polygons.

Note that, because the ASPs comprising several types of convex polygons shown in the study 
of this appendix are based on \mbox{Tile$(1, 1)$} shown in \cite{Smith_2024a} and assumes that 
\mbox{Tile$(1, 1)$} is the chiral aperiodic monotile shown in \cite{Smith_2024b}, it is not confirmed 
that these sets are actually ASPs. Thus, because \cite{Smith_2024b} is still a preview version, 
it does not hold if \mbox{Tile$(1, 1)$} is not a chiral aperiodic monotile.

%%%%%%%%%%%%%%%%%%%%%%%%%%%%%%%%%%%%%%%%%%%%%%%%%%%%%%%%%%%%%%%%%%%%%%

\subsection{Conditions for treating anterior and posterior sides as the same type}
\label{appD_4}

The ASP created by $S_t (\delta, \epsilon)$ that satisfies the restriction on the 
number of types of congruent diagrams can be regarded as ``an ASP comprising 
three types of convex polygons" from the viewpoint that does not 
distinguish between anterior and posterior sides.

If $S_t (\delta, \epsilon)$ is subject to the same rules as the monohedral tiling, that is, 
``the use of reflected tiles is allowed in the tiling," then $S_t (\delta, \epsilon)$ 
is not ASP. \mbox{Tile$(1, 1)$} can generate a periodic tiling if the use of reflected 
tiles is allowed during the tiling generation process \cite{Smith_2024a}. 
In other words, if the reflected $CP(X_{14})$ is allowed in tiling, the set comprises of 
six types of convex polygons that can generate periodic tiling.

On the contrary, $S_t (\delta, \epsilon)$ considered the following ``tile 
combination specification condition" as a method to be able to argue that 
``ASP comprises of three types of convex polygons" from the viewpoint that 
does not distinguish between anterior and posterior sides.

\begin{nameth2}
A set of convex polygons must contain two or more convex polygons that 
cannot be reflected during the tiling generation process. Convex polygons 
that cannot be reflected during the tiling generation process do not exhibit line 
symmetry \rm{(}\it{i.e., they can be distinguished between the anterior and 
posterior sides}\rm{)}.
\end{nameth2}

If ``must contain two or more convex polygons that cannot be reflected during the 
tiling generation process" of the above condition is satisfied and the set 
contains one type of convex polygon, then it is a convex polygonal monotile, 
which does not require the use of reflected tiles for tiling, belongs to the known Type 
families, and can generate periodic tiling. Therefore, we added ``contains 
two or more convex polygons" to the condition.

Moreover, ``convex polygons that cannot be reflected during the tiling generation 
process are assumed not to have the property of line symmetry" is important. 
As mentioned previously, the set of six convex polygons with $S_t (\delta, \epsilon)$ 
and reflected $CP(X_{14})$ can form periodic tiling; however, it may also be 
possible to form periodic tiling in other cases. For example, for the case of 
$S_t (\delta = 105^ \circ , \epsilon = 75^ \circ )$ in Figure~\ref{fig06AD}, 
if $CP(X_{14})$, $CP(X_{3})$, and reflected $CP(X_{14})$ are used, 
periodic tiling is possible, as shown in Figure~\ref{fig09AD}. In this case, the merged 
diagram of $CP(X_{14})$ and reflected $CP(X_{14})$ can be regarded as forming 
a convex octagon. In other words, Figure~\ref{fig09AD} can be regarded as periodic 
tiling using anterior and posterior convex octagons and convex pentagons 
without any distinction between the anterior and posterior sides (two types 
of convex pentagons and octagons that are not monotiles from a viewpoint 
that does not distinguish between the anterior and posterior sides).

\renewcommand{\figurename}{{\small Figure}}
\begin{figure}[htbp]
\centering\includegraphics[width=14.5cm,clip]{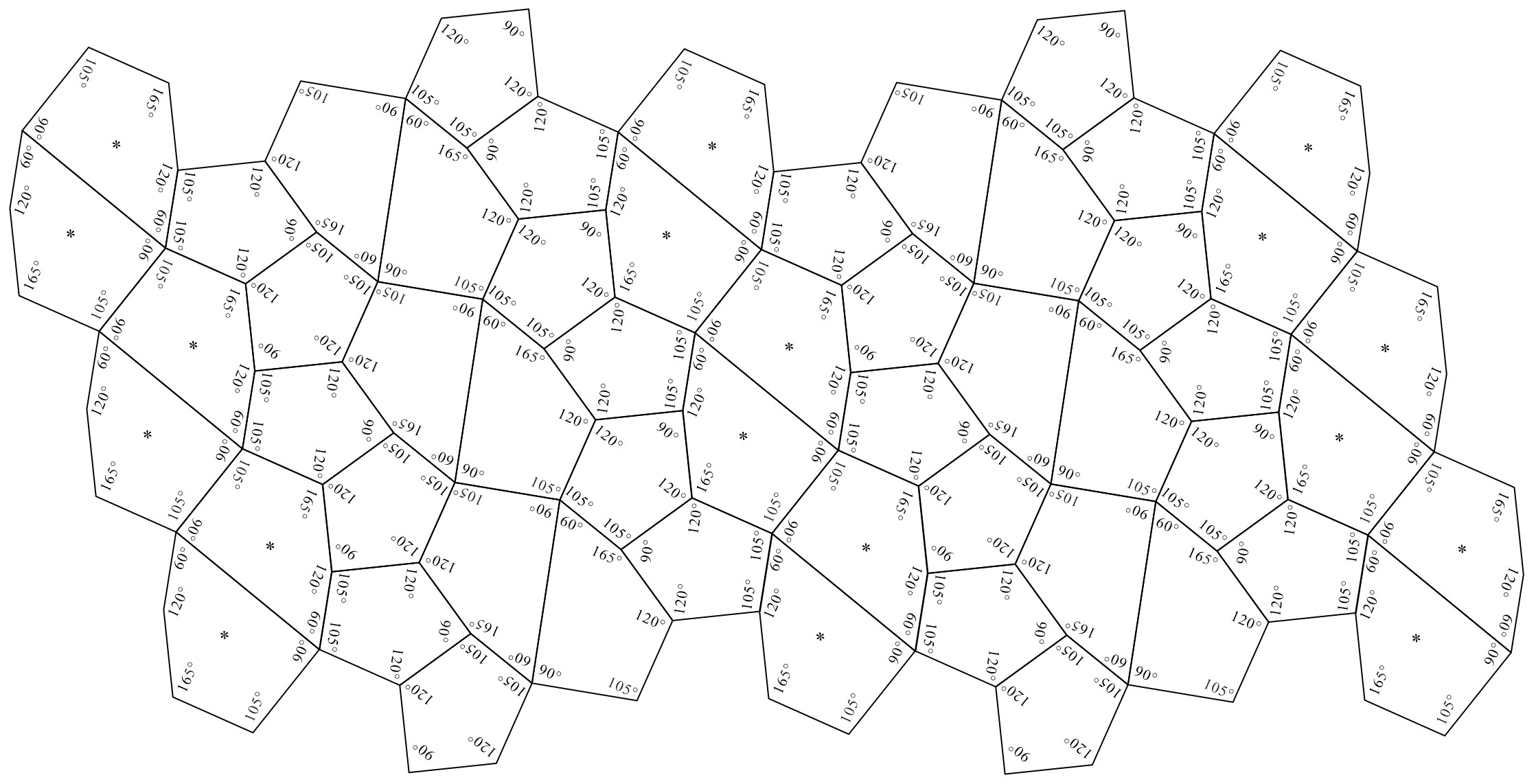} 
\caption{{\small 
Periodic tiling with $CP(X_{14})$, $CP(X_{3})$, and reflected 
$CP(X_{14})$ as shown in Figure~\ref{fig06AD}.
The reflected $CP(X_{14})$ are given an asterisk mark ``*.''
}
\label{fig09AD}
}
\end{figure}

%%%%%%%%%%%%%%%%%%%%%%%%%%%%%%%%%%%%%%%%%%%%%%%%%%%%%%%%%%%%%%%%%%%%%%

\subsection{Supplement 1}
\label{appD_5}

In addition to Figure~\ref{fig03AD}, there is another method for dividing \mbox{Tile$(1, 1)$} 
such that the number of convex polygonal types is less than or equal to 
three from a viewpoint that does not distinguish between the anterior and 
posterior sides, as shown in Figure~\ref{fig10AD}. However, the sets of convex 
polygons created by the divisions in Figure~\ref{fig10AD} always contain convex 
pentagonal monotiles. The convex pentagonal monotile belonging to both the 
Type 2 and Type 4 families with line symmetry can form periodic tiling (see \ref{appA}).
For the division that makes the convex pentagonal monotiles belonging to the Type 13 
family in Figure~\ref{fig10AD}, the tile set of the four types of convex polygons 
contains the anterior and posterior sides of convex pentagonal monotiles belonging 
to the Type 13 family. (For the convex pentagonal monotiles belonging to the Type 13 
family in Figure~\ref{fig10AD}, the tile set of the convex polygons created by the 
division contains the reflected tile. Thus, the convex pentagons can form the 
representative tiling of the Type 13 family.) Therefore, the tile sets of convex 
polygons shown in Figure~\ref{fig10AD} are not ASP.

\renewcommand{\figurename}{{\small Figure}}
\begin{figure}[htbp]
\centering\includegraphics[width=11cm,clip]{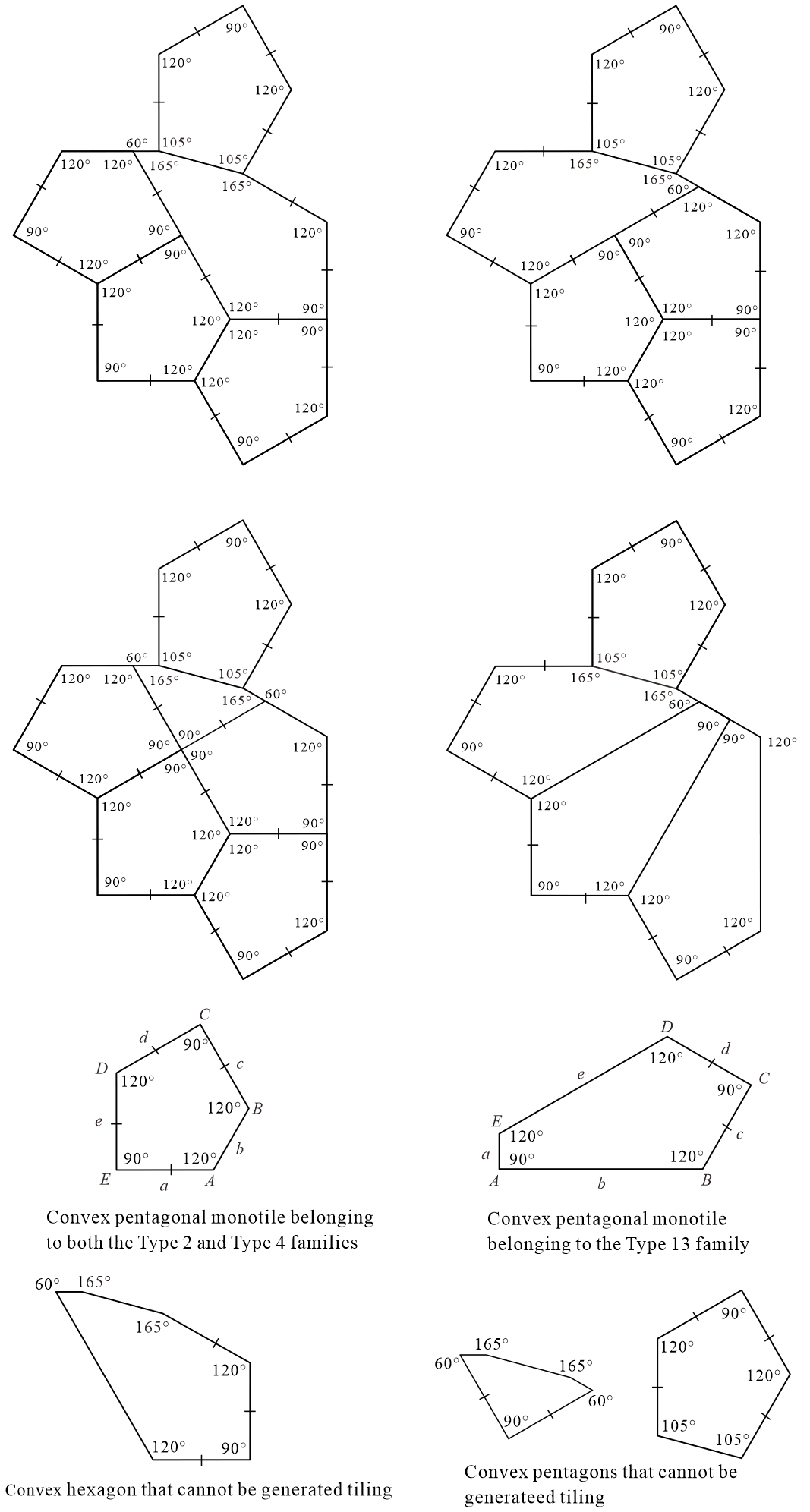} 
\caption{{\small 
How to divide \mbox{Tile$(1, 1)$} such that the number of convex polygonal 
types is less than or equal to three from the viewpoint that does not 
distinguish between the anterior and posterior sides.
}
\label{fig10AD}
}
\end{figure}

Using the division method that makes the convex pentagonal monotiles belong 
to the Type 13 family in Figure~\ref{fig10AD}, we created a set with three types of 
convex pentagons and one type of convex hexagon, as shown in Figure~\ref{fig11AD}(a). 
If the set in Figure~\ref{fig11AD}(a) is an ASP, it is an $ASPcvx(4,4)$ and does not 
satisfy the restriction on the number of types of congruent diagrams.

If the restriction on the number of types of congruent diagrams does not 
have to be satisfied, there are various methods for dividing \mbox{Tile$(1, 1)$} into 
four, five, or more types of convex polygons that are likely to be ASP. 
Examples of these divisions are shown in Figure~\ref{fig12AD}. There is no guarantee 
that the sets of convex polygons in these divides are ASP. However, it can 
be inferred that there are various methods for dividing a \mbox{Tile$(1, 1)$} into 
four or more convex polygons, such that it is an ASP.

\renewcommand{\figurename}{{\small Figure}}
\begin{figure}[htbp]
\centering\includegraphics[width=10cm,clip]{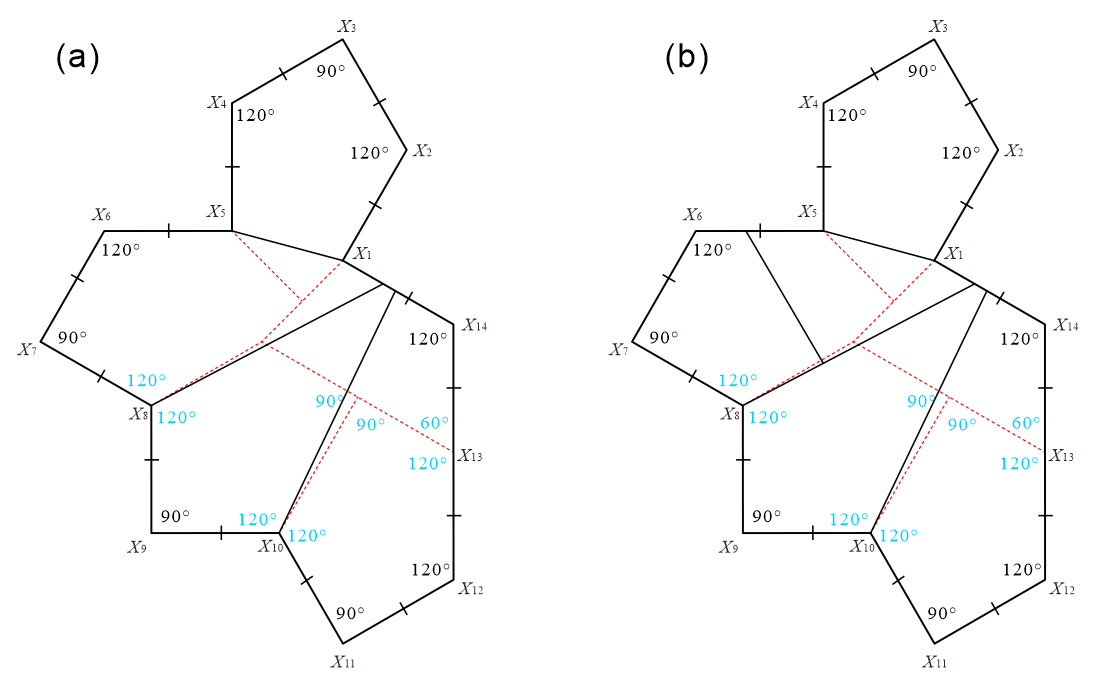} 
\caption{{\small 
Example of a method for dividing \mbox{Tile$(1, 1)$} that does not satisfy the 
restriction on the number of types of congruent diagrams, even if $ASPcvx(4)$ 
or $ASPcvx(5)$.
}
\label{fig11AD}
}
\end{figure}

\renewcommand{\figurename}{{\small Figure}}
\begin{figure}[htbp]
\centering\includegraphics[width=10cm,clip]{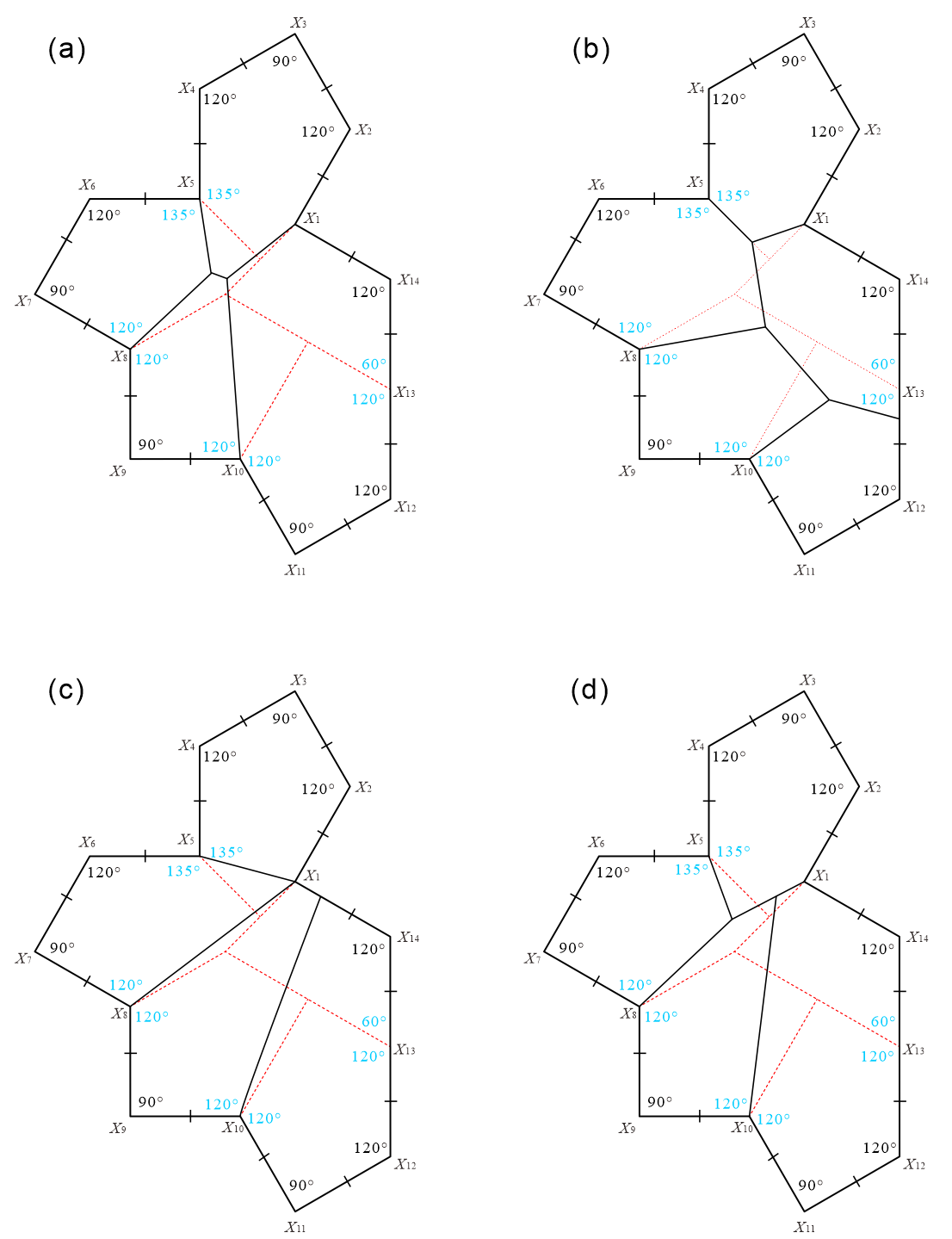} 
\caption{{\small 
Example of dividing \mbox{Tile$(1, 1)$} into four or five convex polygons that could be ASP.
}
\label{fig12AD}
}
\end{figure}

A convex polygon with six or more edges can be successfully divided into two 
or more convex polygons with five or more edges. Figure~\ref{fig11AD}(b) shows an 
example of dividing the convex hexagon in Figure~\ref{fig11AD}(a) into two convex 
pentagons. Another example is the division of the interior of a convex 
pentagon into six convex pentagons. Using this division, it is possible to 
increase the number of convex polygon types in the ASP. Conversely, merging 
several convex polygons into a single convex polygon may be possible. 
However, changing the type of convex polygon in a set by dividing or merging 
is not essential for the search and discussion of ASP. Therefore, in this study, 
we only discuss whether it is possible to construct an ASP comprising 
fewer types of convex polygons.

%%%%%%%%%%%%%%%%%%%%%%%%%%%%%%%%%%%%%%%%%%%%%%%%%%%%%%%%%%%%%%%%%%%%%%

\subsection{Supplement 2}
\label{appD_6}

As shown in Figure~\ref{fig13AD}, we found that there is another method for dividing 
\mbox{Tile$(1, 1)$} such that the number of convex polygonal types is less than or 
equal to three from a viewpoint that does not distinguish between the anterior 
and posterior sides. In Figure~\ref{fig13AD}, the point $Q$ moves on the edge $X_{5}X_{6}$ 
and the point $T$ moves on the edge $X_{12}X_{13}$ of \mbox{Tile$(1, 1)$}. 
In Figure~\ref{fig13AD}, the convex pentagon $CP_{2}(X_{3})$ with 
the points ``$X_{3}$, $X_{4}$, $X_{5}$, $X_{1}$, $X_{2}$" as vertices is fixed shape. 
Additionally, the convex pentagon $CP_{2}(X_{7})$ with the points 
``$X_{7}$, $X_{8}$, $R$, $Q$, $X_{6}$" as vertices, the convex pentagon 
$CP_{2}(X_{9})$ with the points ``$X_{9}$, $X_{10}$, $S$, $R$, $X_{8}$" as 
vertices, and the convex pentagon $CP_{2}(X_{11})$ with the points 
``$X_{11}$, $X_{12}$, $T$, $S$, $X_{10}$" as vertices are congruent convex 
pentagons. However, $CP_{2}(X_{7})$ and $CP_{2}(X_{11})$ have a reflection 
relation with $CP_{2}(X_{9})$. The convex polygon $CP_{2}(X_{14})$ with the points 
``$X_{14}$, $X_{1}$, $X_{5}$, $Q$, $R$, $S$, $T$" as vertices is a heptagon, 
but it is a hexagon if $\zeta = 120^ \circ $ and $\eta = 90^ \circ $.

If the five convex polygons in Figure~\ref{fig13AD} are not convex pentagonal 
or hexagonal monotiles, we consider that a set of four types of convex 
polygons is a candidate of $ASPcvx(4, 3)$.

\begin{itemize}
\item Candidates for $ASPcvx(4, 3)$
\begin{itemize}
\item[$\circ$] \textbf{Case 5}. When $\zeta = 120^ \circ $ and $\eta = 90^ \circ $ 
in Figures~\ref{fig13AD}, there is a ``set of three types of convex pentagons ($CP_{2}(X_{3})$, 
$CP_{2}(X_{7})$, $CP_{2}(X_{9})$, $CP_{2}(X_{11})$) and one type of convex hexagon 
($CP_{2}(X_{14})$)," and, from the viewpoint that does not distinguish between anterior and 
posterior sides, it is a ``set of two types of convex pentagons and one type of convex 
hexagon."
\end{itemize}

\begin{itemize}
\item[$\circ$] \textbf{Case 6}. When $\zeta \ne 120^ \circ $ and $\eta \ne 90^ \circ $ 
in Figures~\ref{fig13AD}, there is a ``set of three types of convex pentagons ($CP_{2}(X_{3})$, 
$CP_{2}(X_{7})$, $CP_{2}(X_{9})$, $CP_{2}(X_{11})$) and one type of convex heptagon 
($CP_{2}(X_{14})$)," and, from the viewpoint that does not distinguish between anterior and 
posterior sides, it is a ``set of two types of convex pentagons and one type of convex 
heptagon."
\end{itemize}
\end{itemize}

\renewcommand{\figurename}{{\small Figure}}
\begin{figure}[htbp]
\centering\includegraphics[width=11.5cm,clip]{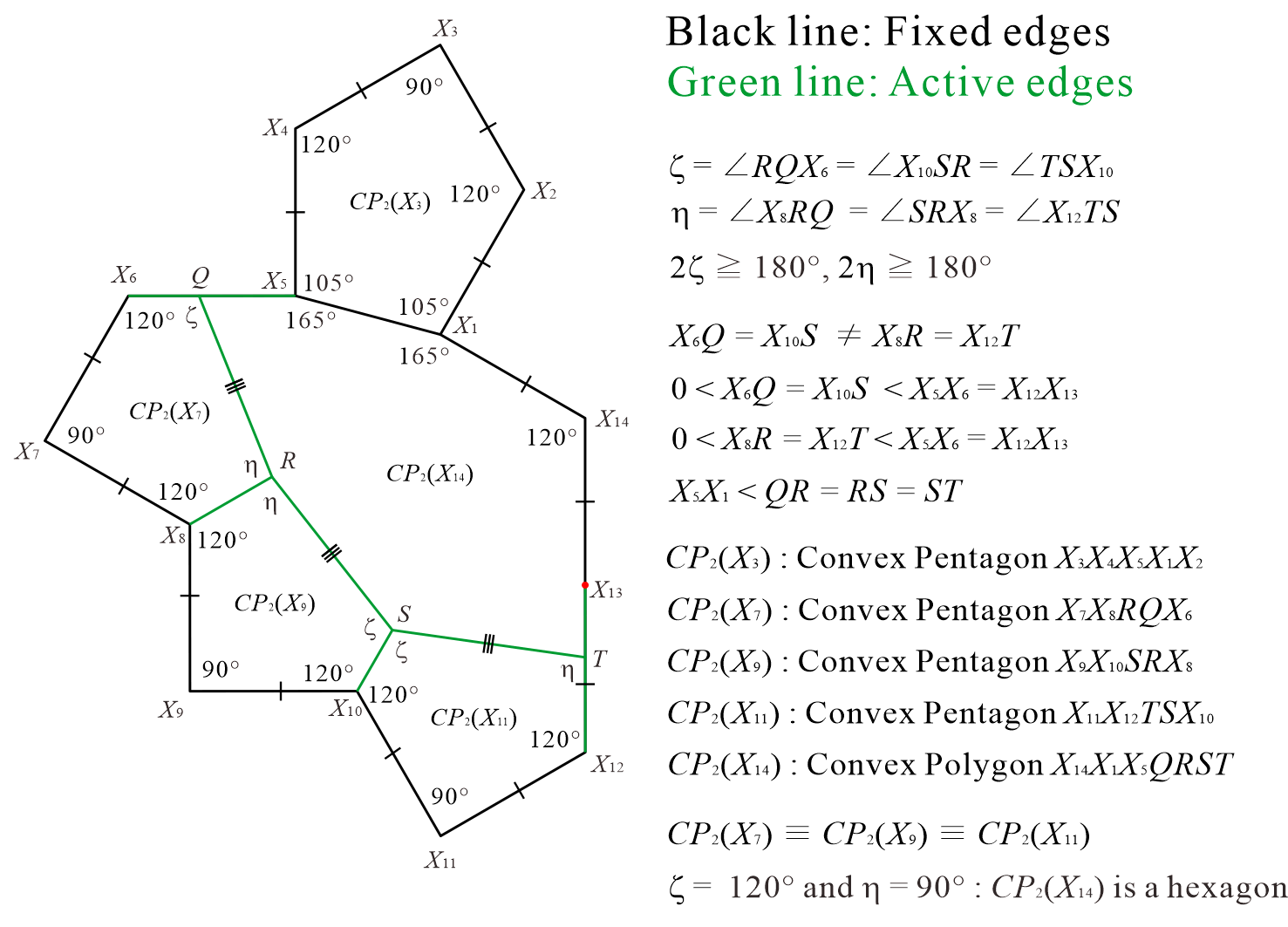} 
\caption{{\small 
How to divide \mbox{Tile$(1, 1)$} containing convex polygons that 
are candidates for $ASPcvx(4, 3)$.
}
\label{fig13AD}
}
\end{figure}

%%%%%%%%%%%%%%%%%%%%%%%%%%%%%%%%%%%%%%%%%%%%%%%%%%%%%%%%%%%%%%%%%%%%%%
%%%%%%%%%%%%%%%%%%%%%%%%%%%%%%%%%%%%%%%%%%%%%%%%%%%%%%%%%%%%%%%%%%%%%%
\end{document}